\PassOptionsToPackage{noend}{algpseudocode}
\PassOptionsToPackage{disable}{todonotes}
\PassOptionsToPackage{compress}{cleveref}
\PassOptionsToPackage{appendix,abbr}{clevethm}
\PassOptionsToPackage{usenames,dvipsnames,table}{xcolor}
\documentclass[10pt, reqno, twoside]{amsart}
\usepackage[foot]{amsaddr}
\PassOptionsToPackage{caption=false}{subfig}
\PassOptionsToPackage{textsize=tiny,textwidth=1.8cm}{todonotes}
\PassOptionsToPackage{noend}{algpseudocode}
\usepackage[%
	eqreset=section,
	thmreset=section,
]{myPreamble}
\usepackage{mathrsfs}
\usepackage{marginnote}
\usepackage{wrapfig}
\usepgfplotslibrary{fillbetween}

\newcommand\cont{{\rm C}}
\newcommand\epicomp[2]{(#1#2)}

	\newcommand\FBE{\varphi_\gamma^{\text{\sc fb}}}

	\newcommand\sDRS{\varphi_1}
	\newcommand\nsDRS{\varphi_2}
	\newcommand\DRE{\varphi_\gamma^{\text{\sc dr}}}
	\renewcommand\FBsmooth{\sDRS}
	\renewcommand\FBnonsmooth{\nsDRS}

	\newcommand\LL{\mathscr L}

	\newcommand\note[1]{%
		\unskip
		\bgroup
		\let\marginpar\marginnote
		\reversemarginpar
		\todo[
			backgroundcolor=MidnightBlue!10,
			linecolor=MidnightBlue!25,
			bordercolor=MidnightBlue!25,
		]{#1}
		\egroup
	}

\showgrayfalse

	\title[DRS and ADMM for nonconvex optimization]{%
		Douglas-Rachford splitting and ADMM for nonconvex optimization: tight convergence results%
	}
	
	\author[A. Themelis]{Andreas Themelis}
	\author[P. Patrinos]{Panagiotis Patrinos}
	\address{%
		\normalfont
		\TheAddressKU.
		\textit{E-mail:}
		\href{mailto:andreas.themelis@esat.kuleuven.be}{andreas.themelis@esat.kuleuven.be},
		\href{mailto:panos.patrinos@esat.kuleuven.be}{panos.patrinos@esat.kuleuven.be}.%
		\newline
		This work was supported by:
		KU Leuven internal funding: StG/15/043 Fonds de la Recherche Scientifique -- FNRS and the Fonds Wetenschappelijk Onderzoek -- Vlaanderen under EOS Project no 30468160 (SeLMA) FWO projects: G086318N; G086518N.%
	}

\begin{document}

	\begin{abstract}
		Although originally designed and analyzed for convex problems, the alternating direction method of multipliers (ADMM) and its close relatives, Douglas-Rachford splitting (DRS) and Peaceman-Rachford splitting (PRS), have been observed to perform remarkably well when applied to certain classes of structured nonconvex optimization problems.
		However, partial global convergence results in the nonconvex setting have only recently emerged.
		In this paper we show how the Douglas-Rachford envelope (DRE), introduced in 2014, can be employed to unify and considerably simplify the theory for devising global convergence guarantees for ADMM, DRS and PRS applied to nonconvex problems under less restrictive conditions, larger prox-stepsizes and over-relaxation parameters than previously known.
		In fact, our bounds are tight whenever the over-relaxation parameter ranges in \((0,2]\).
		The analysis of ADMM uses a universal primal equivalence with DRS that generalizes the known duality of the algorithms.
	\end{abstract}
	\maketitle

%


	\section{Introduction}
		First introduced in \cite{douglas1956numerical} for finding numerical solutions of heat differential equations, the \DEF{Douglas-Rachford splitting} (DRS) is now considered a textbook algorithm in convex optimization or, more generally, in monotone inclusion problems.
As the name suggests, DRS is a \emph{splitting scheme}, meaning that it works on a problem decomposition by addressing each component separately, rather than operating on the whole problem which is typically too hard to be tackled directly.
In optimization, the objective to be minimized is \emph{split} as the sum of two functions, resulting in the following canonical framework addressed by DRS:
\begin{equation}\label{eq:P}
	\minimize_{s\in\R^p} \varphi(s)\equiv\sDRS(s)+\nsDRS(s).
\end{equation}
Here, \(\func{\sDRS,\nsDRS}{\R^p}{\Rinf}\) are proper, lower semicontinuous (lsc), extended-real-valued functions (\(\Rinf\coloneqq\R\cup\set\infty\) denotes the extended-real line).
Starting from some \(s\in\R^p\), one DR-iteration applied to \eqref{eq:P} with \emph{stepsize} \(\gamma>0\) and \emph{relaxation} parameter \(\lambda>0\) amounts to
\[
	\tag{DRS}\label{DRS}
	\begin{cases}[rl]
		u\phantom{^+}
	{}\in{} &
		\prox_{\gamma\sDRS}(s)
	\\
		v\phantom{^+}
	{}\in{} &
		\prox_{\gamma\nsDRS}(2u-s)
	\\
		s^+
	{}={} &
		s+\lambda(v-u).
	\end{cases}
\]
The case \(\lambda=1\) corresponds to the classical DRS, whereas for \(\lambda=2\) the scheme is also known as Peaceman-Rachford splitting (PRS).
If \(s\) is a \emph{fixed point} for the DR-iteration --- that is, such that \(s^+=s\)
--- then it can be easily seen that \(u\) satisfies the first-order necessary condition for optimality in problem \eqref{eq:P}.
When both \(\sDRS\) and \(\nsDRS\) are convex functions, the condition is also sufficient
and \ref{DRS} iterations are known to converge for any \(\gamma>0\) and \(\lambda\in(0,2)\).

Closely related to DRS and possibly even more popular is the \DEF{alternating direction method of multipliers} (ADMM), first appeared in \cite{glowinski1975approximation,gabay1976dual}, see also \cite{glowinski2014alternating} for a recent historical overview.
ADMM addresses linearly constrained optimization problems
\begin{equation}\label{eq:CP}
	\minimize_{(x,z)\in\R^m\times\R^n}{
		f(x)+g(z)
	}
\quad
	\stt Ax+Bz=b,
\end{equation}
where \(\func{f}{\R^m}{\Rinf}\), \(\func{g}{\R^n}{\Rinf}\), \(A\in\R^{p\times m}\), \(B\in\R^{p\times n}\), and \(b\in\R^p\).
ADMM is an iterative scheme based on the following recursive steps
\[\tag{ADMM}\label{ADMM}
	\renewcommand\arraystretch{1.2}
	\begin{cases}[l >{{}}c<{{}} l]
		y^{\nicefrac+2}
	&=&
		y-\beta(1-\lambda)(Ax+Bz-b)
	\\
		x^+
	&\in&
		\argmin\LL_\beta({}\cdot{},z,y^{\nicefrac+2})
	\\
		y^+
	&=&
		y^{\nicefrac+2} + \beta(Ax^++Bz-b)
	\\
		z^+
	&\in&
		\argmin\LL_\beta(x^+,{}\cdot{},y^+).
	\end{cases}
\]
Here, \(\beta>0\) is a \DEF{penalty} parameter, \(\lambda>0\) is a possible \DEF{relaxation} parameter, and
\begin{align}\label{eq:L}
	\LL_\beta(x,z,y)
{}\coloneqq{} &
	f(x)+g(z)+\innprod{y}{Ax+Bz-b}+\tfrac\beta2\|Ax+Bz-b\|^2
\end{align}
is the \(\beta\)-augmented Lagrangian of \eqref{eq:CP} with \(y\in\R^p\) as Lagrange equality multiplier.
It is well known that for convex problems \ref{ADMM} is simply \ref{DRS} applied to a dual formulation \cite{gabay1983chapter}, and its convergence properties for \(\lambda=1\) and arbitrary penalty parameters \(\beta>0\) are well documented in the literature, see \eg \cite{boyd2011distributed}.
Recently, \ref{DRS} and \ref{ADMM} have been observed to perform remarkably well when applied to certain classes of structured nonconvex optimization problems and partial or case-specific convergence results have also emerged.

		\subsection{Contributions}
			Our contributions can be summarized as follows.
\begin{enumerate}[leftmargin=\widthof{\hspace*{\labelsep}3)}]
\item\emph{New tight convergence results for nonconvex \ref{DRS}}.
	We provide novel convergence results for \ref{DRS} applied to nonconvex problems with one function being Lipschitz-dif\-fer\-en\-tiable (\Cref{thm:DRS}).
	Differently from the results in the literature, we make no a priori assumption on the existence of accumulation points and we consider all relaxation parameters \(\lambda\in(0,4)\), as opposed to \(\lambda\in\set{1,2}\).
	Moreover, our results are tight for all \(\lambda\in(0,2]\) (\Cref{thm:DRStight}).
	\Cref{fig:DRS,fig:PRS} highlight the extent of the improvement with respect to the state of the art.
	\item\emph{Primal equivalence of \ref{DRS} and \ref{ADMM}}.
		We prove the equivalence of \ref{DRS} and \ref{ADMM} for arbitrary problems and relaxation parameters, so extending the well-known duality of the algorithms holding in the convex case and the recently observed primal equivalence when \(\lambda=1\).
\item\emph{New convergence results for \ref{ADMM}}.
	Thanks to the equivalence with \ref{DRS}, not only do we provide new convergence results for the \ref{ADMM} scheme, but we also offer an elegant unifying framework that greatly simplifies and generalizes the theory in the literature, is based on less restrictive assumptions, and provides explicit bounds for stepsizes and possible other coefficients.
	A comparison with the state of the art is shown in \Cref{fig:ADMM}.
\item\emph{A continuous and exact merit function for \ref{DRS} and \ref{ADMM}}.
	Our results are based on the \emph{Douglas-Rachford Envelope} (DRE), first introduced in \cite{patrinos2014douglas} for convex problems and here generalized.
	The DRE extends the known properties of the Moreau envelope and its connections to the proximal point algorithm, to composite functions as in \eqref{eq:P} and \eqref{eq:CP}.
	In particular, we show that the DRE serves as an exact, continuous and real-valued (as opposed to extended-real-valued) merit function for the original problem, computable with quantities obtained in the iterations of \ref{DRS} (or \ref{ADMM}).
\end{enumerate}
Finally, we propose out-of-the-box implementations of \ref{DRS} and \ref{ADMM} where the stepsize \(\gamma\) and the penalty parameter \(\beta\) are adaptively tuned, so that no prior knowledge of quantities such as Lipschitz moduli is needed.

		\subsection{Comparisons \& related work}
			We now compare our results with a selection of recent related works which, to the best of our knowledge, represent the state of the art for generality and contributions.

\subsubsection{ADMM}
	A primal equivalence of \ref{DRS} and \ref{ADMM} has been observed in \cite[Rem. 3.14]{bauschke2015projection} when \(A=-B=\I\) and \(\lambda=1\).
	In \cite[Thm. 1]{yan2016self} the equivalence is extended to arbitrary matrices; although limited to convex problems, the result is easily extendable.
	Our generalization to any relaxation parameter (and nonconvex problems) is largely based on this result and uses the same problem reformulation proposed therein.
	The relaxation considered in this paper corresponds to that introduced in \cite{eckstein1992douglas}; it is worth mentioning that another type of relaxation has been proposed, corresponding to \(\lambda=1\) in \eqref{ADMM} but with a different steplength for the \(y\)-update: that is, with \(\beta\) replaced by \(\theta\beta\) for some \(\theta>0\).
	The known convergence results for \(\theta\in(0,\frac{1+\sqrt 5}{2})\) in the convex case, see \cite[\S5]{glowinski2013numerical}, were recently extended to nonconvex problems and for \(\theta\in(0,2)\) in \cite{goncalves2017convergence}.
	
	In \cite{wang2018global} convergence of ADMM is studied for problems of the form
	\[
		\minimize_{\bm x=(x_0\ldots x_p),z}{
			g(\bm x)
			{}+{}
			{\textstyle\sum_{i=0}^pf_i(x_i)}
			{}+{}
			h(z)
		}
	\quad\stt~
		\bm A\bm x+Bz=0.
	\]
	Although addressing a more general class of problem than \eqref{eq:CP}, when specialized to the standard two-function formulation analyzed in this paper it relies on numerous assumptions.
	These include Lipschitz continuous minimizers of all ADMM subproblems (in particular, uniqueness of their solution).
	For instance, the requirements rule out interesting cases involving discrete variables or rank constraints.

	In \cite{hong2016convergence} a class of nonconvex problems with more than two functions is presented and variants of ADMM with deterministic and random updates are discussed.
	The paper provides a nice theory and explicit bounds for the penalty paramenter in ADMM, which agree with ours in best- and worst-case scenerarios, but are more restrictive otherwise (cf. \cref{fig:ADMM} for a more detailed comparison).
	The main limitation of the proposed approach, however, is that the theory only allows for functions either convex or smooth, differently from ours where the nonsmooth term can basically be anything.
	Once again, many interesting applications are not covered.

	The work \cite{li2015global} studies a proximal ADMM where a possible Bregman divergence term in the second block update is considered.
	By discarding the Bregman term so as to recover the original ADMM scheme, the same bound on the stepsize as in \cite{hong2016convergence} is found.
	Another proximal variant is proposed in \cite{goncalves2017convergence}, under less restrictive assumptions related to the concept of smoothness relative to a matrix that we will introduce in \Cref{defin:Bsmooth}.
	When matrix \(B\) has full-column rank, the proximal term can be discarded and their method reduces to the classical ADMM.
	
	The problem addressed in \cite{guo2017convergence} is fully covered by our analysis, as they consider \ref{ADMM} for \eqref{eq:CP} where \(f\) is \(L\)-Lipschitz continuously differentiable and \(B\) is the identity matrix.
	Their bound \(\beta>2L\) for the penalty parameter is more conservative than ours; in fact, the two coincide only in a worst-case scenario.
	
\subsubsection{Douglas-Rachford splitting}
	Few exceptions apart \cite{li2016douglas,li2017peaceman}, advances in nonconvex \ref{DRS} theory are problem specific and only provide local convergence results, at best.
	These mainly focus on feasibility problems, where the goal is to find points in the intersection of nonempty closed sets \(A\) and \(B\) subjected to some regularity conditions.
	This is done by applying \ref{DRS} to the minimization of the sum of \(\sDRS=\indicator_A\) and \(\nsDRS=\indicator_B\), where \(\indicator_C\) is the \emph{indicator function} of a set \(C\) (see \cref{sec:Notation}).
	The minimization subproblems in \ref{DRS} then reduce to (set-valued) projections onto either set, regardless of the stepsize parameter \(\gamma>0\).
	This is the case of \cite{bauschke2014local}, for instance, where \(A\) and \(B\) are finite unions of convex sets.
	Local linear convergence when \(A\) is affine, under some conditions on the (nonconvex) set \(B\), are shown in \cite{hesse2013nonconvex,hesse2014alternating}.

	Although this particular application of \ref{DRS} does not comply with our requirements, as \(\sDRS\) fails to be Lipschitz differentiable, however replacing \(\indicator_A\) with \(\sDRS=\tfrac12\dist_A^2\) yields an equivalent problem which fits into our framework when \(A\) is a convex set.
	In terms of \ref{DRS} iterations, this simply amounts to replacing \(\proj_A\), the projection onto set \(A\), with a ``relaxed'' version \(\proj_{A,t}\coloneqq(1-t)\id+t\proj_A\) for some \(t\in(0,1)\).
	Then, it can be easily verified that for any \(\alpha,\beta\in(0,+\infty]\) one \ref{DRS}-step applied to
	\begin{equation}\label{eq:P-MARP}
	\smash{
		\minimize_{s\in\R^p}{
			\tfrac\alpha2\dist_A^2(s)
			{}+{}
			\tfrac\beta2\dist_B^2(s)
		}
	}
	\end{equation}
	results in
	\begin{equation}\label{eq:MARP}
	\smash{
		s^+
	{}\in{}
		(1-\nicefrac\lambda2)s
		{}+{}
		\nicefrac\lambda2
		\proj_{B,q}\proj_{A,p}s
	}
	\end{equation}
	for
	\(
		p
	{}={}
		\tfrac{2\alpha\gamma}{1+\alpha\gamma}
	\)
	and
	\(
		q
	{}={}
		\tfrac{2\beta\gamma}{1+\beta\gamma}
	\).
	Notice that \eqref{eq:MARP} is the \(\nicefrac\lambda2\)-relaxation of the ``method of alternating \((p,q)\)-relaxed projections'' (\((p,q)\)-MARP) \cite{bauschke2014method}.
	The (non-relaxed) \((p,q)\)-MARP is recovered by setting \(\lambda=2\), that is, by applying PRS to \eqref{eq:P-MARP}.
	Local linear convergence of MARP was shown when \(A\) and \(B\), both possibly nonconvex, satisfy some constraint qualifications, and also global convergence when some other requirements are met.
	When set \(A\) is convex, then \(\tfrac\alpha2\dist_A^2\) is convex and \(\alpha\)-Lipschitz differentiable; our theory then ensures convergence of the \emph{fixed-point residual} and subsequential convergence of the iterations \eqref{eq:MARP} for any \(\lambda\in(0,2)\), \(p\in(0,1)\) and \(q\in(0,1]\), without any requirements on the (nonempty closed) set \(B\).
	Here, \(q=1\) is obtained by replacing \(\tfrac\beta2\dist_B^2\) with \(\indicator_B\), which can be interpreted as the hard penalization obtained by letting \(\beta=\infty\).
	Although the non-relaxed MARP is not covered due to the non-strong convexity of \(\dist_A^2\), however \(\lambda\) can be set arbitrarily close to \(2\).

	\begin{figure}[tbh]%
		\centering
		\begin{subfigure}{0.325\linewidth}
			{{%
			\pgfkeys{/pgf/images/include external/.code={\includegraphics[width=\linewidth]{#width=\linewidth}}}%
			\tikzsetnextfilename{DRS}%
			\input{./TeX/Tikz/DRS.tex}%
		}}%
			\caption{}%
			\label{fig:DRS}%
		\end{subfigure}
		\begin{subfigure}{0.325\linewidth}
			{{%
			\pgfkeys{/pgf/images/include external/.code={\includegraphics[width=\linewidth]{#width=\linewidth}}}%
			\tikzsetnextfilename{PRS}%
			\input{./TeX/Tikz/PRS.tex}%
		}}%
			\caption{}%
			\label{fig:PRS}%
		\end{subfigure}
		\begin{subfigure}{0.325\linewidth}
			{{%
			\pgfkeys{/pgf/images/include external/.code={\includegraphics[width=\linewidth]{#width=\linewidth}}}%
			\tikzsetnextfilename{ADMM}%
			\input{./TeX/Tikz/ADMM.tex}%
		}}%
			\caption{}%
			\label{fig:ADMM}%
		\end{subfigure}
		\caption[%
			Comparison of stepsize ranges ensuring convergence of DRS, PRS and ADMM between our bounds and the best in the literature.%
		]{%
			Maximum stepsize \(\gamma\) ensuring convergence of DRS (\cref{fig:DRS}) and PRS (\cref{fig:PRS}), and maximum inverse of the penalty paramenter \(\nicefrac1\beta\) in ADMM (\cref{fig:ADMM}); comparison between our bounds (blue plot) and \cite{li2016douglas} for DRS, \cite{li2017peaceman} for PRS and \cite{goncalves2017convergence,guo2017convergence,hong2016convergence,li2015global,wang2018global} for ADMM.
			On the \(x\)-axis the ratio between hypoconvexity parameter \(\sigma\) and the Lipschitz modulus \(L\) of the gradient of the smooth function.
			On the \(y\)-axis, the supremum of stepsize \(\gamma\) such that the algorithms converge.
			For ADMM, the analysis is made for a common framework: 2-block ADMM with no Bregman or proximal terms, Lipschitz-differentiable \(f\), \(A\) invertible and \(B\) identity; \(L\) and \(\sigma\) are relative to the transformed problem.
			Notice that, due to the proved analogy of DRS and ADMM, our theoretical bounds coincide in \cref{fig:DRS,fig:ADMM}.%
		}%
		\label{fig:Comparisons}%
	\end{figure}

	The work \cite{li2016douglas} presents the first general analysis of global convergence of DRS (non-relaxed) for fully nonconvex problems where one function is Lipschitz differentiable.
	In \cite{li2017peaceman} PRS is also considered under the additional requirement that the smooth function is strongly convex with strong-convexity/Lipschitz moduli ratio of at least \(\nicefrac23\).
	For sufficiently small (explicitly computable) stepsizes one iteration of DRS or PRS yields a sufficient decrease on an augmented Lagrangian, and the generated sequences remain bounded when the cost function has bounded level sets.

	Other than completing the analysis to all relaxation parameters \(\lambda\in(0,4)\), as opposed to \(\lambda\in\set{1,2}\), we improve their results by showing convergence for a considerably larger range of stepsizes and, in the case of PRS, with no restriction on the strong convexity modulus of the smooth function.
	We also show that our bounds are optimal whenever \(\lambda\in(0,2]\).
	The extent of the improvement is evident in the comparisons outlined in \Cref{fig:Comparisons}.
	Thanks to the lower boundedness of the DRE, as opposed to the lower unbounded augmented Lagrangian, we show that the vanishing of the fixed-point residual occurs without coercivity assumptions.

		\subsection{Organization of the paper}
			The paper is organized as follows.
\Cref{sec:Background} introduces some notation and offers a brief recap of the needed theory.
In \Cref{sec:DRE}, after formally stating the needed assumptions for the \ref{DRS} problem formulation \eqref{eq:P} we introduce the DRE and analyze in detail its key properties.
Based on these properties, in \Cref{sec:DRS} we prove convergence results of \ref{DRS} and show the tightness of our findings by means of suitable counterexamples.
In \Cref{sec:ADMM} we deal with \ref{ADMM} and show its equivalence with \ref{DRS}; based on this, convergence results for \ref{ADMM} are derived from the ones already proven for \ref{DRS}.
\Cref{sec:Conclusion} concludes the paper.
For the sake of readability, some proofs and auxiliary results are deferred to the Appendix.

	\section{Background}
		\label{sec:Background}%

		\subsection{Notation}
			\label{sec:Notation}%
			The extended-real line is \(\Rinf=\R\cup\set{\infty}\).
The positive and negative parts of \(r\in\R\) are defined respectively as
\(
	[r]_+
{}\coloneqq{}
	\max\set{0,r}
\)
and
\(
	[r]_-
{}\coloneqq{}
	\max\set{0,-r}
\),
so that \(r=[r]_+-[r]_-\).
We adopt the convention that \(\nicefrac10=\infty\).

The open and closed balls centered in \(x\) and with radius \(r\) are denoted by \(\ball xr\) and \(\cball xr\), respectively.
With \(\id\) we indicate the identity function \(x\mapsto x\) defined on a suitable space, and with \(\I\) the identity matrix of suitable size.
For a nonzero matrix \(M\in\R^{p\times n}\) we let \(\sigma_+(M)\) denote its smallest nonzero singular value.

For a set \(E\) and a sequence \(\seq{x^k}\) we write \(\seq{x^k}\subset E\) to indicate that \(x^k\in E\) for all \(k\in\N\).
We say that \(\seq{x^k}\subset\R^n\) is \DEF{summable} if \(\sum_{k\in\N}\|x^k\|\) is finite, and \DEF{square-summable} if \(\seq{\|x^k\|^2}\) is summable.


The \DEF{domain} of an extended-real-valued function \(\func{h}{\R^n}{\Rinf}\) is the set
\(
	\dom h
{}\coloneqq{}
	\set{x\in\R^n}[
		h(x)<\infty
	]
\),
while its \DEF{epigraph} is the set
\(
	\epi h
{}\coloneqq{}
	\set{(x,\alpha)\in\R^n\times\R}[
		h(x)\leq\alpha
	]
\).
\(h\) is said to be \DEF{proper} if \(\dom h\neq\emptyset\), and \DEF{lower semicontinuous (lsc)} if \(\epi h\) is a closed subset of \(\R^{n+1}\).
For \(\alpha\in\R\), \(\lev_{\leq\alpha}h\) is the \DEF{\(\alpha\)-level set} of \(h\), \ie
\(
	\lev_{\leq\alpha}h
{}\coloneqq{}
	\set{x\in\R^n}[
		h(x)\leq\alpha
	]
\).
We say that \(h\) is \emph{level bounded} if \(\lev_{\leq\alpha}h\) is bounded for all \(\alpha\in\R\).

 
 

We use the notation \(\ffunc H{\R^n}{\R^m}\) to indicate a point-to-set mapping \(\func H{\R^n}{\mathcal P(\R^m)}\), where \(\mathcal P(\R^m)\) is the power set of \(\R^m\) (the set of all subsets of \(\R^m\)).
The \DEF{graph} of \(H\) is the set
\(
	\graph H
{}\coloneqq{}
	\set{(x,y)\in\R^n\times\R^m}[
		y\in H(x)
	]
\).
 

We denote by \(\ffunc{\hat\partial h}{\R^n}{\R^n}\) the \DEF{regular subdifferential} of \(h\), where
\begin{equation}\label{eq:hatpartial}
	v\in\hat\partial h(\bar x)
\quad\Leftrightarrow\quad
	\liminf_{\limsubstack{x&\to&\bar x\\x&\neq&\bar x}}{
		\frac{h(x)-h(\bar x)-\innprod{v}{x-\bar x}}{\|x-\bar x\|}
	}
{}\geq{}
	0.
\end{equation}
A necessary condition for local minimality of \(x\) for \(h\) is \(0\in\hat\partial h(x)\), see \cite[Thm. 10.1]{rockafellar2011variational}.
The (limiting) \DEF{subdifferential} of \(h\) is \(\ffunc{\partial h}{\R^n}{\R^n}\), where
\(
	v\in\partial h(\bar x)
\)
iff there exists a sequence \(\seq{x^k,v^k}\subseteq\graph\hat\partial h\) such that
\[
	\lim_{k\to\infty}(x^k,h(x^k),v^k)
{}={}
	(x,h(x),v).
\]
The set of \DEF{horizon subgradients} of \(h\) at \(x\) is \(\partial^\infty h(x)\), defined as \(\partial h(x)\) except that \(v^k\to v\) is meant in the ``cosmic'' sense, namely \(\lambda_kv^k\to v\) for some \(\lambda_k\searrow 0\).

		\subsection{Smoothness and hypoconvexity}
			\label{sec:smooth}%
			The class of functions \(\func{h}{\R^n}{\R}\) that are \(k\) times continuously differentiable is denoted as \(\cont^k(\R^n)\).
We write \(h\in\cont^{1,1}(\R^n)\) to indicate that \(h\in\cont^1(\R^n)\) and that \(\nabla h\) is Lipschitz continuous with modulus \(L_h\).
To simplify the terminology, we will say that such an \(h\) is \DEF{\(L_h\)-smooth}.
It follows from \cite[Prop. A.24]{bertsekas2016nonlinear} that if \(h\) is \(L_h\)-smooth, then
\(
	|h(y)-h(x)-\innprod{\nabla h(x)}{y-x}|
{}\leq{}
	\tfrac{L_h}{2}\|y-x\|^2
\)
for all \(x,y\in\R^n\).
In particular, there exists \(\sigma_h\in[-L_h,L_h]\) such that \(h\) is \DEF{\(\sigma_h\)-hypoconvex}, in the sense that \(h-\tfrac{\sigma_h}{2}\|{}\cdot{}\|^2\) is a convex function.
Thus, every \(L_h\)-smooth and \(\sigma_h\)-hypoconvex function \(h\) satisfies
\begin{equation}\label{eq:LipBound}
	\tfrac{\sigma_h}{2}
	\|y-x\|^2
{}\leq{}
	h(y)-h(x)-\innprod{\nabla h(x)}{y-x}
{}\leq{}
	\tfrac{L_h}{2}
	\|y-x\|^2
\quad
	\forall x,y\in\R^n.
\end{equation}
By applying \cite[Thm. 2.1.5]{nesterov2003introductory} to the (convex) function \(\psi=h-\tfrac\sigma2\|{}\cdot{}\|^2\) we obtain that this is equivalent to having
\begin{equation}\label{eq:innprod}
	\sigma_h
	\|y-x\|^2
{}\leq{}
	\innprod{\nabla h(y)-\nabla h(x)}{y-x}
{}\leq{}
	L_h
	\|y-x\|^2
\quad
	\forall x,y\in\R^n.
\end{equation}
Note that \(\sigma\)-hypoconvexity generalizes the notion of (strong) convexity by allowing negative strong convexity moduli.
In fact, if \(\sigma=0\) then \(\sigma\)-hypoconvexity reduces to convexity, while for \(\sigma>0\) it denotes \(\sigma\)-strong convexity.
\begin{lem}[Subdifferential characterization of smoothness]\label{thm:C11subdiff}%
	Let \(\func h{\R^n}{\R}\) be such that \(\partial h(x)\neq\emptyset\) for all \(x\in\R^n\), and suppose that there exist \(L\geq0\) and \(\sigma\in[-L,L]\) such that
	\begin{equation}\label{eq:ineqC11}
		\sigma
		\|x_1-x_2\|^2
	{}\leq{}
		\innprod{v_1-v_2}{x_1-x_2}
	{}\leq{}
		L
		\|x_1-x_2\|^2
	\quad
		\forall x_i\in\R^n,~
		v_i\in\partial h(x_i),~
		i=1,2.
	\end{equation}
	Then, \(h\in\cont^{1,1}(\R^n)\) is \(L\)-smooth and \(\sigma\)-hypoconvex.
	\begin{proof}
		See \Cref{proof:thm:C11subdiff}.
	\end{proof}
\end{lem}

\begin{thm}[Lower bounds for smooth functions]\label{thm:rho}%
	Let \(h\in C^{1,1}(\R^n)\) be \(L_h\)-smooth and \(\sigma_h\)-hypoconvex.
	Then, for all \(x,y\in\R^n\) it holds that
	\[
		h(y)
	{}\geq{}
		h(x)
		{}+{}
		\innprod{\nabla h(x)}{y-x}
		{}+{}
		\rho(y,x),
	\]
	where
	\begin{enumerate}
	\item\label{thm:hypo:rho0}%
		either~
		\(
			\rho(y,x)
		{}={}
			\tfrac{\sigma_h}{2}
			\|y-x\|^2
		\),
	\item\label{thm:hypo:rho1}%
		or~
		\(
			\rho(y,x)
		{}={}
			\tfrac{\sigma_hL_h}{2(L_h+\sigma_h)}
			\|y-x\|^2
			{}+{}
			\tfrac{1}{2(L_h+\sigma_h)}
			\|\nabla h(y)-\nabla h(x)\|^2
		\),~
		provided that~
		\(
			-L_h<\sigma_h\leq 0
		\).
	\end{enumerate}
	Clearly, all inequalities remain valid if one replaces \(L_h\) with any \(L\geq L_h\) and \(\sigma_h\) with any \(\sigma\in[-L,\sigma_h]\).
	\begin{proof}
		See \Cref{proof:thm:rho}.
	\end{proof}
\end{thm}

		\subsection{Proximal mapping}
			\label{sec:prox}%
			The \DEF{proximal mapping} of \(\func{h}{\R^n}{\Rinf}\) with parameter \(\gamma>0\) is \(\ffunc{\prox_{\gamma h}}{\R^n}{\dom h}\) defined as
\begin{equation}
	\prox_{\gamma h}(x)
{}\coloneqq{}
	\argmin_{w\in\R^n}\set{
		h(w) + \tfrac{1}{2\gamma}\|w-x\|^2
	}.
\end{equation}
We say that a function \(h\) is \DEF{prox-bounded} if \(h+\tfrac{1}{2\gamma}\|{}\cdot{}\|^2\) is lower bounded for some \(\gamma>0\).
The supremum of all such \(\gamma\)
is the \emph{threshold of prox-boundedness of \(h\)}, denoted as \(\gamma_h\).
If \(h\) is lsc, then \(\prox_{\gamma h}\) is nonempty- and compact-valued over \(\R^n\) for all \(\gamma\in(0,\gamma_h)\) \cite[Thm. 1.25]{rockafellar2011variational}.
Consequently, the value function of the minimization problem defining the proximal mapping, namely the \DEF{Moreau envelope} with stepsize \(\gamma\in(0,\gamma_h)\), denoted by \(\func{h^\gamma}{\R^n}{\R}\) and defined as
\begin{equation}
	\fillwidthof[r]{
		\prox_{\gamma h}(x)
	}{
		h^\gamma(x)
	}
{}\coloneqq{}
	\fillwidthof[l]{
		\argmin_{w\in\R^n}\set{
			h(w) + \tfrac{1}{2\gamma}\|w-x\|^2
		}.
	}{
		\inf_{w\in\R^n}\set{
			h(w) + \tfrac{1}{2\gamma}\|w-x\|^2
		},
	}
\end{equation}
is everywhere finite and, in fact, strictly continuous \cite[Ex. 10.32]{rockafellar2011variational}.
Moreover, the necessary optimality conditions of the problem defining \(\prox_{\gamma g}\) together with \cite[Thm. 10.1 and Ex. 8.8]{rockafellar2011variational} imply that
\begin{equation}\label{eq:prox:subdiff}
	\tfrac1\gamma(x-\bar x)
{}\in{}
	\hat\partial h(\bar x)
\quad
	\forall\bar x\in\prox_{\gamma h}(x).
\end{equation}
When \(h\in\cont^{1,1}(\R^n)\), its proximal mapping and Moreau envelope enjoy many favorable properties which we summarize next.
\begin{prop}[Proximal properties of smooth functions]\label{thm:proxf}%
	Let \(h\in \cont^{1,1}(\dom h)\) be \(L_h\)-smooth, hence \(\sigma_h\)-hypoconvex for some \(\sigma_h\in[-L_h,L_h]\).
	Then, \(h\) is prox-bounded with \(\gamma_{h}\geq\nicefrac{1}{[\sigma_{h}]_-}\) and for all \(\gamma<\nicefrac{1}{[\sigma_{h}]_-}\) the following hold:
	\begin{enumerate}
	\item\label{thm:proxfEquiv}
		\(\prox_{\gamma h}\) is single valued, and for all \(s\in\R^n\) it holds that \(u=\prox_{\gamma h}(s)\) iff
		\(
			s=u+\gamma\nabla h(u)
		\).
	\item\label{thm:proxfLip}
		\(\prox_{\gamma h}\) is \((\tfrac{1}{1+\gamma L_{h}})\)-strongly monotone and \((1+\gamma\sigma_h)\)-cocoercive, in the sense that
		\[
			\innprod{u-u'}{s-s'}
		{}\geq{}
			\tfrac{1}{1+\gamma L_h}
			\|s-s'\|^2
		\quad\text{and}\quad
			\innprod{u-u'}{s-s'}
		{}\geq{}
			(1+\gamma\sigma_h)
			\|u-u'\|^2
		\]
		for all \(s,s'\in\R^n\), where \(u=\prox_{\gamma h}(s)\) and \(u'=\prox_{\gamma h}(s')\).
		In particular,
		\begin{equation}\label{eq:biLip}
			\tfrac{1}{1+\gamma L_h}
			\|s-s'\|
		{}\leq{}
			\|u-u'\|
		{}\leq{}
			\tfrac{1}{1+\gamma\sigma_h}
			\|s-s'\|.
		\end{equation}
		Thus, \(\prox_{\gamma h}\) is a \(\tfrac{1}{1+\gamma\sigma_h}\)-Lipschitz and invertible mapping, and its inverse \(\id+\gamma\nabla h\) is \((1+\gamma L_h)\)-Lipschitz continuous.
	\item\label{thm:MoreaufC1}
		\(h^\gamma\in \cont^{1,1}(\R^n)\) is \(L_{h^\gamma}\)-smooth and \(\sigma_{h^\gamma}\)-hypoconvex, with
		\(
			L_{h^\gamma}
		{}={}
			\max\set{
				\tfrac{L_{h}}{1+\gamma L_{h}},~
				\tfrac{[\sigma_{h}]_-}{1+\gamma\sigma_{h}}
			}
		\)
		and
		\(
			\sigma_{h^\gamma}
		{}={}
			\tfrac{\sigma_{h}}{1+\gamma\sigma_{h}}
		\).
		Moreover,
		\(
			\nabla h^\gamma(s)
		{}={}
			\tfrac{1}{\gamma}(s-\prox_{\gamma h}(s))
		\)
		and
		\(
			\nabla h(\prox_{\gamma h}(s))
		{}={}
			\tfrac{1}{\gamma}\bigl(s-\prox_{\gamma h}(s)\bigr)
		\).
	\end{enumerate}
	\begin{proof}
		See \Cref{proof:thm:proxf}.
	\end{proof}
\end{prop}

	\section{Douglas-Rachford envelope}
		\label{sec:DRE}%
		We now list the blanket assumptions for the functions in problem \eqref{eq:P}.
\begin{ass}[Requirements for the DRS formulation \eqref{eq:P}]\label{ass:DRS}%
	The following hold
	\begin{enumerate}
	\item\label{ass:DRS:f}%
		\(\sDRS\in\cont^{1,1}(\R^n)\) is \(L_{\sDRS}\)-smooth, hence \(\sigma_{\sDRS}\)-hypoconvex for some \(\sigma_{\sDRS}\in[-L_{\sDRS},L_{\sDRS}]\).
	\item\label{ass:DRS:g}%
		\(\nsDRS\) is proper and lsc.
	\item\label{ass:DRS:LB}%
		Problem \eqref{eq:P} has a solution, that is, \(\argmin\varphi\neq\emptyset\).
	\end{enumerate}
\end{ass}
\begin{rem}[Feasible stepsizes for \ref{DRS}]\label{thm:proxbounded}%
	Under \Cref{ass:DRS}, both \(\sDRS\) and \(\nsDRS\) are prox-bounded with threshold at least \(\nicefrac{1}{L_{\sDRS}}\), and in particular \ref{DRS} iterations are well defined for all \(\gamma\in(0,\nicefrac{1}{L_{\sDRS}})\).
	That \(\gamma_{\sDRS}\geq\nicefrac{1}{L{\sDRS}}\) follows from \Cref{thm:proxf}, having \(\nicefrac{1}{[\sigma_{\sDRS}]_-}\geq\nicefrac{1}{L_{\sDRS}}\).
	As for \(\nsDRS\), for all \(s\in\R^p\) it holds that
	\[
		\inf\varphi
	{}\leq{}
		\sDRS(s)+\nsDRS(s)
	{}\overrel*[\leq]{\eqref{eq:LipBound}}{}
		\sDRS(0)+\innprod{\nabla\sDRS(0)}{s}+\tfrac{L_{\sDRS}}{2}\|s\|^2
		{}+{}
		\nsDRS(s),
	\]
	hence, for all \(\gamma<\nicefrac{1}{L_{\sDRS}}\) the function \(s\mapsto\nsDRS(s)+\tfrac{1}{2\gamma}\|s\|^2\) is lower bounded.
\end{rem}
Starting from \(s\in\R^p\), let us consider variables \((u,v)\) generated by a \ref{DRS} step under \Cref{ass:DRS}.
As first noted in \cite{patrinos2014douglas}, from the relation \(s=u+\gamma\nabla\sDRS(u)\) (see \cref{thm:proxfEquiv}) it follows that
\begin{equation}\label{eq:vFBS}
	v\in\FB u
\end{equation}
is the result of a \DEF{forward-backward step} at \(u\), amounting to
\begin{align}\label{eq:DRSv}
	v
{}\in{}
	\argmin_{w\in\R^p}{} &
	\set{
		\nsDRS(w)
		{}+{}
		\smashunderbracket{
			\sDRS(u)
			{}+{}
			\innprod{\nabla\sDRS(u)}{w-u}
			{}+{}
			\tfrac{1}{2\gamma}\|w-u\|^2
		}{}
	},
\shortintertext{%
	see \eg \cite{bolte2014proximal,themelis2018forward} for an extensive discussion on nonconvex forward-backward splitting (FBS).
	This shows that \(v\) is the result of the minimization of a \emph{majorization model} for the original function \(\varphi=\sDRS+\nsDRS\), where the smooth function \(\sDRS\) is replaced by the quadratic upper bound emphasized by the under-bracket in \eqref{eq:DRSv}.
	First introduced in \cite{patrinos2014douglas} for convex problems, the \DEF{Douglas-Rachford envelope} (DRE) is the function \(\func{\DRE}{\R^n}{\R}\) defined as%
}
\label{eq:DRE}
	\DRE(s)
{}\coloneqq{}
	\min_{w\in\R^p}{} &
	\set{
		\nsDRS(w)
		{}+{}
		\sDRS(u)
		{}+{}
		\innprod{\nabla\sDRS(u)}{w-u}
		{}+{}
		\tfrac{1}{2\gamma}\|w-u\|^2
	}.
\end{align}
Namely, rather than the \emph{minimizer} \(v\), \(\DRE(s)\) is the \emph{value} of the minimization problem \eqref{eq:DRSv} defining the \(v\)-update in \eqref{DRS}.
The expression \eqref{eq:DRE} emphasizes the close connection that the DRE has with the \DEF{forward-backward envelope} (FBE) as in \cite{themelis2018forward}, here denoted \(\FBE\), namely
\begin{equation}\label{eq:FBEu}
	\DRE(s)
{}={}
	\FBE(u),
\quad
	\text{where \(u=\prox_{\gamma\sDRS}(s)\).}
\end{equation}
The FBE is an exact penalty function for FBS, which was initially proposed for convex problems in \cite{patrinos2013proximal} and later extended and further analyzed in \cite{stella2017forward,themelis2018forward,liu2017further}.
In this section we will see that, under \Cref{ass:DRS}, the DRE serves a similar role with respect to \ref{DRS} which will be key for establishing (tight) convergence results in the nonconvex setting.
Another useful intepretation of the DRE is obtained by plugging the minimizer \(w=v\) in \eqref{eq:DRE}.
This leads to
\begin{equation}\label{eq:DRE=L}
	\DRE(s)
{}={}
	\LL_{\nicefrac1\gamma}(u,v,\gamma^{-1}(u-s)),
\end{equation}
where \(u\) and \(v\) come from the \ref{DRS} iteration and
\begin{equation}\label{eq:DRS-L}
	\LL_\beta(x,z,y)
{}\coloneqq{}
	\sDRS(x)+\nsDRS(z)+\innprod{y}{x-z}+\tfrac\beta2\|x-z\|^2
\end{equation}
is the \(\beta\)-augmented Lagrangian relative to the equivalent problem formulation
\begin{equation}\label{eq:DRS-constrained}
	\minimize_{x,z\in\R^p}\sDRS(x)+\nsDRS(z)
\quad
	\stt x-z=0.
\end{equation}
This expression also emphasizes that evaluating \(\DRE(s)\) requires the same operations as performing one \ref{DRS} update
\(
	s\mapsto(u,v)
\).

		\subsection{Properties}
			\label{sec:Properties}%
			Building upon the connection with the FBE emphasized in \eqref{eq:FBEu}, in this section we highlight some important properties enjoyed by the DRE.
We start by observing that \(\DRE\) is a strictly continuous function for \(\gamma<\nicefrac{1}{L_{\sDRS}}\), owing to the fact that so is the FBE \cite[Prop. 4.2]{themelis2018forward}, and that \(\prox_{\gamma\sDRS}\) is Lipschitz continuous as shown in \cref{thm:proxfLip}.
\begin{prop}[Strict continuity]\label{thm:C0}%
	Suppose that \Cref{ass:DRS} is satisfied.
	For all \(\gamma<\nicefrac{1}{L_{\sDRS}}\) the DRE \(\DRE\) is a real-valued and strictly continuous function.
\end{prop}
Next, we investigate on the fundamental connections relating the DRE \(\DRE\) and the cost function \(\varphi\).
\note{loc minim equiv needed? Maybe for linear cvg?}%
We show, for  \(\gamma\) small enough and up to an (invertible) change of variable, that infima and minimizers of the two functions coincide, as well as equivalence of level boundedness of \(\varphi\) and \(\DRE\).
\begin{prop}[Sandwiching property]\label{thm:sandwich}%
	Suppose that \Cref{ass:DRS} is satisfied.
	Let \(\gamma<\nicefrac{1}{L_{\sDRS}}\) be fixed, and consider \(u,v\) generated by one \ref{DRS} iteration starting from \(s\in\R^p\).
	Then,
	\begin{enumerate}
	\item\label{thm:DREleq}
		\(\DRE(s)\leq\varphi(u)\).
	\item\label{thm:DREgeq}
		\(
			\varphi(v)
		{}\leq{}
			\DRE(s)
			{}-{}
			\tfrac{1-\gamma L_{\sDRS}}{2\gamma}\|u-v\|^2
		\).
	\end{enumerate}
	\begin{proof}
		\ref{thm:DREleq} is easily inferred from definition \eqref{eq:DRE} by considering \(w=u\).
		Moreover, it follows from \cite[Prop. 4.3]{themelis2018forward} and the fact that \(v\in\FB u\), cf. \eqref{eq:vFBS}, that
		\(
			\varphi(v)
		{}\leq{}
			\FBE(u)
			{}-{}
			\tfrac{1-\gamma L_{\sDRS}}{2\gamma}\|u-v\|^2
		\).
		\ref{thm:DREgeq} then follows from \eqref{eq:FBEu}.
	\end{proof}
\end{prop}
\begin{thm}[Minimization and level-boundedness equivalence]\label{thm:DREequiv}%
	Suppose that \Cref{ass:DRS} is satisfied.
	For any \(\gamma<\nicefrac{1}{L_{\sDRS}}\) the following hold:
	\begin{enumerate}
	\item\label{thm:inf}%
		\(\inf\varphi=\inf\DRE\).
	\item\label{thm:argmin}%
		\(\argmin\varphi=\prox_{\gamma\sDRS}\bigl(\argmin\DRE\bigr)\).
	\item\label{thm:boundedlev}%
		\(\varphi\) is level bounded iff so is \(\DRE\).
	\end{enumerate}
	\begin{proof}
		It follows from \cite[Thm. 4.4]{themelis2018forward} that the FBE satisfies \(\inf\varphi=\inf\FBE\) and \(\argmin\varphi=\argmin\FBE\).
		The similar properties \ref{thm:inf} and \ref{thm:argmin} of the DRE then follow from the identity \(\DRE=\FBE\circ\prox_{\gamma\sDRS}\), cf. \eqref{eq:FBEu}, and the fact that \(\prox_{\gamma\sDRS}\) is invertible, as shown in \cref{thm:proxf}.
		
		We now show \ref{thm:boundedlev}.
		Denote
		\(
			\varphi_\star
		{}\coloneqq{}
			\inf\varphi
		{}={}
			\inf\DRE
		\),
		which is finite by assumption.
		\begin{proofitemize}
		\item
			Suppose that \(\DRE\) is level bounded, and let \(u\in\lev_{\leq\alpha}\varphi\) for some \(\alpha>\varphi_\star\).
			Then, \(s\coloneqq u+\gamma\nabla\sDRS(u)\) is such that \(\prox_{\gamma\sDRS}(s)=u\), as shown in \cref{thm:proxfEquiv}.
			Thus, from \cref{thm:sandwich} it follows that \(s\in\lev_{\leq\alpha}\DRE\).
			In particular,
			\(
				\lev_{\leq\alpha}\varphi
			{}\subseteq{}
				[I+\gamma\nabla\sDRS]
				(\lev_{\leq\alpha}\DRE),
			\)
			and since \(\I+\gamma\nabla\sDRS\) is Lipschitz continuous and \(\lev_{\leq\alpha}\DRE\) is bounded by assumption, it follows that \(\lev_{\leq\alpha}\varphi\) is also bounded.
		\item
			Suppose now that \(\DRE\) is not level bounded.
			Then, there exists \(\alpha>\varphi_\star\) together with a sequence \(\seq{s_k}\) satisfying
			\(
				s_k
			{}\in{}
				\lev_{\leq\alpha}\DRE
				{}\setminus{}
				\ball 0k
			\)
			for all \(k\in\N\).
			Let \(u_k\coloneqq\prox_{\gamma\sDRS}(s_k)\), so that \(s_k=u_k+\gamma\nabla\sDRS(u_k)\) (cf. \cref{thm:proxfEquiv}), and let \(v_k\in\FB{u_k}\).
			From \cref{thm:DREgeq} it then follows that \(v_k\in\lev_{\leq\alpha}\varphi\), and that
			\[
				\alpha
				{}-{}
				\varphi_\star
			{}\geq{}
				\DRE(s_k)
				{}-{}
				\varphi_\star
			{}\geq{}
				\DRE(s_k)
				{}-{}
				\varphi(v_k)
			{}\geq{}
				\tfrac{1-\gamma L_{\sDRS}}{2\gamma}
				\|u_k-v_k\|^2.
			\]
			Therefore,
			\(
				\|u_k-v_k\|^2
			{}\leq{}
				\tfrac{
					2\gamma(\alpha-\varphi_\star)
				}{
					1-\gamma L_{\sDRS}
				}
			\)
			and
			\begin{align*}
				\|v_k\|
			{}\geq{} &
				\|u_k-u_0\|
				{}-{}
				\|u_0\|
				{}-{}
				\|u_k-v_k\|
			\smash{{}\overrel[\geq]{\ref{thm:proxfLip}}{}}
				\tfrac{1}{1+\gamma L_{\sDRS}}
				\|s_k-s_0\|
				{}-{}
				\|u_0\|
				{}-{}
				\|u_k-v_k\|
			\\
			{}\geq{} &
				\tfrac{
					k-\|s_0\|
				}{
					1+\gamma L_{\sDRS}
				}
				{}-{}
				\|u_0\|
				{}-{}
				\sqrt{
					\tfrac{
						2\gamma(\alpha-\varphi_\star)
					}{
						1-\gamma L_{\sDRS}
					}
				}
			{}\to{}
				+\infty
			\quad\text{as \(k\to\infty\).}
			\end{align*}
			This shows that \(\lev_{\leq\alpha}\varphi\) is also unbounded.
			\qedhere
		\end{proofitemize}
	\end{proof}
\end{thm}

	\section{Convergence of Douglas-Rachford splitting}
		\label{sec:DRS}%
		Closely related to the DRE, the augmented Lagrangian \eqref{eq:DRS-L} was used in \cite{li2016douglas} under the name of \emph{Douglas-Rachford merit function} to analyze \ref{DRS} for the special case \(\lambda=1\).
It was shown that for sufficiently small \(\gamma\) there exists \(c>0\) such that the iterates generated by \ref{DRS} satisfy
\begin{equation}\label{eq:LPong}
	\LL_{\nicefrac1\gamma}(u^{k+1},v^{k+1},\eta^{k+1})
{}\leq{}
	\LL_{\nicefrac1\gamma}(u^k,v^k,\eta^k)
	{}-{}
	c\|u^k-u^{k+1}\|^2
\quad\text{with \(\eta^k=\gamma^{-1}(v^k-s^k)\),}
\end{equation}
to infer that \(\seq{u^k}\) and \(\seq{v^k}\) have same accumulation points, all of which are stationary for \(\varphi\).
In \cite{li2017peaceman}, where also the case \(\lambda=2\) is addressed with a slightly different penalty function, it was then shown that the sequence remains bounded and thus accumulation points exist in case \(\varphi\) is level bounded.
We now generalize the decrease property \eqref{eq:LPong} shown in \cite{li2016douglas,li2017peaceman} by considering arbitrary relaxation parameters \(\lambda\in(0,4)\) (as opposed to \(\lambda\in\set{1,2}\)) and providing tight ranges for the stepsize \(\gamma\) whenever \(\lambda\in(0,2]\).
Thanks to the lower boundedness of \(\DRE\), it will be possible to show that the \ref{DRS} residual vanishes without any coercivity assumption.
\grayout{%
	We are only interested in the case \(\gamma<\nicefrac{1}{L_{\sDRS}}\), for otherwise the DRE may fail to be lower bounded, as \Cref{es:DREunbounded} demonstrates.
	Morever, it will be shown in \Cref{sec:DRS} that the bound \(\gamma<\nicefrac{1}{L_{\sDRS}}\) is necessary for ensuring the convergence of \ref{DRS}, unless the generality of \Cref{ass:fg} is sacrificed.
}%

\begin{thm}[Sufficient decrease on the DRE]\label{thm:DRE:SD}%
	Suppose that \Cref{ass:DRS} is satisfied, and consider one \ref{DRS} update
	\(
		s
	{}\mapsto{}
		(u,v,s^+)
	\)
	for some stepsize
	\(
		\gamma
	{}<{}
		\min\set{
			\tfrac{2-\lambda}{2[\sigma_{\sDRS}]_-},\,
			\tfrac{1}{L_{\sDRS}}
		}
	\)
	and relaxation \(\lambda\in(0,2)\).
	Then,
	\begin{equation}\label{eq:DRE:SD}
		\DRE(s)-\DRE(s^+)
	{}\geq{}
		\tfrac{c}{(1+\gamma L_{\sDRS})^2}\,
		\|s-s^+\|^2,
	\end{equation}
	where, denoting \(p_{\sDRS}\coloneqq\nicefrac{\sigma_{\sDRS}}{L_{\sDRS}}\in[-1,1]\), \(c\) is a strictly positive constant defined as%
	\footnote{%
		A one-line expression for the constant is
		\(
			c
		{}={}
			\tfrac{2-\lambda}{2\lambda\gamma}
			{}-{}
			\min\set{
				\tfrac{[p_{\sDRS}]_-}{\lambda},~
				L_{\sDRS}\max\set{
					\tfrac{[\sigma_{\sDRS}]_-}{2(1-[p_{\sDRS}]_-)},\,
					\tfrac12-\tfrac{\gamma L_{\sDRS}}{\lambda}
				}
			}
		\).
	}
	\begin{equation}\label{eq:DRE:c0}
		c
	{}={}
		\frac{2-\lambda}{2\lambda\gamma}
		{}-{}
		\begin{ifcases}
			L_{\sDRS}\max\set{
				\tfrac{[p_{\sDRS}]_-}{2(1-[p_{\sDRS}]_-)},\,
				\tfrac12-\tfrac{\gamma L_{\sDRS}}{\lambda}
			}
		&
				p_{\sDRS}
			{}\geq{}
				\tfrac\lambda2-1
		\\
			\tfrac{[\sigma_{\sDRS}]_-}{\lambda}
		\otherwise[otherwise.]
		\end{ifcases}
	\end{equation}
	If \(\sDRS\) is strongly convex, then \eqref{eq:DRE:SD} also holds for
	\begin{equation}\label{eq:lambda>=2}
		2
	{}\leq{}
		\lambda
	{}<{}
		\tfrac{4}{1+\sqrt{1-p_{\sDRS}}}
	\quad\text{and}\quad
		\tfrac{
			p_{\sDRS}\lambda
			{}-{}
			\delta
		}{
			4\sigma_{\sDRS}
		}
	{}<{}
		\gamma
	{}<{}
		\tfrac{
			p_{\sDRS}\lambda
			{}+{}
			\delta
		}{
			4\sigma_{\sDRS}
		},
	\end{equation}
	where
	\(
		\delta
	{}\coloneqq{}
		\sqrt{(p_{\sDRS}\lambda)^2-8p_{\sDRS}(\lambda-2)}
	\),
	in which case
	\begin{equation}\label{eq:DRE:c1}
		c
	{}={}
		\tfrac{2-\lambda}{2\lambda\gamma}
		{}+{}
		\tfrac{\sigma_{\sDRS}}{\lambda}
		(\tfrac12-\tfrac{\gamma L_{\sDRS}}{\lambda}).
	\end{equation}
	\begin{proof}
		\newcommand\x\xi%
		Let \((u^+,v^+)\) be generated by one \ref{DRS} iteration starting at \(s^+\).
		Then,
		\begin{align*}
			\DRE(s^+)
		{}={} &
			\min_{w\in\R^n}\set{
				\sDRS(u^+)
				{}+{}
				\nsDRS(w)
				{}+{}
				\innprod{\nabla\sDRS(u^+)}{w-u^+}
				{}+{}
				\tfrac{1}{2\gamma}
				\|w-u^+\|^2
			}
		\end{align*}
		and the minimum is attained at \(w=v^+\).
		Therefore, letting \(\rho\) be as in \cref{thm:rho},
		{\mathtight[0.9]%
			\begin{align*}
				\DRE(s^+)
			{}\leq{} &
				\sDRS(u^+)
				{}+{}
				\innprod{\nabla\sDRS(u^+)}{v-u^+}
				{}+{}
				\nsDRS(v)
				{}+{}
				\tfrac{1}{2\gamma}
				\|u^+-v\|^2
			\\
			{}={} &
				\smashunderbracket[0.5pt]{
					\sDRS(u^+)
					{}+{}
					\innprod{\nabla\sDRS(u^+)}{u-u^+}
				}{}
				{}+{}
				\innprod{\nabla\sDRS(u^+)}{v-u}
				{}+{}
				\nsDRS(v)
				{}+{}
				\tfrac{1}{2\gamma}
				\|u^+-v\|^2
			\\
			{}\overrel[\leq]{\llap{\cref{thm:rho}}}{} &
				\smashoverbracket[0.5pt]{
					\sDRS(u)
					{}-{}
					\rho(u,u^+)
				}{}
				{}+{}
				\innprod{\nabla\sDRS(u^+)}{v-u}
				{}+{}
				\nsDRS(v)
				{}+{}
				\tfrac{1}{2\gamma}
				\|u^+-v\|^2
			\\
			{}={} &
				\sDRS(u)
				{}-{}
				\rho(u,u^+)
				{}+{}
				\innprod{\nabla\sDRS(u)}{v-u}
				{}+{}
				\nsDRS(v)
				{}+{}
				\tfrac{1}{2\gamma}
				\|u^+-v\|^2
			\\
			&
				{+{}}
				\innprod{\nabla\sDRS(u^+)-\nabla\sDRS(u)}{v-u}
			\\
			{}={} &
				\DRE(s)
				{}-{}
				\rho(u,u^+)
				{}+{}
				\innprod{\nabla\sDRS(u^+)-\nabla\sDRS(u)}{v-u}
				{}+{}
				\tfrac{1}{2\gamma}
				\|u-u^+\|^2
				{}+{}
				\tfrac1\gamma
				\innprod{u^+-u}{u-v}.
			\end{align*}%
		}%
		Since
		\(
			u-v
		{}={}
			\tfrac1\lambda(s-s^+)
		{}={}
			\tfrac1\lambda(u-u^+)
			{}+{}
			\tfrac\gamma\lambda(\nabla\sDRS(u)-\nabla\sDRS(u^+))
		\),
		as it follows from \cref{thm:proxfEquiv}, it all simplifies to
		\begin{equation}\label{eq:DRE:SD0}
			\DRE(s)
			{}-{}
			\DRE(s^+)
		{}\geq{}
			\tfrac{2-\lambda}{2\gamma\lambda}
			\|u-u^+\|^2
			{}-{}
			\tfrac\gamma\lambda
			\|\nabla\sDRS(u^+)-\nabla\sDRS(u)\|^2
			{}+{}
			\rho(u,u^+).
		\end{equation}
		It will suffice to show that
		\[
			\DRE(s)-\DRE(s^+)
		{}\geq{}
			c
			\|u-u^+\|^2;
		\]
		inequality \eqref{eq:DRE:SD} will then follow from the \(\tfrac{1}{1+\gamma L_{\sDRS}}\)-strong monotonicity of \(\prox_{\gamma\sDRS}\), see \cref{thm:proxfLip}.
		We now proceed by cases.
		\begin{proofitemize}
		\item\textbf{Case 1: \(\lambda\in(0,2)\).}
			
			\noindent
			Let \(\sigma\coloneqq-[\sigma_{\sDRS}]_-=\min\set{\sigma_{\sDRS},0}\) and \(L\geq L_{\sDRS}\) be such that \(L+\sigma>0\); the value of such an \(L\) will be fixed later.
			Then, \(\sigma\leq0\) and \(\sDRS\) is \(L\)-smooth and \(\sigma\)-hypoconvex.
			We may thus choose \(\rho(u,u^+)\) as in \cref{thm:hypo:rho1} with these values of \(L\) and \(\sigma\).
			Inequality \eqref{eq:DRE:SD0} then becomes
			\[
				\tfrac{\DRE(s)-\DRE(s^+)}{L}
			{}\geq{}
				\left(
					\tfrac{2-\lambda}{2\lambda\x}
					{}+{}
					\tfrac{p}{2(1+p)}
				\right)
				\|u^+-u\|^2
				{}+{}
				\tfrac{1}{L^2}
				\left(
					\tfrac{1}{2(1+p)}
					{}-{}
					\tfrac{\x}{\lambda}
				\right)
				\|\nabla\sDRS(u^+)-\nabla\sDRS(u)\|^2,
			\]
			where \(\x\coloneqq\gamma L\) and \(p\coloneqq\nicefrac\sigma L\in(-1,0]\).
			Since \(\nabla\sDRS\) is \(L_{\sDRS}\)-Lipschitz continuous, the claim holds provided that the constant
			\begin{equation}\label{eq:c}
				\frac cL
			{}={}
				\begin{ifcases}
					\tfrac{2-\lambda}{2\lambda\x}
					{}+{}
					\tfrac{p}{2(1+p)}
				&
					0
					{}<{}
					\tfrac{1}{2(1+p)}
					{}-{}
					\tfrac{\x}{\lambda}
				\\[3pt]
					\tfrac{2-\lambda}{2\lambda\x}
					{}+{}
					\tfrac{p}{2(1+p)}
					{}+{}
					\tfrac{L_{\sDRS}^2}{L^2}
					\left(
						\tfrac{1}{2(1+p)}
						{}-{}
						\tfrac{\x}{\lambda}
					\right)
				\otherwise[otherwise,]
				\end{ifcases}
			\end{equation}
			is strictly positive.
			Now, let us consider two subcases:
			\begin{itemize}
			\item\textbf{Case 1a:} \(0<\lambda\leq 2(1+\nicefrac{\sigma}{L_{\sDRS}})\).
			
				Then,
				\(
					\sigma\geq -\frac{2-\lambda}{2}L_{\sDRS}>-L_{\sDRS}
				\)
				and we can take \(L=L_{\sDRS}\).
				Consequently, \(p=\nicefrac{\sigma}{L_{\sDRS}}\), \(\x=\gamma L_{\sDRS}\), and \eqref{eq:c} becomes
				\begin{equation}\label{eq:cbar}
					\frac{c}{L_{\sDRS}}
				{}={}
					\frac{2-\lambda}{2\lambda\gamma L_{\sDRS}}
					{}+{}
					\begin{ifcases}
						\tfrac{p}{2(1+p)}
					&
						\gamma
						{}<{}
						\tfrac{\lambda}{2(1+p)}
					\\[3pt]
						\tfrac{1}{2}
						{}-{}
						\tfrac{\gamma L_{\sDRS}}{\lambda}
					\otherwise.
					\end{ifcases}
				\end{equation}
				Let us verify that in this case any \(\gamma\) such that \(\gamma<\nicefrac{1}{L_{\sDRS}}\) yields a strictly positive coefficient \(c\).
				If \(0<\gamma L_{\sDRS}<\frac{\lambda}{2(1+p)}\leq 1\), then
				\[
					\tfrac{c}{L_{\sDRS}}
				{}={}
					\tfrac{2-\lambda}{2\lambda\gamma L_{\sDRS}}
					{}+{}
					\tfrac{p}{2(1+p)}
				{}>{}
					\tfrac{2-\lambda}{2\lambda}
					{}+{}
					\tfrac{p}{\lambda}
				{}={}
					\tfrac{1+p}{\lambda}
					{}-{}
					\tfrac{1}{2}
				{}\geq{}
					0,
				\]
				where in the inequality we used the fact that \(\lambda<2\) and \(p\leq0\).
				If instead \(\frac{\lambda}{2(1+p)}<\gamma L_{\sDRS}<1\), then
				\[
					\tfrac{c}{L_{\sDRS}}
				{}={}
					\tfrac{2-\lambda}{2\lambda\gamma L_{\sDRS}}
					{}+{}
					\tfrac{1}{2}
					{}-{}
					\tfrac{\gamma L_{\sDRS}}{\lambda}
				{}>{}
					\tfrac{2-\lambda}{2\lambda}
					{}+{}
					\tfrac{1}{2}
					{}-{}
					\tfrac{1}{\lambda}
				{}={}
					0.
				\]
				Either way, the sufficient decrease constant \(c\) is strictly positive.
				Since \(\sigma=-[\sigma_{\sDRS}]_-\) and
				\[
					\tfrac{2-\lambda}{2\lambda\gamma}
					{}+{}
					\tfrac{\sigma}{2(1+p)}
				{}\leq{}
					\tfrac{2-\lambda}{2\lambda\gamma}
					{}+{}
					\tfrac{L_{\sDRS}}{2}
					{}-{}
					\tfrac{\gamma L_{\sDRS}^2}{\lambda}
				\quad\Leftrightarrow\quad
					\gamma
				{}\leq{}
					\tfrac{\lambda}{2(L_{\sDRS}+\sigma)},
				\]
				from \eqref{eq:cbar} we conclude that \(c\) is as in \eqref{eq:DRE:SD}.
			\item\textbf{Case 1b:} \(2(1+\nicefrac{\sigma}{L_{\sDRS}})<\lambda<2\).
				
				Necessarily \(\sigma<0\), for otherwise the range of \(\lambda\) would be empty.
				In particular, \(\sigma=\sigma_{\sDRS}\), and the lower bound on \(\lambda\) can be expressed as
				\(
					\sigma_{\sDRS}<-\frac{2-\lambda}{2}L_{\sDRS}
				\).
				Consequently, \(L\coloneqq\tfrac{-2\sigma_{\sDRS}}{2-\lambda}\) is strictly larger than \(L_{\sDRS}\), and in particular \(\sigma+L=\sigma_{\sDRS}+L>0\).
				The ratio of \(\sigma\) and \(L\) is thus \(p=\tfrac\lambda2-1\), and \eqref{eq:c} becomes
				\begin{equation}\label{eq:c2}
					c
				{}={}
					\tfrac{2-\lambda}{2\lambda\gamma}
					{}+{}
					\begin{ifcases}
						\tfrac{\sigma_{\sDRS}}{\lambda}
					&
						\gamma
						{}<{}
						\tfrac{2-\lambda}{-2\sigma_{\sDRS}}
					\\[3pt]
						\tfrac{\sigma_{\sDRS}}{\lambda}-\tfrac{\gamma L_{\sDRS}^2}{\lambda}+\tfrac{2-\lambda}{-2\sigma_{\sDRS}\lambda}L_{\sDRS}^2
					\otherwise.
					\end{ifcases}
				\end{equation}
				Let us show that, when
				\(
					\gamma
				{}<{}
					\tfrac{2-\lambda}{-2\sigma_{\sDRS}}
				{}={}
					\tfrac1L
				\),
				also in this case the sufficient decrease constant \(c\) is strictly positive.
				We have
				\[
					\tfrac cL
				{}={}
					\tfrac{2-\lambda}{2\lambda\gamma L}
					{}+{}
					\tfrac{\sigma_{\sDRS}}{\lambda}
					\tfrac1L
				{}>{}
					\tfrac{2-\lambda}{2\lambda}
					{}+{}
					\tfrac{\sigma_{\sDRS}}{\lambda}
					\tfrac{2-\lambda}{-2\sigma_{\sDRS}}
				{}={}
					0,
				\]
				hence the claim.
				This concludes the proof for the case \(\lambda\in(0,2)\).
			\end{itemize}
		\item\textbf{Case 2: \(\lambda\geq2\).}

			\noindent
			In this case we need to assume that \(\sDRS\) is strongly convex, that is, that \(\sigma_{\sDRS}>0\).
			Instead of considering a single expression of \(\rho\), we will rather take a convex combination of those in \cref{thm:hypo:rho0,thm:hypo:rho1}, namely
			\[
				\rho(u,u^+)
			{}={}
				(1-\alpha)\tfrac{\sigma_{\sDRS}}{2}\|u-u^+\|^2
				{}+{}
				\alpha\tfrac{1}{2L_{\sDRS}}\|\nabla\sDRS(u)-\nabla\sDRS(u^+)\|^2
			\]
			for some \(\alpha\in[0,1]\) to be determined.
			\eqref{eq:DRE:SD0} then becomes
			\[
				\tfrac{
					\DRE(s)-\DRE(s^+)
				}{
					L_{\sDRS}
				}
			{}\geq{}
				\left(
					\tfrac{2-\lambda}{2\lambda\xi}
					{}+{}
					\tfrac{(1-\alpha)p}{2}
				\right)
				\|u-u^+\|^2
				{}+{}
				\tfrac{1}{L_{\sDRS}^2}
				\left(
					\tfrac{\alpha}{2}
					{}-{}
					\tfrac{\xi}{\lambda}
				\right)
				\|\nabla\sDRS(u)-\nabla\sDRS(u^+)\|^2,
			\]
			where \(\x\coloneqq\gamma L_{\sDRS}\) and \(p\coloneqq\nicefrac{\sigma_{\sDRS}}{L_{\sDRS}}\in(0,1]\).
			By restricting \(\xi\in(0,1)\), since \(\lambda\geq2\) one can take
			\(
				\alpha
			{}\coloneqq{}
				\nicefrac{2\xi}{\lambda}
			{}\in{}
				(0,1)
			\)
			to make the coefficient multiplying the gradient norm vanish.
			We then obtain
			\begin{equation}\label{eq:clambda2}
				\tfrac{c}{L_{\sDRS}}
			{}={}
				\tfrac{2-\lambda}{2\lambda\xi}
				{}+{}
				\tfrac{(\lambda-2\xi)p}{2\lambda}.
			\end{equation}
			Imposing \(c>0\) results in the following second-order equation in variable \(\xi\),
			\begin{equation}\label{eq:xi}
				2p\xi^2-p\lambda\xi+(\lambda-2)
			{}<{}
				0.
			\end{equation}
			The discriminant is
			\(
				\Delta
			{}\coloneqq{}
				(p\lambda)^2
				{}-{}
				8p(\lambda-2)
			\),
			which, for \(\lambda\geq2\), is strictly positive iff
			\[
				2
			{}\leq{}
				\lambda
			{}<{}
				\tfrac{4}{1+\sqrt{1-p}}
			~~\vee~~
				\lambda
			{}>{}
				\tfrac{4}{1-\sqrt{1-p}}.
			\]
			Denoting
			\(
				\delta
			{}\coloneqq{}
				\sqrt\Delta
			{}={}
				\sqrt{
					(p\lambda)^2
					{}-{}
					8p(\lambda-2)
				}
			\),
			the solution to \eqref{eq:xi} is
			\(
				\tfrac{
					p\lambda
					{}-{}
					\delta
				}{
					4p
				}
			{}<{}
				\x
			{}<{}
				\tfrac{
					p\lambda
					{}+{}
					\delta
				}{
					4p
				}
			\).
			However, the case \(\lambda\geq 4\) has to be discarded, as
			\(
				\tfrac{p\lambda-\delta}{4p}
			{}>{}
				1
			\)
			in this case, contradicting the fact that \(p\leq1\).
			To see this, suppose \(\lambda\geq4\).
			Then,
			\begin{align*}
				\tfrac{
					p\lambda
					{}-{}
					\delta
				}{
					4p
				}
			{}<{}
				1
			~~\Leftrightarrow~~ &
				p(\lambda-4)
			{}<{}
				\delta
			\\
			~~\Leftrightarrow~~ &
				p^2(\lambda-4)^2
			{}<{}
				\Delta
			{}={}
				(p\lambda)^2
				{}-{}
				8p(\lambda-2)
			\\
			~~\Leftrightarrow~~ &
				p(2-\lambda)
			{}<{}
				2-\lambda,
			\end{align*}
			hence \(p>1\), which contradicts the fact that \(\sigma_{\sDRS}\leq L_{\sDRS}\).
			Thus, the only feasible ranges are the ones given in \eqref{eq:lambda>=2}, hence the claimed sufficient decrease constant \(c\), cf. \eqref{eq:clambda2}.
			\qedhere
		\end{proofitemize}
	\end{proof}
\end{thm}

\begin{rem}[Simpler bounds for \ref{DRS}]\label{rem:DRS}%
	By using the (more conservative) estimate \(\sigma_{\sDRS}=0\) when the smooth function \(\sDRS\) is convex, and \(\sigma_{\sDRS}=-L_{\sDRS}\) otherwise, the range of \(\gamma\) can be simplified as follows in case \(\lambda\in(0,2]\):
	\begin{qedalign*}
		\lambda\in(0,2)
	~&
		\begin{cases}[l @{~\text{ and }~} l @{~~~} l]
			\gamma<\frac{1}{L_{\sDRS}}
		&
			c
			{}={}
			\frac{2-\lambda}{2\lambda\gamma}
			{}-{}
			L_{\sDRS}[\nicefrac12-\gamma L_{\sDRS}]_+
		&
			\text{if \(\sDRS\) is convex,}
		\\[5pt]
			\gamma<\frac{2-\lambda}{2L_{\sDRS}}
		&
			c
			{}={}
			\frac{2-\lambda}{2\lambda\gamma}
			{}-{}
			\frac{L_{\sDRS}}{\lambda}
		&
			\text{otherwise.}
		\end{cases}
	\\[5pt]
		\lambda=2
	~&
		\begin{cases}[l c l @{~~~} l]
			\gamma<\frac{1}{L_{\sDRS}}
		&~\text{ and }~~&
			c
			{}={}
			\frac{\sigma_{\sDRS}}{2}(1-\gamma L_{\sDRS})
		&
			\text{if \(\sDRS\) is strongly convex,}
		\\[5pt]
			\hphantom{
				\gamma<\frac{2-\lambda}{2L_{\sDRS}}
			}
		&\emptyset&
			\hphantom{
				c
				{}={}
				\frac{2-\lambda}{2\lambda\gamma}
				{}-{}
				L_{\sDRS}[\nicefrac12-\gamma L_{\sDRS}]_+
			}
		&
			\text{otherwise.}
		\end{cases}
	\end{qedalign*}
\end{rem}


\grayout{%
	\begin{lem}\label{thm:DRS:EB}%
		Suppose that \Cref{ass:DRS} is satisfied, and let \(\gamma<\nicefrac{1}{L_{\sDRS}}\) be fixed.
		Then
		\[
			\dist(0,\hat\partial\varphi(v))
		{}\leq{}
			\tfrac{1-\gamma\sigma_{\sDRS}}{\gamma}\|u-v\|
		\qquad
			\text{for all \(u\in\R^p\) and \(v\in\FB u\).}
		\]
		In particular, this is true for any pair \((u,v)\) generated by \ref{DRS} with stepsize \(\gamma\).
		\begin{proof}
			Let \(M(w,u)\) be the function in the minimization subproblem of \eqref{eq:DRE}, so that \(v\in\argmin_wM(w,u)\).
			The optimality of \(v\) then implies that \(0\in\hat\partial_wM(v,u)\), where \(\hat\partial_w\) denotes the subdifferential of the section \(w\mapsto M(w,u)\).
			Using the subdifferential calculus rule of \cite[Ex. 8.8]{rockafellar2011variational}, we obtain
			\(
				0
			{}\in{}
				\hat\partial_wM(v,u)
			{}={}
				\nabla\sDRS(u)
				{}+{}
				\hat\partial\nsDRS(v)
				{}+{}
				\tfrac1\gamma(v-u)
			\).
			Denoting
			\(
				\delta\coloneqq\Fw{}
			\),
			the inclusion can be expressed as
			\(
				\nabla\delta(u)
				{}-{}
				\nabla\delta(v)
			{}\in{}
				\hat\partial\varphi(v)
			\).
			Notice that \(\delta\) is \(L_\delta\)-smooth, with \(L_\delta=\tfrac{1-\gamma\sigma_f}{\gamma}\).
			In particular,
			\[
				\dist(0,\hat\partial\varphi(v))
			{}\leq{}
				\|\nabla\delta(u)-\nabla\delta(v)\|
			{}\leq{}
				\tfrac{1-\gamma\sigma_{\sDRS}}{\gamma}\|u-v\|,
			\]
			hence the claimed inequality.
			The proof then follows by invoking \eqref{eq:vFBS}.
		\end{proof}
	\end{lem}
}%
\begin{thm}[Subsequential convergence]\label{thm:DRS}%
	Suppose that \Cref{ass:DRS} is satisfied, and consider a sequence \(\seq{s^k,u^k,v^k}\) generated by \ref{DRS} with stepsize \(\gamma\) and relaxation \(\lambda\) as in \Cref{thm:DRE:SD}, starting from \(s^0\in\R^p\).
	The following hold:%
	\begin{enumerate}
	\item\label{thm:DRS:r}
		The residual \(\seq{u^k-v^k}\) vanishes with rate
		\(
			\min_{i\leq k}\|u^i-v^i\|
		{}={}
			o(\nicefrac{1}{\sqrt k})
		\).
	\item\label{thm:DRS:uv}
		\(\seq{u^k}\) and \(\seq{v^k}\) have same cluster points, all of which are stationary for \(\varphi\) and on which \(\varphi\) has same value, this being the limit of \(\seq{\DRE(s^k)}\).
		In fact, for each \(k\) one has
		\(
			\dist(0,\hat\partial\varphi(v^k))
		{}\leq{}
			\frac{1-\gamma\sigma_{\sDRS}}{2\gamma}
			\|u^k-v^k\|
		\).
	\item\label{thm:DRS:bounded}
		If \(\varphi\) has bounded level sets, then the sequence \(\seq{s^k,u^k,v^k}\) is bounded.
	\end{enumerate}
	\begin{proof}
		To avoid trivialities, we assume that a fixed point is not found in a finite number of iterations, hence that \(v^k\neq u^k\) for all \(k\)'s.
		\begin{proofitemize}
		\item\ref{thm:DRS:r}
			Let \(c=c(\gamma,\lambda)\) be as in \cref{thm:DRE:SD}.
			Telescoping the inequality \eqref{eq:DRE:SD} yields
			\begin{align*}
			\textstyle
				\tfrac{c\lambda^2}{(1+\gamma L_{\sDRS})^2}
				\sum_{k\in\N}{
					\|u^k-v^k\|^2
				}
			{}\leq{} &
			\textstyle
				\sum_{k\in\N}{
				\bigl[
					\DRE(s^k)
					{}-{}
					\DRE(s^{k+1})
				\bigr]}
			{}\leq{}
				\DRE(s^0)-\inf\DRE.
			\end{align*}
			Since \(\inf\DRE=\inf\varphi>-\infty\) and \(\DRE\) is real valued (cf. \cref{thm:C0,thm:DREequiv}), it follows that \(\seq{u^k-v^k}\) is square summable, hence the claimed rate of convergence.
			Moreover, since \(\DRE(s^k)\) is decreasing it admits a (finite) limit, be it \(\varphi_\star\).
		\item\ref{thm:DRS:uv}
			Since \(\seq{u^k-v^k}\to 0\), necessarily \(\seq{u^k}\) and \(\seq{v^k}\) have same cluster points.
			Suppose that \(\seq{u^k}[k\in K]\to u'\) for some \(K\subseteq\N\) and \(u'\in\R^p\).
			Then, \(\seq{v^k}[k\in K]\to u'\), and since \(s^k=u^k+\nabla\sDRS(u^k)\) (cf. \cref{thm:proxfEquiv}), continuity of \(\nabla\sDRS\) implies that
			\(
				\seq{s^k}[k\in K]
			{}\to{}
				s'
			{}={}
				u'+\gamma\nabla\sDRS(u')
			\).
			From \cref{thm:proxfEquiv} we infer that \(u'=\prox_{\gamma\sDRS}(s')\).
			
			Similarly, \(\seq{\Fw{u^k}}[k\in K]\to \Fw{u'}\), and the outer semicontinuity of \(\prox_{\gamma\nsDRS}\) \cite[Ex. 5.23(b)]{rockafellar2011variational} combined with \eqref{eq:vFBS} implies that
			\[
				u'
			{}={}
				\lim_{K\ni k\to\infty}{
					\hspace*{-5pt}
					v^k
				}
			{}\in{}
				\limsup_{K\ni k\to\infty}{
					\FB{u^k}
				}
			{}\subseteq{}
				\FB{u'}.
			\]
			From \eqref{eq:prox:subdiff} we then have that
			\(
				-\nabla\sDRS(u')
			{}\in{}
				\hat\partial\nsDRS(u')
			\),
			hence \(0\in\hat\partial\varphi(u')\), as it follows from \cite[Ex. 8.8]{rockafellar2011variational}.
			Finally, since \(v^k\to u'\),
			\[
				\varphi(u')
			{}\leq{}
				\liminf_{K\ni k\to\infty}{
					\varphi(v^k)
				}
			{}\leq{}
				\limsup_{K\ni k\to\infty}{
					\varphi(v^k)
				}
			{}\leq{}
				\limsup_{K\ni k\to\infty}{
					\DRE(s^k)
				}
			{}={}
				\DRE(s')
			{}\leq{}
				\varphi(u'),
			\]
			where the first inequality is due to lower semicontinuity of \(\varphi\), the third and the last to the sandwiching property (\cref{thm:sandwich}), and the equality to the continuity of \(\DRE\) (\cref{thm:C0}).
			This shows that \(\seq{\varphi(u^k)}[k\in K]\to\varphi(u')=\DRE(s')\), and since \(\seq{\DRE(s^k)}\to\varphi_\star\), then necessarily \(\varphi(u')=\DRE(s')=\varphi_\star\) independently of the cluster point \(u'\).
			The last assert follows from the optimality condition of \(v^k\) as in \eqref{eq:vFBS}, namely,
			\(
				\tfrac1\gamma(u^k-v^k)
				{}-{}
				\nabla\sDRS(u^k)
			{}\in{}
				\hat\partial\nsDRS(v^k)
			\)
			due to \eqref{eq:prox:subdiff}, together with \(\sigma_{\sDRS}\)-hypoconvexity of \(\sDRS\).
		\item\ref{thm:DRS:bounded}.
			Suppose that \(\varphi\) has bounded level sets.
			Then, it follows from \cref{thm:boundedlev} that so does \(\DRE\), and since
			\(
				s^k
			{}\in{}
				\lev_{\leq\DRE(s^0)}\DRE
			\)
			for all \(k\in\N\), then the sequence \(\seq{s^k}\) is bounded.
			Due to Lipschitz continuity of \(\prox_{\gamma\sDRS}\) (cf. \cref{thm:proxfLip}), also \(\seq{u^k}\) is bounded.
			In turn, since \(v^k-u^k\to 0\) we conclude that also \(\seq{v^k}\) is bounded.
		\qedhere
		\end{proofitemize}
	\end{proof}
\end{thm}

The Kurdyka-\L ojasiewicz (KL) property is a powerful tool to establish global convergence (as opposed to subsequential convergence) of descent methods, see \cite{attouch2013convergence}, and semialgebraic functions comprise a wide class of functions that enjoy this property.
It was first observed in \cite{li2016douglas} that the augmented Lagrangian decreases along iterates generated by non-relaxed \ref{DRS}, cf. \eqref{eq:LPong}, and global convergence was thus established when \(\sDRS\) and \(\nsDRS\) are semialgebraic functions and the sequence remains bounded.
The latter requirement was later shown to hold in \cite{li2017peaceman} when \(\varphi\) has bounded level sets, as \Cref{thm:boundedlev} confirms.
Due to the equivalence of the DRE and the augmented Lagrangian evaluated at points generated by \ref{DRS}, cf. \eqref{eq:DRE=L}, by invoking \Cref{thm:DRE:SD} we can extend their result to the tight ranges we provided.
\begin{thm}[{Global convergence of \ref{DRS} \cite[Thm. 2]{li2016douglas}}]\label{thm:DRSglobal}%
	Suppose that \Cref{ass:DRS} is satisfied, that \(\varphi\) is level bounded, and that \(\sDRS\) and \(\nsDRS\) are semialgebraic.
	Then, the sequences \(\seq{u^k}\) and \(\seq{v^k}\) generated by \ref{DRS} with \(\gamma\) and \(\lambda\) as in \Cref{thm:DRS} converge to (the same) stationary point for \(\varphi\).
\end{thm}

		\subsection{Adaptive variant}
			\label{sec:Adaptive}%
			As described in \Cref{rem:DRS}, when the hypoconvexity modulus \(\sigma_{\sDRS}\) is not known one can always consider \(\sigma_{\sDRS}=-L_{\sDRS}\); in case \(\sDRS\) is convex, the tighter estimate \(\sigma_{\sDRS}=0\) is also feasible.
In particular, for any \(\lambda\in(0,2)\) the knowledge of \(L_{\sDRS}\) is enough for determining ranges of \(\gamma\), although possibly conservative, that comply with \Cref{thm:DRE:SD} and thus make \ref{DRS} iterations convergent.

When also the Lipschitz constant \(L_{\sDRS}\) is not readily available, it is however possible to adjust the stepsize \(\gamma\) along the iterations without losing the convergence properties of \Cref{thm:DRS}.
This can be done by selecting an initial estimate \(\gamma\) for the stepsize, and reduce it whenever a sufficient decrease condition is violated.
Due to the fact that \(\gamma\) may be larger than the unknown threshold \(\nicefrac{1}{[\sigma_{\sDRS}]_-}\), below which \(\prox_{\gamma\sDRS}\) is ensured to be single valued (cf. \cref{thm:proxf}), the DRE may fail to be a well-defined function of \(s\).
For this reason, we resort to the augmented Lagrangian interpretation given in \eqref{eq:DRE=L}.

At each iteration, the stepsize \(\gamma\) is reduced whenever a sufficient decrease condition on the augmented Lagrangian is violated.
This can happen only a finite number of times, since for \(\gamma\) small enough \eqref{eq:DRE=L} holds and the sufficient decrease property as stated in \Cref{thm:DRE:SD} applies.
It may also be the case that \(\gamma\) remains high and lower boundedness cannot be inferred from \Cref{thm:inf}.
To prevent the augmented Lagrangian from dropping arbitrarily low, we may thus enforce a bound similar to that of \Cref{thm:DREgeq} so as to keep it above \(\inf\varphi\).
This is a feasible requirement, since as soon as \(\gamma\) falls below \(\nicefrac{1}{L_{\sDRS}}\) the statement of \Cref{thm:DREgeq} applies.

The procedure is summarized in \Cref{alg:DRSadaptive}.
Note that, apart from the re-evaluation of \(u^k\) and \(v^k\) whenever \(\gamma\) is decreased, the adaptive variant comes at the additional cost of computing \(\sDRS(u^k)\), \(\sDRS(v^k)\), and \(\nsDRS(v^k)\) at each iteration, needed for the test at \cref{step:DRSadaptive:if}.
\begin{algorithm}[htb]
	\algcaption[\ref{DRS} with adaptive stepsize]{%
		\ref{DRS} with adaptive stepsize.\\
		\(\LL_\beta\) is the augmented Lagrangian as defined in \eqref{eq:DRS-L}.\\
		\(\ffunc{\text{\ref{DRS}}_{\gamma,\lambda}}{\R^p}{\R^p\times\R^p\times\R^p}\) maps \(s\in\R^p\) to a triplet \((u,v,s^+)\) as in \eqref{DRS}.%
	}%
	\label{alg:DRSadaptive}%
	\begin{algorithmic}[1]
	\setlength\baselineskip{1.25\baselineskip}
	\Require
		\(s^0\in\R^p\),~
		\(L>0\),~
		\(\lambda\in(0,2),\gamma,c\) as in \cref{rem:DRS} with \(L\) in place of \(L_{\sDRS}\)
	\Initialize
		\(
			(u^0,v^0,s^1)
		{}\in{}
			{\rm\ref{DRS}}_{\gamma,\lambda}(s^0)
		\),~
		\(
			\mathcal L_0
		{}={}
			\LL_{\nicefrac1\gamma}(u^0,v^0,\gamma^{-1}(u^0-s^0))
		\)
	\item[{\bf For} ~\(k=1,2,\ldots\)~ {\bf do}]
	\State\label{step:DRSadaptive:1}%
		\(
			(u^k,v^k,s^{k+1})
		{}\in{}
			{\rm\ref{DRS}}_{\gamma,\lambda}(s^k)
		\)
	\Statex
		\(
			\mathcal L_k
		{}={}
			\LL_{\nicefrac1\gamma}(u^k,v^k,\gamma^{-1}(u^k-s^k))
		\)
	\If{~
		\(
			\mathcal L_k
		{}>{}
			\mathcal L_{k-1}
			{}-{}
			\tfrac{c\lambda^2}{(1+\gamma L)^2}
			\|v^{k-1}-u^{k-1}\|^2
		\)
		~~{\bf or}~~
		\(
			\varphi(v^k)
		{}>{}
			\mathcal L_k
		\)
	~}\label{step:DRSadaptive:if}%
		\State%
			\(\gamma\gets\nicefrac\gamma2\),~
			\(c\gets2c\),~
			\(L\gets2L\)
		\Statex\hspace*{\algorithmicindent}%
			\(
				(u^{k-1},v^{k-1},s^k)
			{}\in{}
				{\rm\ref{DRS}}_{\gamma,\lambda}(s^{k-1})
			\)
		\Statex\hspace*{\algorithmicindent}%
			\(
				\mathcal L_{k-1}
			{}\gets{}
				\LL_{\nicefrac1\gamma}(u^{k-1},v^{k-1},\gamma^{-1}(u^{k-1}-s^{k-1}))
			\)
			and go back to \cref{step:DRSadaptive:1}
	\EndIf{}
\end{algorithmic}
\end{algorithm}

\grayout{%
Due to the possible multi valuedness of \(\prox_{\gamma\sDRS}\), the DRE may fail to be well defined for \(\gamma\geq\nicefrac{1}{[\sigma_{\sDRS}]_-}\). 
Nevertheless, the following \emph{local} version of \Cref{thm:sandwich} allows to extend the convergence results of \Cref{thm:DRS} to this adaptive variant as well.
\begin{thm}[Sandwich property: the general case]\label{thm:LLsandwich}%
	Suppose that \Cref{ass:DRS} is satisfied, and let \(\gamma<\min\set{\gamma_{\sDRS},\gamma_{\nsDRS}}\).
	For all \(s\in\R^p\), \(u\in\prox_{\gamma\sDRS}(s)\), and \(v\in\prox_{\gamma\nsDRS}(2u-s)\) the following hold:
	\begin{enumerate}
	\item\label{thm:LLleq}
		\(
			\LL_{\nicefrac1\gamma}\bigl(u,v,\gamma^{-1}(u-s)\bigr)
		{}\leq{}
			\varphi(u)
		\).
	\item\label{thm:LLgeq}
		Let \(L=L(s,u,v)\) be such that
		\(
			\sDRS(v)
		{}\leq{}
			\sDRS(u)
			{}+{}
			\tfrac1\gamma\innprod{s-u}{v-u}
			{}+{}
			\tfrac L2\|v-u\|^2
		\);
		then,
		\[
			\varphi(v)
		{}\leq{}
			\LL_{\nicefrac1\gamma}\bigl(u,v,\gamma^{-1}(u-s)\bigr)
			{}-{}
			\tfrac{1-\gamma L}{2\gamma}\|v-u\|^2.
		\]
	\end{enumerate}
	Moreover, if \(\gamma<\nicefrac{1}{[\sigma_{\sDRS}]_-}\), then \(u\) is uniquely determined and
	\[
		\LL_{\nicefrac1\gamma}\bigl(u,v,\gamma^{-1}(u-s)\bigr)
	{}={}
		\FBE(u)
	{}={}
		\DRE(s).
	\]
	\begin{proof}
		Since \(v\in\prox_{\gamma\nsDRS}(2u-s)\), we have
		\begin{align*}
			\nsDRS(u)
			{}+{}
			\tfrac{1}{2\gamma}\|s-u\|^2
		{}\geq{} &
			\nsDRS(v)
			{}+{}
			\tfrac{1}{2\gamma}\|v-2u+s\|^2
		\\
		{}={} &
			\nsDRS(v)
			{}+{}
			\tfrac1\gamma
			\innprod{s-u}{v-u}
			{}+{}
			\tfrac{1}{2\gamma}\|v-u\|^2
			{}+{}
			\tfrac{1}{2\gamma}\|s-u\|^2.
		\end{align*}
		By adding
		\(
			\sDRS(u)
			{}-{}
			\tfrac{1}{2\gamma}\|s-u\|^2
		\),
		inequality \ref{thm:LLleq} follows.
		Let now \(L\) be as in \ref{thm:LLgeq}.
		Then,
		\begin{align*}
			\LL_{\nicefrac1\gamma}\bigl(u,v,\gamma^{-1}(u-s)\bigr)
		{}={} &
			\sDRS(u)
			{}+{}
			\nsDRS(v)
			{}+{}
			\tfrac1\gamma\innprod{s-u}{v-u}
			{}+{}
			\tfrac{1}{2\gamma}\|v-u\|^2
		\\
		{}\geq{} &
			\sDRS(v)
			{}-{}
			\tfrac1\gamma\innprod{s-u}{v-u}
			{}+{}
			\nsDRS(v)
			{}+{}
			\tfrac1\gamma\innprod{s-u}{v-u}
			{}+{}
			\tfrac{1-\gamma L}{2\gamma}\|v-u\|^2
		\\
		{}={} &
			\varphi(v)
			{}+{}
			\tfrac{1-\gamma L}{2\gamma}\|v-u\|^2,
		\end{align*}
		which is \ref{thm:LLgeq}.			
		Finally, if \(\gamma<\nicefrac{1}{[\sigma_{\sDRS}]_-}\), then it follows from \cref{thm:proxfEquiv} that \(u\) is unique and characterized by the identity \(s=u+\gamma\nabla\sDRS(u)\).
		Thus, \(v\in\FB u\) and
		\begin{align*}
			\LL_{\nicefrac1\gamma}\bigl(u,v,\gamma^{-1}(u-s)\bigr)
		{}={} &
			\sDRS(u)
			{}+{}
			\nsDRS(v)
			{}+{}
			\innprod{\nabla\sDRS(u)}{v-u}
			{}+{}
			\tfrac{1}{2\gamma}\|v-u\|^2,
		\\
		{}={} &
			\min_{w\in\R^p}\set{
				\sDRS(u)
				{}+{}
				\nsDRS(w)
				{}+{}
				\tfrac{1}{2\gamma}\bigl\|w-\bigl(\Fw u\bigr)\bigr\|^2
				{}-{}
				\tfrac\gamma2\|\nabla\sDRS(u)\|^2
			},
		\end{align*}
		which is exactly \(\DRE(s)=\FBE(u)\), cf. \eqref{eq:DRE} and \cite[Rem. 4.1]{themelis2018forward}.
	\end{proof}
\end{thm}
}%
\begin{thm}[Subsequential convergence of adaptive \ref{DRS}]
	Suppose that \Cref{ass:DRS} is satisfied, and consider the iterates generated by \Cref{alg:DRSadaptive}.
	The following hold:
	\begin{enumerate}
	\item\label{thm:DRSadaptive:r}
		The residual \(\seq{u^k-v^k}\) vanishes with rate
		\(
			\min_{i\leq k}\|u^i-v^i\|
		{}={}
			o(\nicefrac{1}{\sqrt k})
		\).
	\item\label{thm:DRSadaptive:uv}
		\(\seq{u^k}\) and \(\seq{v^k}\) have same cluster points, all of which are stationary for \(\varphi\) and on which \(\varphi\) has same value, this being the limit of \(\seq{\mathcal L_k}\).
	\end{enumerate}
	\begin{proof}
		Note that the sufficient decrease constant in \cref{rem:DRS} satisfies \(c(\nicefrac\gamma2,2L)=2c(\gamma,L)\).
		Therefore, if \(L\geq L_{\sDRS}\) at iteration \(k\), then it follows from \eqref{eq:DRS-L} that \(\mathcal L_k=\DRE(s^k)\), and from \cref{thm:sandwich} we infer that the condition at \cref{step:DRSadaptive:if} is never passed.
		Therefore, starting from iteration \(k\) the stepsize \(\gamma\) is never decreased, and the algorithm reduces to plain (nonadaptive) \ref{DRS}.
		Either way, \(\gamma\) is decreased only a finite number of times; by possibly discarding the first iterates, without loss of generality we may assume that \(\gamma\) is constant (although possibly larger than or equal to \(\nicefrac{1}{L_{\sDRS}}\)).
		The iterates generated by \cref{alg:DRSadaptive} then satisfy
		\[
			\varphi(v^k)
		{}\leq{}
			\mathcal L_k
		\quad\text{and}\quad
			\mathcal L_{k+1}
		{}\leq{}
			\mathcal L_k
			{}-{}
			c'\|u^k-v^k\|^2
		\]
		for some constant \(c'>0\).
		In particular, \(\seq{\mathcal L_k}\) is lower bounded (by \(\inf\varphi\)), and by telescoping the second inequality we obtain that \(\seq{\|u^k-v^k\|}\) is square summable, hence the claimed rate.
		
		Since \(u^k-v^k\to0\), necessarily \(u^k\) and \(v^k\) have same cluster points.
		Suppose that a subsequence \(\seq{u^k}[k\in K]\) converges to a point \(u'\); then, so does \(\seq{v^k}[k\in K]\).
		Moreover, it follows from \cite[Ex. 10.2]{rockafellar2011variational} that \(s^k=u^k+\gamma\nabla\sDRS(u^k)\) (due to the fact that \(\gamma\) may be larger than \(\nicefrac{1}{[\sigma_{\sDRS}]_-}\), differently from the characterization given in \cref{thm:proxfEquiv} this condition is only necessary).
		Thus, for all \(k\)'s it holds that
		\(
			v^k
		{}\in{}
			\prox_{\gamma\sDRS}(2u^k-s^k)
		{}={}
			\FB{u^k}
		\).
		From the continuity of \(\nabla\sDRS\) and the outer semicontinuity of \(\prox_{\gamma\nsDRS}\), cf. \cite[Ex. 5.23(b)]{rockafellar2011variational}, it follows that the limit \(u'\) of \(\seq{v^k}[k\in K]\) satisfies
		\(
			u'
		{}\in{}
			\FB{u'}
		\),
		and the same reasoning as in the proof of \cref{thm:DRS:uv} shows that \(0\in\hat\partial\varphi(u')\).
		
		Finally, since \(\varphi(v^k)\leq\mathcal L_k\leq\mathcal L_0\), if \(\varphi\) is level bounded, then necessarily \(\seq{v^k}\) is bounded, hence so are \(\seq{u^k}\) and \(\seq{s^k}\) (since \(u^k-v^k\to 0\) and \(s^k=u^k+\gamma\nabla\sDRS(u^k)\)).
	\end{proof}
\end{thm}
Global convergence of adaptive \ref{DRS} again falls as a consequence of \cite[Thm. 2]{li2016douglas}.
\begin{thm}[Global convergence of adaptive \ref{DRS}]%
	Suppose that \Cref{ass:DRS} holds, that \(\varphi\) is level bounded, and that \(\sDRS\) and \(\nsDRS\) are semialgebraic.
	Then, the sequences \(\seq{u^k}\) and \(\seq{v^k}\) generated by adaptive \ref{DRS} (\cref{alg:DRSadaptive}) converge to (the same) stationary point of \(\varphi\).%
\end{thm}

		\subsection{Tightness of the results}
			\label{sec:Tightness}%
			When both \(\sDRS\) and \(\nsDRS\) are convex and \(\sDRS+\nsDRS\) attains a minimum, well-known results of monotone operator theory guarantee that for any \(\lambda\in(0,2)\) and \(\gamma>0\) the residual \(u^k-v^k\) generated by \ref{DRS} iterations vanishes (see \eg \cite[Cor. 28.3]{bauschke2017convex}).
In fact, the whole sequence \(\seq{u^k}\) converges and \(\sDRS\) needs not even be differentiable in this case.
On the contrary, when \(\nsDRS\) is nonconvex then the bound \(\gamma<\nicefrac{1}{L_{\sDRS}}\) plays a crucial role, as the next example shows.
\begin{thm}[Necessity of \(\gamma<\nicefrac{1}{L_{\sDRS}}\)]\label{thm:gamma}%
	For any \(L>0\) and \(\sigma\in[-L,L]\) there exist \(\func{\sDRS,\nsDRS}{\R^p}{\Rinf}\) satisfying the following properties
	\begin{enumeratprop}
	\item\label{thm:gamma:s}
		\(\sDRS\) is \(L\)-smooth and \(\sigma\)-hypoconvex;
	\item\label{thm:gamma:ns}
		\(\nsDRS\) is proper and lsc;
	\item\label{thm:gamma:argmin}
		\(\argmin(\sDRS+\nsDRS)\neq\emptyset\);
	\item
		for all \(s^0\in\R^p\), \(\gamma\geq\nicefrac1L\), and \(\lambda>0\), the sequence \(\seq{s^k}\) generated by \ref{DRS} iterations with stepsize \(\gamma\) and relaxation \(\lambda\) starting from \(s^0\) satisfies \(\|s^k-s^{k+1}\|\not\to0\) as \(k\to\infty\).
	\end{enumeratprop}
	\begin{proof}
		Fix \(t>1\), and let \(\varphi=\sDRS+\nsDRS\), where \(\nsDRS=\indicator_{\set{\pm1}}\) and
		\begin{equation}\label{eq:phi1_es}
			\sDRS(x)
		{}={}
			\begin{ifcases}
				\tfrac L2x^2
			&
				x\leq t
			\\[3pt]
				\tfrac L2x^2
				{}-{}
				\tfrac{L-\sigma}{2}(x-t)^2
			\otherwise.
			\end{ifcases}
		\end{equation}
		Notice that \(\dom\varphi=\set{\pm1}\), and therefore \(\pm1\) are the unique stationary points of \(\varphi\) (in fact, they are also global minimizers).
		It can be easily verified that \(\sDRS\) and \(\nsDRS\) satisfy properties \ref{thm:gamma:s}, \ref{thm:gamma:ns} and \ref{thm:gamma:argmin}.
		Moreover, \(\prox_{\gamma\sDRS}\) is well defined iff \(\gamma<\nicefrac{1}{[\sigma]_-}\), in which case
		\begin{equation}\label{eq:phi1_es_prox}
			\prox_{\gamma\sDRS}(s)
		{}={}
			\begin{ifcases}
				\tfrac{s}{1+\gamma L} & s\leq t(1+\gamma L)
			\\[3pt]
				\tfrac{s-\gamma(L-\sigma)t}{1+\gamma\sigma}
			\otherwise[otherwise,]
			\end{ifcases}
		\qquad\text{and}\qquad
			\prox_{\gamma\nsDRS}
		{}={}
			\sign,
		\end{equation}
		where \(\sign(0)=\set{\pm1}\).
		Let now \(s^0\in\R^p\), \(\nicefrac1L\leq\gamma<\nicefrac{1}{[\sigma]_-}\), and \(\lambda>0\) be fixed, and consider a sequence \(\seq{s^k}\) generated by \ref{DRS} with stepsize \(\gamma\) and relaxation \(\lambda\), starting at \(s^0\).
		To arrive to a contradiction, suppose that \(\|s^k-s^{k+1}\|=\lambda\|u^k-v^k\|\to0\) as \(k\to\infty\).
		For any \(k\in\N\) we have \(v^k=-\sign(s^k)\) if \(s^k\leq t(1+\gamma L)\), resulting in
		\[
			u^k-v^k
		{}\in{}
			\begin{ifcases}
				\tfrac{s^k}{1+\gamma L}+\sign(s^k) & s^k\leq t(1+\gamma L)
			\\[3pt]
				\tfrac{s^k}{1+\gamma\sigma}
				{}-{}
				\tfrac{\gamma(L-\sigma)t}{1+\gamma\sigma}
				{}-{}
				v^k
				\otherwise[otherwise,]
			\end{ifcases}
		\]
		where \(v^k\) is either \(1\) or \(-1\) in the second case.
		Since \(u^k-v^k\to 0\), then
		\[
			\min\set{
				\bigl|
					\tfrac{s^k}{1+\gamma L}
					{}+{}
					\sign(s^k)
				\bigr|,\,
				\bigl|
					\tfrac{s^k}{1+\gamma\sigma}
					{}-{}
					\tfrac{L-\sigma}{1+\gamma\sigma}
					\gamma t
					{}-{}
					1
				\bigr|,\,
				\bigl|
					\tfrac{s^k}{1+\gamma\sigma}
					{}-{}
					\tfrac{L-\sigma}{1+\gamma\sigma}
					\gamma t
					{}+{}
					1
				\bigr|
			}
		{}\to{}
			0.
		\]
		Notice that the first element in the set above is always larger than \(1\), and therefore eventually \(s^k\) will be always close to either
		\(
			(L-\sigma)
			\gamma t
			{}+{}
			(1+\gamma\sigma)
		\)
		or
		\(
			(L-\sigma)
			\gamma t
			{}-{}
			(1+\gamma\sigma)
		\),
		both of which are strictly smaller than \(t(1+\gamma L)\) (since \(t>1\)).
		Therefore, eventually \(s^k\leq t(1+\gamma L)\) and the residual will then be
		\(
			u^k-v^k
		{}={}
			\tfrac{s^k}{1+\gamma L}
			{}+{}
			\sign(s^k)
		\)
		which is bounded away from zero, contradicting the fact that \(u^k-v^k\to0\).
	\end{proof}
\end{thm}
\begin{thm}[Necessity of \(0<\lambda<2(1+\gamma\sigma)\)]\label{thm:lambda}%
	For any \(L>0\) and \(\sigma\in[-L,L]\) there exist \(\func{\sDRS,\nsDRS}{\R^p}{\Rinf}\) satisfying the following properties
	\begin{enumeratprop}
	\item\label{thm:lambda:s}
		\(\sDRS\) is \(L\)-smooth and \(\sigma\)-hypoconvex;
	\item\label{thm:lambda:ns}
		\(\nsDRS\) is proper, lsc, and strongly convex;
	\item\label{thm:lambda:argmin}
		\(\argmin(\sDRS+\nsDRS)\neq\emptyset\);
	\item
		for all \(s^0\), \(0<\gamma<\nicefrac1L\), and \(\lambda>2(1+\gamma\sigma)\), the sequence \(\seq{s^k}\) generated by \ref{DRS} with stepsize \(\gamma\) and relaxation \(\lambda\) starting from \(s^0\) satisfies \(\|s^k-s^{k+1}\|\not\to0\) as \(k\to\infty\) (unless \(s^0\) is a fixed point for \ref{DRS}).
	\end{enumeratprop}
	\begin{proof}
		Let \(\varphi=\sDRS+\nsDRS\), where \(\sDRS\) is as in \eqref{eq:phi1_es} with \(t=1\), and \(\nsDRS=\indicator_{\set p}\) for some \(p>1\).
		Clearly, properties \ref{thm:lambda:s}, \ref{thm:lambda:ns}, and \ref{thm:lambda:argmin} are satisfied.
		Let \(\gamma<\nicefrac 1L\), \(\lambda\geq2(1+\gamma\sigma)\).
		Starting from \(s^0\neq(1+\gamma\sigma)p+\gamma(L-\sigma)\) (so that \(u^0\neq p\)), consider \ref{DRS} with stepsize \(\gamma\) and relaxation \(\lambda\).
		To arrive to a contradiction, suppose that the residual vanishes.
		Since \(v^k=\prox_{\gamma\nsDRS}(2u^k-s^k)=p\), necessarily \(u^k\to p\); therefore, eventually \(u^k>1\) and in particular
		\[
			u^{k+1}
			{}+{}
			\gamma\tfrac{L-\sigma}{1+\gamma\sigma}
		{}={}
			\tfrac{1}{1+\gamma\sigma}
			s^{k+1}
		{}={}
			\tfrac{1}{1+\gamma\sigma}
			(s^k+\lambda(p-u^k))
		{}={}
			u^k
			{}+{}
			\gamma\tfrac{L-\sigma}{1+\gamma\sigma}
			{}+{}
			\tfrac{\lambda}{1+\gamma\sigma}(p-u^k),
		\]
		where the identity \(s^k=(1+\gamma\sigma)u^k+\gamma(L-\sigma)\) was used, cf. \eqref{eq:phi1_es_prox}.
		Therefore,
		\[
			\bigl|
				u^{k+1}-p
			\bigr|
		{}={}
			\bigl|
				1-\tfrac{\lambda}{1+\gamma\sigma}
			\bigr|
			\bigl|
				u^k-p
			\bigr|
		{}\geq{}
			\bigl|
				u^k-p
			\bigr|,
		\]
		where the inequality is due to the fact that \(\lambda\geq2(1+\gamma\sigma)\).
		Since \(u^0\neq p\) due to the choice of \(s^0\), apparently \(\seq{u^k}\) is bounded away from \(p\), hence the contradiction.
	\end{proof}
\end{thm}

\vspace*{0pt}
\par\noindent
Let us draw some conclusions:
\begin{itemize}[leftmargin=*]
\item
	The nonsmooth function \(\nsDRS\) is (strongly) convex in \Cref{thm:lambda}, therefore even for fully convex formulations
	the bound \(0<\lambda<2(1+\gamma\sigma_{\sDRS})\) needs be satisfied.
\item
	If \(\lambda>2\) (which is feasible only if \(\sDRS\) is strongly convex, \ie if \(\sigma_{\sDRS}>0\)), then, regardless of whether also \(\nsDRS\) is (strongly) convex or not, we obtain that \emph{the stepsize must be lower bounded as
	\(
		\gamma
	{}>{}
		\tfrac{\lambda-2}{2\sigma_{\sDRS}}
	\).}
	In the more general setting of \(\sigma\)-strongly monotone operators in Hilbert spaces, hence \(\sigma\geq0\), the similar bound
	\(
		\lambda
	{}<{}
		\min\set{
			2(1+\gamma\sigma),
			2+\gamma\sigma+\nicefrac{1}{\gamma\sigma}
		}
	\)
	has been recently established in \cite{monteiro2018complexity}.
\item
	Combined with the bound \(\gamma<\nicefrac{1}{L_{\sDRS}}\) shown in \Cref{thm:gamma}, we infer that (at least when \(\nsDRS\) is nonconvex) necessarily \(0<\lambda<2(1+\nicefrac{\sigma_{\sDRS}}{L_{\sDRS}})\) and consequently \(\lambda\in(0,4)\).
\end{itemize}
\begin{thm}[Tightness]\label{thm:DRStight}%
	Unless the generality of \Cref{ass:DRS} is sacrificed, when \(\lambda\in(0,2)\) or \(\sDRS\) is not strongly convex the bound
	\(
		\gamma
	{}<{}
		\min\set{
			\tfrac{1}{L_{\sDRS}},\,
			\tfrac{2-\lambda}{2[\sigma_{\sDRS}]_-}
		}
	\)
	is tight for ensuring convergence of \ref{DRS}.
	Similarly, PRS (\ie \ref{DRS} with \(\lambda=2\)) is ensured to converge iff \(\sDRS\) is strongly convex and
	\(
		\gamma
	{}<{}
		\nicefrac{1}{L_{\sDRS}}
	\).
\end{thm}

	\section{Alternating direction method of multipliers}
		\label{sec:ADMM}%
			While the classical interpretation of \ref{ADMM} as \ref{DRS} applied to the dual formulation is limited to convex problems, it has been recently observed that the two schemes are in fact related through a primal equivalence, when \(\lambda=1\).
A proof of this fact can be found in \cite[Rem. 3.14]{bauschke2015projection} when \(A=-B=\I\); in turn, \cite[Thm. 1]{yan2016self} shows that there is no loss of generality in limiting the analysis to this case.
Patterning the arguments of \cite{yan2016self} in the next subsection we will show that the equivalence can be further extended to any relaxation parameter \(\lambda\).
To this end, we introduce the notion of \emph{image function}, also known as \emph{epi-composition} or \emph{infimal post-composition} \cite{auslender2002asymptotic,bauschke2017convex,rockafellar2011variational}.
\begin{defin}[Image function]\label{defin:epicomp}%
	Given \(\func{h}{\R^n}{\Rinf}\) and \(C\in\R^{m\times n}\), the \DEF{image function} \(\func{\epicomp Ch}{\R^m}{[-\infty,+\infty]}\)
	is defined as
	\[
		\epicomp Ch(s)
	{}\coloneqq{}
		\inf_{x\in\R^n}\set{h(x)}[Cw=s].
	\]
\end{defin}
We now list some useful properties of the image function; the proofs are deferred to \Cref{sec:proof:ADMM}.
\begin{prop}\label{thm:lsc}
	Let \(\func{h}{\R^n}{\Rinf}\) and \(C\in\R^{p\times n}\).
	Suppose that for some \(\beta>0\) the set-valued mapping \(\ffunc{X_\beta}{\R^p}{\R^n}\), defined by
	\(
		X_\beta(s)
	{}\coloneqq{}
		\argmin_{x\in\R^n}\set{
			h(x)+\tfrac\beta2\|Cx-s\|^2
		}
	\),
	is nonempty for all \(s\in\R^p\).
	Then,
	\begin{enumerate}
	\item\label{thm:epiproper}%
		The image function \(\epicomp Ch\) is proper.
	\item\label{thm:epiexact}%
		\(
			\epicomp Ch(Cx_\beta)
		{}={}
			h(x_\beta)
		\)
		for all \(s\in\R^p\) and
		\(
			x_\beta
		{}\in{}
			X_\beta(s)
		\).
	\item\label{thm:AX}%
		\(
			\prox_{\nicefrac{\epicomp Ch}{\beta}}\supseteq CX_\beta
		\).
	\note{Can't we have the other inclusion?}%
\grayout{%
			and in particular the set-valued map \(s\mapsto CX_\beta(s)\) is locally bounded;
			\item\label{thm:penalty}%
				\(
					\tfrac\beta2\|Cx_\beta-\bar s\|^2\to 0
				\)
				as \(\beta\to\infty\), for any \(x_\beta\in X_\beta(\bar s)\);
			\item\label{thm:epilsc}%
				\(\epicomp Ch\) is lsc at \(\bar s\) iff
				\(
					h(x_\beta)\to\epicomp Ch(\bar s)
				\)
				as
				\(
					\beta\to\infty
				\),
				where \(x_\beta\in X_\beta(\bar s)\).
}%
	\end{enumerate}
\end{prop}
\begin{prop}\label{thm:epipartial}%
	For a function \(\func{h}{\R^n}{\Rinf}\) and \(C\in\R^{p\times n}\), let \(\ffunc X{\R^p}{\R^n}\) be defined as
	\(
		X(s)
	{}\coloneqq{}
		\argmin_{x\in\R^n}\set{h(x)}[Cx=s]
	\).
	Then, for all \(\bar s\in C\dom h\) and \(\bar x\in X(\bar s)\) it holds that
	\[
		\trans C\hat\partial\epicomp Ch(\bar s)
	{}\subseteq{}
		\hat\partial h(\bar x).
	\]
\end{prop}

\begin{prop}[Strong convexity of the image function]\label{thm:epicompStrCvx}%
	Suppose that \(\func h{\R^n}{\Rinf}\) is proper, lsc, and \(\sigma_h\)-strongly convex.
	Then, for every \(C\in\R^{p\times n}\) the image function \(\epicomp Ch\) is \(\sigma_{\epicomp Ch}\)-strongly convex with \(\sigma_{\epicomp Ch}=\nicefrac{\sigma_h}{\|C\|^2}\).
\end{prop}

\grayout{%
	In order to proceed to the next result, we first need to introduce the following important notion for parametric minimization.
	\begin{defin}[Locally uniform level boundedness {\cite[Def. 1.16]{rockafellar2011variational}}]\label{defin:ULB}%
		We say that a function \(\func{M}{\R^m\times\R^n}{\Rinf}\) with values \(M(w,x)\) is \DEF{level bounded in \(w\) locally uniformly in \(x\)} if for all \(\alpha\in\R\) and \(\bar x\in\R^n\) there exists \(\varepsilon>0\) such that the set
		\[
			\set{(w,x)\in\R^m\times\R^n}[
				M(w,x)\leq\alpha,~
				\|x-\bar x\|\leq\varepsilon
			]
		\]
		is bounded.
	\end{defin}

	\begin{thm}\label{thm:epicompLB}%
		Let \(\func{h}{\R^n}{\Rinf}\) be lsc and \(C\in\R^{p\times n}\).
		Suppose that for some \(\beta>0\) the function \(h+\tfrac\beta2\|C{}\cdot{}-s\|^2\) is level bounded for all \(s\in\R^p\).
		Then, the following hold:
		\begin{enumerate}
		\item
			\(\epicomp Ch\) is proper and lsc.
		\item
			For all \(s\in C\dom h\) the set of minimizers \(X(s)\) as in \eqref{eq:X(s)} is nonempty; moreover, \(X\) is locally bounded, and it is osc with respect to \(\epicomp Ch\)-attentive convergence: for all \(\bar s\in C\dom h\)
			\[
				\limsup_{k\to\infty}X(s_k)
			{}\subseteq{}
				X(\bar s)
			\]
			whenever \(\bigl(s^k,\epicomp Ch(s^k)\bigr)\to\bigl(\bar s,\epicomp Ch(\bar s)\bigr)\) as \(k\to\infty\).
		\item
			For all \(\bar s\in C\dom h\) and \(\bar x\in X(\bar s)\) one has
			\[
				\trans C\partial\epicomp Ch(\bar s)
			{}\subseteq{}
				\bigcup_{\bar x\in X(\bar s)}{
					\partial h(\bar x)
				}.
			\]
		\end{enumerate}
		\begin{proof}
			The level boundedness condition ensures that 
			\(
				H(x,s)
			{}\coloneqq{}
				h(x)+\indicator_{\set0}(Cx-s)
			\)
			is level bounded in \(x\), locally uniformly in \(s\), cf. \cref{defin:ULB}.
			The first two claims then follow from \cite[Thm. 1.32]{rockafellar2011variational}.
			
			Let \(\bar v\in\partial\epicomp Ch(\bar s)\) be fixed.
			Then, there exits a sequence \(\seq{s^k,v^k}\subseteq\graph\hat\partial h\) such that
			\(
				\bigl(
					s^k,
					\epicomp Ch(s^k),
					v^k
				\bigr)
			{}\to{}
				\bigl(
					\bar s,
					\epicomp Ch(\bar s),
					\bar v
				\bigr)
			\)
			as \(k\to\infty\).
			For each \(k\in\N\) let \(x^k\in X(s^k)\); then, \(\seq{x^k}\) is bounded and all its accumulation points belong to \(X(\bar s)\); thus, up to possibly extracting, \(x^k\to\bar x\) for some \(\bar x\in X(\bar s)\) as \(k\to\infty\).
			Then,
			\[
				\trans C\bar v
			{}={}
				\lim_{k\to\infty}{
					\trans Cv^k
				}
			{}\overrel*[\in]{\ref{thm:epipartial}}[2pt]{}
				\limsup_{k\to\infty}{
					\hat\partial h(x^k)
				}
			{}\subseteq{}
				\partial h(\bar x),
			\]
			where the last inclusion follows from the definition of \(\partial h\) and the fact that
			\[
				h(x^k)
			{}={}
				\epicomp Ch(s^k)
			{}\to{}
				\epicomp Ch(\bar s)
			{}={}
				h(\bar x).
			\]
			The claimed inclusion then follows from the arbitrarity of \(\bar v\in\partial\epicomp Ch(\bar s)\).
		\end{proof}
	\end{thm}

	\begin{lem}\label{thm:epicomp:convex}%
		Let \(\func h{\R^n}{\Rinf}\) be convex and \(C\in\R^{p\times n}\) be surjective.
		Then, \(\epicomp Ch\) is convex, and as long as the set of minimizers \(X(\bar s)\) is nonempty (see \eqref{eq:X(s)}), it holds that
		\[
			\partial\epicomp Ch(\bar s)
		{}={}
			\set{y}[\trans Ay\in\partial h(\bar x)],
		\]
		where \(\bar x\) is any element of \(X(\bar s)\).
		In particular, if \(h\) is differentiable at some point in \(X(\bar s)\), then \(\epicomp Ch\) is differentiable at \(\bar s\).
		\begin{proof}
			See \cite[Thm. D.4.5.1 and Cor. D.4.5.2]{hiriarturruty2012fundamentals}.
		\end{proof}
	\end{lem}
}%

		\subsection{A universal equivalence of DRS and ADMM}
			Let us eliminate the linear coupling between \(x\) and \(z\) in the \ref{ADMM} problem formulation \eqref{eq:CP}, so as to bring it into \ref{DRS} form \eqref{eq:P}.
To this end, let us introduce a slack variable \(s\in\R^p\) and rewrite \eqref{eq:CP} as
\begin{align*}
	\minimize_{x\in\R^m,z\in\R^n,s\in\R^p}{}
&
	f(x)+g(z)
\quad
	\stt Ax=s,~Bz=b-s.
\intertext{%
	Since the problem is independent of the order of minimization \cite[Prop. 1.35]{rockafellar2011variational}, we may minimize first with respect to \((x,z)\) to arrive to
}
	\minimize_{s\in\R^p}{}
&
	\inf_{x\in\R^m}\set{f(x)}[Ax=s]
	{}+{}
	\inf_{z\in\R^n}\set{g(z)}[Bz=b-s].
\intertext{%
	The two parametric infima define two image functions, cf. \Cref{defin:epicomp}: indeed, \ref{ADMM} problem formulation \eqref{eq:CP} can be expressed as
}
\numberthis\label{eq:CP2P}
	\minimize_{s\in\R^p}{}
&
	\epicomp Af(s)
	{}+{}
	\epicomp Bg(b-s),
\end{align*}
which is exactly \eqref{eq:P} with \(\sDRS=\epicomp Af\) and \(\nsDRS=\epicomp Bg(b-{}\cdot{})\).
Apparently, unless \(A\) and \(B\) are injective the correspondence between variable \(s\) in \eqref{eq:CP2P} and variables \(x,z\) in \eqref{eq:CP} may fail to be one to one, as \(s\) is associated to sets of variables \(x\in X(s)\) and \(z\in Z(s)\) defined as
\[
	X(s)
{}\coloneqq{}
	\argmin_{x\in\R^m}\set{f(x)}[Ax=s]
\quad\text{and}\quad
	Z(s)
{}\coloneqq{}
	\argmin_{z\in\R^n}\set{g(z)}[Bz=b-s].
\]

\begin{thm}[Primal equivalence of \ref{DRS} and \ref{ADMM}]\label{thm:ADMM_DRS}%
	Starting from a triplet \((x,y,z)\in\R^m\times\R^p\times\R^n\), consider an \ref{ADMM}-update applied to problem \eqref{eq:CP} with relaxation \(\lambda\) and large enough penalty \(\beta>0\) so that any \ref{ADMM} minimization subproblem has solutions.
	Let
	\begin{equation}\label{eq:ADMM2DRS}
		\begin{cases}[l >{{}}c<{{}} l]
			s
		&\coloneqq&
			Ax-\nicefrac y\beta
		\\
			u
		&\coloneqq&
			Ax
		\\
			v
		&\coloneqq&
			b-Bz
		\end{cases}
	\quad\text{and, similarly,}\quad
		\begin{cases}[l >{{}}c<{{}} l]
			s^+
		&\coloneqq&
			Ax^+-\nicefrac{y^+}\beta
		\\
			u^+
		&\coloneqq&
			Ax^+
		\\
			v^+
		&\coloneqq&
			b-Bz^+.
		\end{cases}
	\end{equation}
	Then, the variables are related as follows:
	\[
		\begin{cases}[l >{{}}c<{{}} l]
			s^+
		&=&
			s+\lambda(v-u)
		\\
			u^+
		&\in&
			\prox_{\gamma\sDRS}(s^+)
		\\
			v^+
		&\in&
			\prox_{\gamma\nsDRS}(2u^+-s^+),
		\end{cases}
	\quad\text{where}\quad
		\begin{cases}[l >{{}}c<{{}} l]
			\sDRS
		&\coloneqq&
			\epicomp Af
		\\
			\nsDRS
		&\coloneqq&
			\epicomp Bg(b-{}\cdot{})
		\\
		\gamma
		&\coloneqq&
			\nicefrac1\beta.
		\end{cases}
	\]
	Moreover,
	\begin{enumerate}
	\item\label{thm:f(x)}
		\(\sDRS(u^+)=\epicomp Af(Ax^+)=f(x^+)\),
	\item\label{thm:ADMM:g(z)}%
		\(\nsDRS(v^+)=\epicomp Bg(Bz^+)=g(z^+)\),
	\item\label{thm:ADMM:-y}%
		\(-y^+\in\hat\partial\sDRS(u^+)=\hat\partial\epicomp Af(Ax^+)\),
	\item\label{thm:ADMM:-Ay}%
		\(
			-\trans Ay^+\in\hat\partial f(x^+)
		\), and
	\item\label{thm:ADMM:-By}%
		\(
			\dist(-\trans By^+,\hat\partial g(z^+))
		{}\leq{}
			\beta\|B\|
			\|Ax^++Bz^+-b\|
		\).
	\end{enumerate}
	If, additionally, \(A\) has full row rank, \(\sDRS\in C^{1,1}(\R^p)\) is \(L_{\sDRS}\)-smooth, and \(\beta>L_{\sDRS}\), then it also holds that
	\begin{enumerate}[resume]
	\item\label{thm:DRE2LL}
		\(\DRE(s^+)=\LL_\beta(x^+,z^+,y^+)\).
	\end{enumerate}
	\begin{proof}
		\begin{subequations}
			Observe first that, as shown in \cref{thm:AX}, it holds that
			\begin{align}
				\label{eq:proxs}
				\prox_{\gamma\sDRS}
			{}\supseteq{} &
				A\argmin\set{
					f+\tfrac{1}{2\gamma}\|A{}\cdot{}-s\|^2
				}.
			\shortintertext{%
				Similarly, with a simple change of variable one obtains that
			}
				\label{eq:proxns}
				\prox_{\gamma\nsDRS}
			{}\supseteq{} &
				b-B\argmin\set{
					g+\tfrac{1}{2\gamma}\|B{}\cdot{}+s-b\|^2
				}.
			\end{align}
		\end{subequations}
		Let \((s,u,v)\) and \((s^+,u^+,v^+)\) be as in \eqref{eq:ADMM2DRS}.
		We have
		\[
			s+\lambda(v-u)
		{}={}
			Ax-\tfrac1\beta y-\lambda(Ax+Bz-b)
		{}={}
			Ax-\tfrac1\beta y^{\nicefrac+2}-(Ax+Bz-b)
		{}={}
			{-{}}\tfrac1\beta y^++Ax^+
		{}={}
			s^+,
		\]
		where in the second and third equality the \ref{ADMM} update rule for \(y^{\nicefrac+2}\) and \(y^+\), respectively, was used.
		Moreover,
		\[
			u^+
		{}={}
			Ax^+
		{}\in{}
			A\argmin\LL_\beta({}\cdot{},z,y^{\nicefrac+2})
		{}\overrel*[\subseteq]{\eqref{eq:proxs}}{}
			\prox_{\nicefrac{\sDRS}{\beta}}(b-Bz-y^{\nicefrac+2}\nicefrac{}\beta)
		{}={}
			\prox_{\nicefrac{\sDRS}{\beta}}(s^+),
		\]
		where the last equality uses the identity
		\(
			b-Bz-y^{\nicefrac+2}\nicefrac{}\beta
		{}={}
			v-\gamma y+(1-\lambda)(u-v)
		{}={}
			s+\lambda(v-u)
		{}={}
			s^+
		\).
		Next, observe that
		\(
			2u^+-s^+
		{}={}
			2Ax^+-(Ax^+-y^+\nicefrac{}\beta)
		{}={}
			Ax^++y^+\nicefrac{}\beta
		\),
		hence
		\[
			v^+
		{}={}
			b-Bz^+
		{}\in{}
			b-B\argmin\LL_\beta(x^+,{}\cdot{},y^+)
		{}\overrel[\subseteq]{\eqref{eq:proxns}}{}
			\prox_{\nicefrac{\nsDRS}{\beta}}(Ax^++y^+\nicefrac{}\beta)
		{}={}
			\prox_{\nicefrac{\nsDRS}{\beta}}(2u^+-s^+).
		\]
		
		Let us now show the numbered claims.
		\begin{proofitemize}
		\item\ref{thm:f(x)} \& \ref{thm:ADMM:g(z)}.
			Follow from \cref{thm:epiexact}.
		\item\ref{thm:ADMM:-y}.
			Since \(u^+\in\prox_{\gamma\sDRS}(s^+)\) and
			\(
				-y^+
			{}={}
				\tfrac1\gamma(s^+-u^+)
			\),
			the claim follows from \eqref{eq:prox:subdiff}.
		\item\ref{thm:ADMM:-Ay}.
			This follows from the optimality conditions of \(x^+\) in the \ref{ADMM}-subproblem defining the \(x\)-update.
			Alternatively, the claim can also be deduced from \ref{thm:ADMM:-y} and \cref{thm:epipartial}.
		\item\ref{thm:ADMM:-By}.
			The optimality conditions in the \ref{ADMM}-subproblem defining the \(z\)-update read
			\[
				0
			{}\in{}
				\hat\partial_z\LL_\beta(x^{k+1},z^{k+1},y^{k+1})
			{}={}
				\trans B(Ax^{k+1}+Bz^{k+1}-b+y^{k+1}\nicefrac{}\beta),
			\]
			and the claim readily follows.
		\item\ref{thm:DRE2LL}.
			Suppose now that \(\sDRS\) is \(L_{\sDRS}\)-smooth (hence \(A\) is surjective, for otherwise \(\sDRS\) has not full domain), and that \(\beta>L_{\sDRS}\).
			Due to smoothness, the inclusion in \ref{thm:ADMM:-y} can be strengthened to
			\(
				\nabla\sDRS(u^+)
			{}={}
				-y^+
			\).
			We may then invoke the expression \eqref{eq:DRE} of the DRE (recall that the minimum is attained at \(v^+\)) to obtain
			\begin{qedalign*}
				\DRE(s^+)
			{}={} &
				\sDRS(u^+)
				{}+{}
				\nsDRS(v^+)
				{}+{}
				\innprod{\nabla\sDRS(u^+)}{v^+-u^+}
				{}+{}
				\tfrac{1}{2\gamma}\|v^+-u^+\|^2
			\\
			{}={} &
				f(x^+)
				{}+{}
				g(z^+)
				{}+{}
				\innprod{y^+}{Ax^++Bz^+-b}
				{}+{}
				\tfrac\beta2\|Ax^++Bz^+-b\|^2
			{}={}
				\LL_\beta(x^+,z^+,y^+).
			\end{qedalign*}
		\end{proofitemize}\let\qed\relax
	\end{proof}
\end{thm}

		\subsection{Convergence of the ADMM}
			In order to extend the theory developed for \ref{DRS} to \ref{ADMM} we shall impose that \(\sDRS\) and \(\nsDRS\) as in \eqref{eq:CP2P} comply with \Cref{ass:DRS}.
This motivates the following blanket requirement.
\begin{ass}[Requirements for the ADMM formulation \eqref{eq:CP}]\label{ass:ADMM}%
	The following hold:
	\begin{enumeratass}
	\item\label{ass:ADMM:lsc}%
		\(\func{f}{\R^m}{\R}\) and \(\func{g}{\R^n}{\Rinf}\) are proper and lsc.
	\item\label{ass:ADMM:beta}%
		\(A\) is surjective, and \(\beta\) is large enough so that the \ref{ADMM} subproblems have solution.%
		\note{\rm If we prove equality in \cref{thm:AX} the requirement on \(\beta\) can be discarded}%
	\item\label{ass:ADMM:f}%
		\(\sDRS\coloneqq\epicomp Af\in\cont^{1,1}(\R^p)\) is \(L_{\epicomp Af}\)-smooth, hence \(\sigma_{\epicomp Af}\)-hypoconvex with \(|\sigma_{\epicomp Af}|\leq L_{\epicomp Af}\).
	\item\label{ass:ADMM:g}%
		\(\nsDRS\coloneqq\epicomp Bg\) is lsc.
	\item\label{ass:ADMM:LB}%
		Problem \eqref{eq:CP} has a solution: \(\argmin\Phi\neq\emptyset\), where \(\Phi(x,z)\coloneqq f(x)+g(z)+\indicator_S(x,z)\) and
		\(
			S\coloneqq\set{(x,z)\in\R^m\times\R^n}[Ax+Bz=b]
		\)
		is the feasible set.
	\end{enumeratass}
\end{ass}
These requirements generalize those in \Cref{ass:DRS} by allowing linear constraints more generic than \(x-z=0\), cf. \eqref{eq:DRS-constrained}.
The assumption of surjectivity of \(A\) is as general as requiring the inclusion \(\range B\subseteq b+\range A\).
In fact, (up to an orthogonal transformation) without loss of generality we may assume that \(A=\binom{A'}{~}\) for some surjective matrix \(A'\in\R^{r\times m}\), where \(r=\rank A\), stacked over a \((p-r)\times n\) zero matrix.
Then, in light of the prescribed range inclusion necessarily \(B=\binom{B'}{~}\) and \(b=\binom{b'}{~}\), for some \(B'\in\R^{r\times n}\) and \(b\in\R^r\).
Then, problem \eqref{eq:CP} can be simplified to the minimization of
\(
		f(x)+g(z)
\)
subject to
\(
	A'x+B'z=b'
\),
which satisfies the needed surjectivity property.
\begin{thm}[Convergence of \ref{ADMM}]\label{thm:ADMM}%
	Suppose that \Cref{ass:ADMM} is satisfied, and let \(\sDRS\), \(\nsDRS\), and \(\Phi\) be as defined therein.
	Starting from \((x^{-1},y^{-1},z^{-1})\in\R^m\times\R^p\times\R^n\), consider a sequence \(\seq{x^k,y^k,z^k}\) generated by \ref{ADMM} with penalty \(\beta=\nicefrac1\gamma\) and relaxation \(\lambda\), where \(\gamma\) and \(\lambda\) are as in \Cref{thm:DRE:SD}.
	The following hold:
	\begin{enumerate}
	\item\label{thm:ADMM:r}
		\(
			\LL_\beta(x^{k+1},z^{k+1},y^{k+1})
		{}\leq{}
			\LL_\beta(x^k,z^k,y^k)
			{}-{}
			\frac{c\lambda^2}{(1+\gamma L_{\epicomp Af})^2}
			\|Ax^k+Bz^k-b\|^2
		\),
		where \(c\) is as in \Cref{thm:DRE:SD}, and the residual \(\seq{Ax^k+Bz^k-b}\) vanishes with
		\(
			\min_{i\leq k}\|Ax^i+Bz^i-b\|
		{}={}
			o(\nicefrac{1}{\sqrt k})
		\).
	\item\label{thm:ADMM:uv}
		all cluster points \((x,z,y)\) of \(\seq{x^k,z^k,y^k}\) satisfy the KKT conditions
		\begin{itemize}
		\item
			\(-\trans Ay\in\partial f(x)\)
		\item
			\(-\trans By\in\partial g(z)\)%
			\note{\rm Can it be extended to \(\hat\partial f\) and \(\hat\partial g\)?}
		\item
			\(Ax+Bz=b\),
		\end{itemize}
		and attain the same cost \(f(x)+g(z)\), this being the limit of \(\seq{\LL_\beta(x^k,z^k,y^k)}\).
	\item\label{thm:ADMM:bounded}
		the sequence \(\seq{Ax^k,y^k,Bz^k}\) is bounded provided that the cost function \(\Phi\) is level bounded.
		If, additionally, \(f\in\cont^{1,1}(\R^m)\), then the sequence \(\seq{x^k,y^k,z^k}\) is bounded.
	\end{enumerate}
	\begin{proof}
		Let \(s^0\coloneqq Ax^0-y^0\nicefrac{}\beta\), and consider the sequence \(\seq{s^k,u^k,v^k}\) generated by \ref{DRS} applied to \eqref{eq:CP2P}, with stepsize \(\gamma\), relaxation \(\lambda\), and starting from \(s^0\).
		Then, for all \(k\in\N\) it follows from \cref{thm:ADMM_DRS} that the variables are related as
		\[
			\begin{cases}[l >{{}}c<{{}} l]
				s^k
			&=&
				Ax^k-y^k\nicefrac{}\beta
			\\
				u^k
			&=&
				Ax^k
			\\
				v^k
			&=&
				b-Bz^k,
			\end{cases}
		\]
		and satisfy
		\[
			\begin{cases}[l >{{}}c<{{}} l]
				\sDRS(u^k)
			&=&
				f(x^k)
			\\
				\nsDRS(v^k)
			&=&
				g(z^k)
			\\
				\DRE(s^k)
			&=&
				\LL_\beta(x^k,z^k,y^k)
			\end{cases}
		\qquad\text{and}\qquad
			\begin{cases}
				y^k=-\nabla\sDRS(u^k)
			\\
				-\trans Ay^k\in\hat\partial f(x^k)
			\\
				\dist(-\trans By^k,\hat\partial g(z^k))\to0.
			\end{cases}
		\]
		\begin{proofitemize}
		\item\ref{thm:ADMM:r}.
			Readily follows from \cref{thm:DRS}.
		\item\ref{thm:ADMM:uv}.
			Suppose that for some \(K\subseteq\N\) the subsequence \(\seq{x^k,y^k,z^k}[k\in K]\) converges to \((x,y,z)\); then, necessarily \(Ax+Bz=b\).
			Moreover,
			\[
				\epicomp Af(Ax)
			{}\leq{}
				f(x)
			{}\leq{}
				\liminf_{K\ni k\to\infty}{
					f(x^k)
				}
			{}={}
				\liminf_{K\ni k\to\infty}{
					\epicomp Af(Ax^k)
				}
			{}={}
				\epicomp Af(Ax),
			\]
			where the second inequality is due to the fact that \(f\) is lsc, and the last one to the fact that \(\epicomp Af\) is continuous.
			Therefore, \(f(x^k)\to f(x)\), and the inclusion \(-\trans Ay^k\in\hat\partial f(x^k)\) in light of the definition of subdifferential results in \(-\trans Ay\in\partial f(x)\).
			In turn, since \(\sDRS(u^k)+\sDRS(v^k)\) converges to \(\sDRS(Ax)+\nsDRS(b-Bz)=\epicomp Af(Ax)+\epicomp Bg(Bz)\) as it follows from \cref{thm:DRS:uv}, a similar reasoning shows that \(g(z^k)\to g(z)\) as \(K\ni k\to\infty\).
			Thus, since \(\dist(-\trans By^k,\hat\partial g(z^k))\to0\), \(g\)-attentive outer semicontinuity of \(\partial g\), see \cite[Prop. 8.7]{rockafellar2011variational}, implies that \(-\trans By\in\partial g(z)\).
			Finally, that \(f(x)+g(z)\) equals the limit of the whole sequence \(\seq{\LL_\beta(x^k,z^k,y^k)}\) then follows from \cref{thm:DRS:uv} through the identity
			\(
				\DRE(s^k)
			{}={}
				\LL_\beta(x^k,z^k,y^k)
			\).
		\item\ref{thm:ADMM:bounded}.
			Once we show that \(\varphi=\sDRS+\nsDRS\) is level bounded, boundedness of the sequence \(\seq{Ax^k,Bz^k,y^k}\) will follow from \cref{thm:DRS:bounded}.
			For \(\alpha\in\R\) we have
			\begin{align*}
				\lev_{\leq\alpha}\varphi
			{}={} &
				\set{s}[{
					\inf_x\set{f(x)}[
						Ax=s
					]
					{}+{}
					\inf_z\set{g(z)}[
						Bz=b-s
					]
					{}\leq{}
					\alpha
				}]
			\\
			{}={} &
				\set{s}[{
					\inf_{x,z}\set{
						f(x)+g(z)
					}[
						Ax=s,\,
						Bz=b-s
					]
					{}\leq{}
					\alpha
				}]
			\\
			{}={} &
				\set{Ax}[
					f(x)+g(z)\leq\alpha,~
					\exists z:Ax+Bz=b
				]
			{}={}
				\set{Ax}[{
					(x,z)\in\lev_{\leq\alpha}\Phi,\,
					\exists z
				}].
			\end{align*}
			Since
			\(
				\|Bz\|
			{}\leq{}
				\|B\|\|z\|
			{}\leq{}
				\|B\|\|(x,z)\|
			\)
			for any \(x,z\), it follows that if \(\lev_{\leq\alpha}\Phi\) is bounded, then so is \(\lev_{\leq\alpha}\varphi\).
			Suppose now that \(f\in\cont^{1,1}(\R^n)\) is \(L_f\)-smooth, and for all \(k\in\N\) let
			\(
				\xi^k
			{}\coloneqq{}
				x^k
				{}-{}
				\trans A(A\trans A)^{-1}
				(Ax^k+Bz^k-b)
			\).
			Then, \(A\xi^k=b-Bz^k\), hence \(f(\xi^k)+g(z^k)=\Phi(\xi^k,z^k)\), and \(\xi^k-x^k\to 0\) as \(k\to\infty\).
			We have
			\begin{align*}
				|
					\Phi(\xi^k,z^k)-(f(x^k)+g(z^k))
				|
			{}={} &
				\bigl|
					f(\xi^k)-f(x^k)
				\bigr|
			{}\leq{}
				|
					\innprod{\nabla f(x^k)}{\xi^k-x^k}
				|
				{}+{}
				\tfrac{L_f}{2}
				\|\xi^k-x^k\|^2
			\\
			{}\leq{} &
				\bigl|
					\innprod{y^k}{Ax^k-A\xi^k}
				\bigr|
				{}+{}
				\tfrac{L_f}{2}
				\|\trans A(A\trans A)^{-1}\|^2
				\|Ax^k+Bz^k-b\|^2,
			\end{align*}
			where in the second inequality the identity \(\nabla f(x^k)=-\trans Ay^k\) was used, cf. \cref{thm:ADMM:-Ay}.
			In particular, \(f(\xi^k)-f(x^k)\to0\) as \(k\to\infty\), and therefore \(\Phi(\xi^k,z^k)\) converges to a finite quantity (the limit of \(\LL_\beta(x^k,z^k,y^k)\)).
			Since \(\Phi\) is level bounded, necessarily \(\seq{\xi^k,z^k}\) is bounded, hence so is \(\seq{x^k}\).
		\qedhere
		\end{proofitemize}
	\end{proof}
\end{thm}
The smoothness condition on \(f\) required in \Cref{thm:ADMM:bounded} is a standing assumption in the (proximal) ADMM analysis of \cite{li2015global}, which, together with the restriction \(A=\I\), ensures that \(\epicomp Af=f\) complies with \Cref{ass:ADMM:f}.
Our requirement of level boundedness of \(\Phi\) to ensure boundedness of the sequences generated by \ref{ADMM} is milder than that of \cite[Thm. 3]{li2015global}, which instead requires coercivity of either \(f\) or \(g\).

\begin{rem}[Simpler bounds for \ref{ADMM}]\label{rem:ADMM}%
	In parallel with the simplifications outlined in \Cref{rem:DRS} for \ref{DRS}, simpler (more conservative) bounds for the penalty parameter \(\beta\) in \ref{ADMM} are, in case \(\lambda\in(0,2]\):
	\begin{align*}
		\lambda\in(0,2)
	~&
		\begin{cases}[l @{~\text{ and }~} l @{~~~} l]
			\beta>L
		&
			c
			{}={}
			\beta\frac{2-\lambda}{2\lambda}
			{}-{}
			L[\nicefrac12-\nicefrac L\beta]_+
		&
			\text{if \(f\) is convex,}
		\\[5pt]
			\beta>\frac{2L}{2-\lambda}
		&
			c
			{}={}
			\beta\frac{2-\lambda}{2\lambda}
			{}-{}
			\frac L\lambda
		&
			\text{otherwise,}
		\end{cases}
	\\[5pt]
		\lambda=2
	~&
		\begin{cases}[l c l @{~~~} l]
			\beta>L
		&~\text{ and }~~&
			c
			{}={}
			\frac{\sigma_f}{2\|A\|^2}(1-\nicefrac L\beta)
		&
			\text{if \(f\) is strongly convex,}
		\\[5pt]
			\hphantom{
				\beta>\frac{2L}{2-\lambda}
			}
		&\emptyset&
			\hphantom{
				c
				{}={}
				\beta\frac{2-\lambda}{2\lambda}
				{}-{}
				L[\nicefrac12-\nicefrac L\beta]_+
			}
		&
			\text{otherwise,}
		\end{cases}
	\end{align*}
	where
	\(
		L
	{}\coloneqq{}
		L_{\epicomp Af}
	\).
	The case \(\lambda=2\) uses \Cref{thm:epicompStrCvx} to infer strong convexity of \(\epicomp Af\) from that of \(f\).
\end{rem}

As a consequence of the Tarski-Seidenberg theorem, functions \(\sDRS\coloneqq\epicomp Af\) and \(\nsDRS\coloneqq\epicomp Bg(b-{}\cdot{})\) are semialgebraic provided \(f\) and \(g\) are, see \eg \cite{bochnak2013real}.
Therefore, sufficient conditions for global convergence of \ref{ADMM} follow from the similar result for \ref{DRS} stated in \Cref{thm:DRSglobal}, through the primal equivalence of the algorithms illustrated in \Cref{thm:ADMM_DRS}.
We should emphasize, however, that the equivalence identifies \(u^k=Bz^k\) and \(v^k=b-Ax^k\); therefore, only convergence of \(\seq{Ax^k,y^k,Bz^k}\) can be deduced, as opposed to that of \(\seq{x^k,y^k,z^k}\).
\begin{thm}[Global convergence of \ref{ADMM}]%
	Suppose \Cref{ass:ADMM} is satisfied, and let \(\Phi\) be as defined therein.
	If \(\Phi\) is level bounded and \(f\) and \(g\) are semialgebraic, then the sequence \(\seq{Ax^k,y^k,Bz^k}\) generated by \ref{ADMM} with \(\beta\) and \(\lambda\) as in \Cref{thm:ADMM} converges.
\end{thm}

		\subsection{Adaptive variant}
			Similar to what done for \ref{DRS}, one can still ensure a sufficient decrease property on the augmented Lagrangian without knowing the exact value of \(L_{\epicomp Af}\), when \(\lambda\in(0,2)\).
However, due to the implicitness of \(\sDRS=\epicomp Af\), enforcing the inequality \(\varphi(v^k)\leq\mathcal L_k\) as in \cref{step:DRSadaptive:if} of \Cref{alg:DRSadaptive}, needed to ensure the lower boundedness of \(\seq{\LL_{\nicefrac1\gamma}(x^k,z^k,y^k)}\), may not be possible.
Indeed, although we may exploit \eqref{eq:ADMM2DRS} and \Cref{thm:ADMM:g(z)} to arrive to
\[
	\varphi(v^k)
{}={}
	\sDRS(v^k)
	{}+{}
	\nsDRS(v^k)
{}={}
	\epicomp Af(b-Bz^k)
	{}+{}
	g(z^k),
\]
the value of \(\epicomp Af(b-Bz^k)\) may not be readily available.
In the following special cases, however, one can bypass the problem.
\begin{proofitemize}
\item\emph{\(A\) is square and with known inverse \(A^{-1}\):}
	then, \(\epicomp Af(b-Bz^k)=f(A^{-1}(b-Bz^k))\).
\item\emph{A constant \(\Phi_{\text{\sc lb}}\leq\inf\Phi\) is known:}
	in this case, one can rather enforce \(\Phi_{\text{\sc lb}}\leq\mathcal L_k\).
\end{proofitemize}

This detail apart, the adaptive variant of \ref{DRS} outlined in \Cref{alg:DRSadaptive} can be easily translated into an adaptive version of \ref{ADMM} in which the penalty \(\beta\) is suitably adjusted.
For the sake of simplicity, we only consider the case \(\lambda=1\), so that the half-update \(y^{\nicefrac+2}\) can be discarded.

\begin{algorithm}[tb]
	\algcaption[\ref{ADMM} with adaptive stepsize]{%
		\ref{ADMM} with adaptive stepsize (\(\lambda=1\) for simplicity).\\
		\(\ffunc{\text{\ref{ADMM}}_\beta}{\R^p\times\R^n}{\R^m\times\R^p\times\R^n}\) maps \((y,z)\) to a triplet \((x^+,y^+,z^+)\) as in \eqref{ADMM} with \(\lambda=1\) (since \(\lambda=1\), the update does not depend on \(x\)).%
	}%
	\label{alg:ADMMadaptive}%
	\begin{algorithmic}[1]
	\setlength\baselineskip{1.25\baselineskip}
	\Require
		\((y^{-1},z^{-1})\in\R^p\times\R^n\),~
		\(L>0\),~
		\(\beta,c\) as in \cref{rem:ADMM} with \(\lambda=1\)%
	\Initialize
		\(
			(x^0,y^0,z^0)
		{}\in{}
			{\rm\ref{ADMM}}_\beta(y^{-1},z^{-1})
		\),~
		\(
			\mathcal L_0
		{}={}
			\LL_\beta(x^0,z^0,y^0)
		\)
	\item[{\bf For} ~\(k=0,1,\ldots\)~ {\bf do}]
	\State\label{step:ADMMadaptive:1}%
		\(
			(x^{k+1},y^{k+1},z^{k+1})
		{}\in{}
			{\rm\ref{ADMM}}_\beta(y^k,z^k)
		\)
	\Statex
		\(
			\mathcal L_{k+1}
		{}={}
			\LL_\beta(x^{k+1},y^{k+1},z^{k+1})
		\)
	\If{~
		\(
			\mathcal L_{k+1}
		{}>{}
			\mathcal L_k
			{}-{}
			\frac{c\lambda^2}{(1+\nicefrac L\beta)^2}
			\|Ax^k+Bz^k-b\|^2
		\)
	~}\label{step:ADMMadaptive:if}%
		\State%
			\(\beta\gets2\beta\),~
			\(c\gets2c\),~
			\(L\gets2L\)
		\Statex\hspace*{\algorithmicindent}%
			\(
				(x^k,y^k,z^k)
			{}\in{}
				{\rm\ref{ADMM}}_\beta(y^{k-1},z^{k-1})
			\)
		\Statex\hspace*{\algorithmicindent}%
			\(
				\mathcal L_k
			{}\gets{}
				\LL_\beta(x^k,y^k,z^k)
			\)
			and go back to \cref{step:ADMMadaptive:1}
	\EndIf{}
\end{algorithmic}
\end{algorithm}

\begin{thm}[Subsequential convergence of adaptive \ref{ADMM}]
	Suppose that \Cref{ass:ADMM} is satisfied, and consider the iterates generated by \Cref{alg:ADMMadaptive}.
	If the sequence \(\seq{\mathcal L_k}\) is lower bounded, then the following hold:
	\begin{enumerate}
	\item\label{thm:ADMMadaptive:uv}
		All cluster points \((x,y,z)\) of \(\seq{x^k,y^k,z^k}\) satisfy the KKT conditions
		\begin{itemize}
		\item
			\(-\trans Ay\in\partial f(x)\)
		\item
			\(-\trans By\in\partial g(z)\)%
		\item
			\(Ax+Bz=b\),
		\end{itemize}
		and attain the same cost \(f(x)+g(z)\), this being the limit of \(\seq{\mathcal L_k}\).
	\item\label{thm:ADMMadaptive:r}
		The residual \(\seq{\|Ax^k+Bz^k-b\|}\) vanishes with rate
		\(
			\min_{i\leq k}{
				\|Ax^i+Bz^i-b\|
			}
		{}\leq{}
			o(\nicefrac{1}{\sqrt k})
		\).
	\end{enumerate}
	In particular, the claims hold if at some iteration the inequality \(\beta>L_{\epicomp Af}\) is satisfied.
	In this case, and if the cost function \(\Phi\) is level bounded, the following also hold:
	\begin{enumerate}[resume]
	\item
		the sequence \(\seq{Ax^k,y^k,Bz^k}\) is bounded.
	\item
		the sequence \(\seq{Ax^k,y^k,Bz^k}\) is convergent if \(f\) and \(g\) are semialgebraic.
	\end{enumerate}
\end{thm}

		\subsection{Sufficient conditions}
			\label{sec:Sufficient}%
			\grayout{%
	Notice that the equivalence of \ref{ADMM} and \ref{DRS} stated in \Cref{thm:ADMM2DRS} holds regardless of the properties of functions and matrices in the problem formulation \eqref{eq:CP}.
	Consequently, trusting that \ref{DRS} is not guaranteed to converge unless \(\sDRS\) and \(\nsDRS\) comply with some minimal requirements, imposing some regularity on the image functions \(\epicomp A\nsADMM\) and \(\epicomp B\sADMM\) for \ref{ADMM} cannot be avoided.
	Notice that, however, the requirements of \Cref{ass:FG} are much weaker than other studies in the literature, as they do not rule out `nasty' functions such as 0-norms and discrete or rank constraints.
	To this end, consider the problem
	\[
		\minimize_{(x,z)\in\R^2\times\R}{
			\|x\|_0
			{}-{}
			\tfrac12z^2
		}
		\quad\stt{}
		Ax-z=0
	\]
	where \(A=[1~0]\).
	Then, it is easy to see that this formulation fits our framework, as \(\epicomp A\nsADMM(s)=\|s\|_0\) and \(\epicomp 1\sADMM(s)=-\tfrac12z^2\) for \(\nsADMM=\|{}\cdot{}\|_0\) and \(\sADMM=-\tfrac12({}\cdot{})^2\).
}%
In this section we provide some sufficient conditions on \(f\) and \(g\) ensuring that \Cref{ass:ADMM} is satisfied.

			\subsubsection{Lower semicontinuity of the image function}
				\begin{prop}[Lsc of \(\epicomp Bg\)]\label{thm:lsc:sufficient}%
	Suppose that \Cref{ass:ADMM:lsc,ass:ADMM:beta} are satisfied.
	Then, \(\epicomp Bg\) is proper.
	Moreover, it is also lsc provided that for all \(\bar z\in\dom g\) the set
	\(
		Z(s)
	{}\coloneqq{}
		\argmin_z\set{g(z)}[Bz=s]
	\)
	is nonempty and \(\dist(0,Z(s))\) is bounded for all \(s\in B\dom g\) close to \(B\bar z\).
	\begin{proof}
		Properness is shown in \cref{thm:epiproper}.
		Suppose that \(\seq{s_k}\subseteq\lev_{\leq\alpha}\epicomp Bg\) for some \(\alpha\in\R\) and that \(s_k\to\bar s\).
		Then, due to the characterization of \cite[Thm. 1.6]{rockafellar2011variational} it suffices to show that \(\bar s\in\lev_{\leq\alpha}\epicomp Bg\).
		The assumption ensures the existence of a bounded sequence \(\seq{z_k}\) such that eventually \(Bz_k=s_k\) and \(\epicomp Bg(s_k)=g(z_k)\).
		By possibly extracting, \(z_k\to\bar z\) and necessarily \(B\bar z=\bar s\).
		Then,
		\[
			\epicomp Bg(\bar s)
		{}\leq{}
			g(\bar z)
		\smash{
			{}\leq{}
				\liminf_{k\to\infty}g(z_k)
			{}={}
				\liminf_{k\to\infty}\epicomp Bg(s_k)
			{}\leq{}
		}
			\alpha,
		\]
		hence \(\bar s\in\lev_{\leq\alpha}\epicomp Bg\).
	\end{proof}
\end{prop}
%
%
The requirement in \Cref{thm:lsc:sufficient} is weaker than Lipschitz continuity of the map \(s\mapsto Z(s)\), which is the standing assumption in \cite{wang2018global} for the analysis of ADMM.
In fact, no uniqueness or boundedness of the sets of minimizers is required, but only the existence of minimizers not arbitrarily far.

\begin{wrapfigure}[11]{r}{0.32\linewidth}%
	\vspace*{0pt}%
	{{%
			\pgfkeys{/pgf/images/include external/.code={\includegraphics[width=\linewidth]{#width=\linewidth}}}%
			\tikzsetnextfilename{epicomp_lsc}%
			\input{./TeX/Tikz/epicomp_lsc.tex}%
		}}%
	\vspace*{0pt}%
\end{wrapfigure}
The pathological behavior occurring when this condition is not met can be well visualized by considering \(\func{g}{\R^2}{\R}\) defined as
\begin{equation}\label{eq:nlsc}
	\hspace*{-4cm}
	g(x,y)
{}={}
	\begin{ifcases}
		-|x|			& |xy|\geq1\\
		1-q(|xy|)(1+|x|)\otherwise,
	\end{ifcases}
\end{equation}
where \(q(t)\) is any function such that \(q(0)=0<q(t)<1=q(1)\) for all \(t\in(0,1)\).
On the right, a graphical representation of the piecewise definition on the positive orthant of \(\R^2\) (the function is mirrored in all other orthants).
On the axes, \(f\) achieves its maximum value, that is, \(1\).
In the gray region \(|xy|\geq1\), \(f(x,y)=-|x|\).
In the white portion, \(f\) is extended by means of a convex combination of \(1\) and \(-|x|\).
Function \(g\) and \(B\coloneqq[1~0]\) are ADMM-feasible, meaning that
\(
	\argmin_{w\in\R^2}\set{g(w)+\tfrac\beta2\|Bw-s\|^2}
{}\neq{}
	\emptyset
\)
for all \(s\in\R\) and \(\beta\) large enough (in fact, for all \(\beta>0\), being \(g({}\cdot{},y)+\tfrac\beta2\|{}\cdot{}-s\|^2\) coercive for any \(y\in\R\)).
However,
\(
	\epicomp Bg(s)
{}={}
	-|s|
\)
if \(s\neq 0\) while
\(
	\epicomp Bg(0)
{}={}
	1
\),
resulting in the lack of lsc at \(s=0\).
Along \(\ker B=\set0\times\R\), by keeping \(x\) constant \(g\) attains minimum at \(\set{(x,y)}[xy\geq1]\) for \(x\neq 0\), which escapes to infinity as \(x\to 0\), and \(g(x,x^{-1})=-|x|\to 0\).
However, if instead \(x=0\) is fixed (as opposed to \(x\to 0\)), then the pathology comes from the fact that \(g(0,{}\cdot{})\equiv 1>0\).
The \emph{interpolating} function \(q\) simply models the transition from a constant function on the axes and a linear function in the regions delimited by the hyperbolae.
For any \(k\in\N\) it can thus be chosen such that \(g\) is \(k\) times continuously differentiable; the choice \(q(t)=\tfrac12(1-\cos{\pi t})\), for instance, makes \(g\in C^1(\R^2)\).
In particular, (high-order) continuous differentiability is not enough for \(\epicomp Bg\) to be lsc.

The next result provides necessary and sufficient conditions ensuring the image function \(\epicomp Bg\) to inherit lower semicontinuity from that of \(g\).
It will be evident that pathological cases such as the one depicted in \eqref{eq:nlsc} may only occur due to the behavior of \(g\) at infinity.

\begin{thm}
	For any lsc function \(\func g{\R^n}{\Rinf}\) and \(B\in\R^{p\times n}\), the image function \(\epicomp Bg\) is lsc iff
	\begin{equation}\label{eq:Aflsc:liminf}
		\liminf_{\limsubstack{\|d\|&\to&\infty\\Bd&\to&0}}{
			g(\bar z+d)
		}
	{}\geq{}
		\inf_{d\in\ker B}g(\bar z+d)
	\qquad
		\forall\bar z\in\dom g.
	\end{equation}
	In particular, for any lsc and level bounded function \(\func g{\R^n}{\Rinf}\) and \(B\in\R^{p\times n}\), \(\epicomp Bg\) is lsc.
	\begin{proof}
		Observe first that the right-hand side in \eqref{eq:Aflsc:liminf} is \(\epicomp Bg(B\bar z)\).
		Suppose now that \eqref{eq:Aflsc:liminf} holds, and given \(\bar s\in\dom\epicomp Bg\) consider a sequence \(\seq{s_k}\subseteq\lev_{\leq\alpha}\epicomp Bg\) for some \(\alpha\in\R\) and such that \(s_k\to\bar s\).
		Then, it suffices to show that \(\bar s\in\lev_{\leq\alpha}\epicomp Bg\).
		Let \(\seq{z_k}\) be such that \(Bz_k=s_k\) and \(g(z_k)\leq\epicomp Bg(s_k)+\nicefrac1k\) for all \(k\in\N\).
		If, up to possibly extracting, there exists \(z\) such that \(z^k\to z\) as \(k\to\infty\), then the claim follows with a similar reasoning as in the proof of \cref{thm:lsc:sufficient}.
		Suppose, instead, that \(t_k\coloneqq\|z_k\|\to\infty\) as \(k\to\infty\), and let
		\(
			d_k
		{}\coloneqq{}
			z_k-\bar z
		\),
		where \(\bar z\in\dom g\) is any such that \(B\bar z=s\) (such a \(\bar z\) exists, being \(\bar s\in\dom\epicomp Bg=B\dom g\)).
		Since
		\(
			Bd_k
		{}={}
			B(z_k-\bar z)
		{}={}
			s_k-\bar s
		{}\to{}
			0
		\),
		we have
		\[
			\epicomp Bg(\bar s)
		{}={}
			\inf_{d\in\ker B}{
				g(\bar z+d)
			}
		{}\leq{}
			\liminf_{k\to\infty}{
				g(\bar z+d_k)
			}
		{}={}
			\liminf_{k\to\infty}{
				g(z_k)
			}
		{}\leq{}
			\liminf_{k\to\infty}{
				\epicomp Bg(s_k)+\tfrac1k
			}
		{}\leq{}
			\alpha,
		\]
		proving that \(\bar s\in\lev_{\leq\alpha}\epicomp Bg\).
		
		To show the converse implication, suppose that \eqref{eq:Aflsc:liminf} does not hold.
		Thus, there exist \(\bar z\in\dom g\) and \(\seq{d^k}\subset\R^n\) such that \(Bd^k\to0\) as \(k\to\infty\), and such that, for some \(\varepsilon>0\),
		\[
			g(\bar z+d^k)
			{}+{}
			\varepsilon
		{}\leq{}
			\inf_{d\in\ker B}g(\bar z+d)
		{}={}
			\epicomp Bg(B\bar z)
		\quad
			\forall k.
		\]
		Then, \(s_k\coloneqq B(\bar z+d^k)\) satisfies \(s_k\to B\bar z\) as \(k\to\infty\), and
		\[
			\epicomp Bg(B\bar z)
		{}\geq{}
			\liminf_{k\to\infty}{
				g(\bar z+d^k)
			}
			{}+{}
			\varepsilon
		{}\geq{}
			\liminf_{k\to\infty}{
				\epicomp Bg(s^k)
			}
			{}+{}
			\varepsilon,
		\]
		hence \(\epicomp Bg\) is not lsc at \(B\bar z\).
	\end{proof}
\end{thm}
The \emph{asymptotic function}
\(
	g_\infty(\bar d)
{}\coloneqq{}
	\liminf_{
		d\to\bar d,\,
		t\to\infty
	}{
		\frac{g(td)}{t}
	}
\)
is a tool used in \cite{auslender2002asymptotic} to analyze the behavior of \(g\) at infinity and derive sufficient properties ensuring lsc of \(\epicomp Bg\).
These all ensure that the set of minimizers \(Z(s)\) as defined in \Cref{thm:lsc:sufficient} is nonempty, although this property is not necessary as long as lower semicontinuity is concerned.
To see this, it suffices to modify \eqref{eq:nlsc} as follows
\[
	g(x,y)
{}={}
	\begin{ifcases}
		-|x|							& |xy|\geq1\\
		e^{-y^2}-q(|xy|)(e^{-y^2}+|x|)	\otherwise,
	\end{ifcases}
\]
that is, by replacing the constant value \(1\) on the \(y\) axis with \(e^{-y^2}\).
Then, \(\epicomp Bg(s)=-|s|\) is lsc, but the set of minimizers
\(
	\argmin_w\set{g(w)}[Bw=0]
{}={}
	\set 0\times\argmin_ye^{-y^2}
\)
is empty at \(s=0\).

\grayout{%
	We conclude the section by providing sufficient conditions for the functions \(\epicomp A\nsADMM\) and \(\epicomp B\sADMM\) to comply with the requirements of \Cref{ass:FG}.
	Namely, it suffices that \(\nsADMM\) satisfies any of the requirements of \Cref{prop:epicomplsc}, and that \(\sADMM\) satisfies any of those of \Cref{prop:epicompsmooth}.
	\begin{prop}[Basic properties of image functions]\label{prop:epicomplsc}%
		Let \(\func{\nsADMM}{\R^n}{\Rinf}\) be a proper and lsc function, and for \(A\in\R^{m\times n}\) consider the (possibly set-valued) mapping
		\(
			X(s)
		{}\coloneqq{}
			\argmin_x\set{\nsADMM(x)}[Ax=s]
		\).
		Then, the image function \(\epicomp A\nsADMM\) is proper and lsc as long as either
		\begin{enumerate}
		\item\label{prop:locbound}
			\(X(\bar s)\) is nonempty for all \(\bar s\in A\dom\nsADMM\) and \DEF{locally bounded} (in the sense of \cite[Def. 5.14]{rockafellar2011variational}):
			for all \(\bar s\) there exist \(\varepsilon,r>0\) such that
			\[
				X(s)\subseteq\ball 0r
			\quad
				\forall s\in\ball{\bar s}{\varepsilon},
			\]
		\item\label{prop:levbound}
			or for some \(\beta>0\)
			\begin{equation}\label{eq:levelbounded}
				\nsADMM
				{}+{}
				\tfrac{\beta}{2}
				\|A{}\cdot{}-\bar s\|^2
			\quad
				\text{is level bounded for all \(\bar s\in A\dom\nsADMM\),}
			\end{equation}
			as it is the case when \(\nsADMM\) is prox-bounded and \(A\) is full-column rank.
		\item\label{prop:cvxepicomp}
			or \(\nsADMM\) is (strongly) convex and \(\range\trans A\cap\relint\range\partial g\neq\emptyset\), in which case \(\epicomp A\nsADMM\) is also (strongly) convex.
		\end{enumerate}
		In any case,
		\(
			X(\bar s)\neq\emptyset
		\)
		for all \(\bar s\in\dom\epicomp A\nsADMM=A\dom\nsADMM\).
		Moreover, in case \ref{prop:levbound} one also has
		\begin{equation}\label{eq:epicompSubgrad}
			\partial\epicomp A\nsADMM(\bar s)
		{}\subseteq{}
			\hspace*{-5pt}
			\bigcup_{x\in X(\bar s)}{\hspace*{-5pt}
				\set{y\in\R^p}[\trans Ay\in\partial\nsADMM(x)]
			}
		\qquad
			\forall\bar s\in\dom\epicomp A\nsADMM.
		\end{equation}
		whereas if \ref{prop:cvxepicomp} holds, then the inclusion is actually an equality. 
		\begin{proof}
			See \Cref{proof:prop:epicomplsc}.
		\end{proof}
	\end{prop}
	Notice that the requirement of local boundedness of \Cref{prop:locbound} is much weaker an assumption than Lipschitz continuity of the map \(s\mapsto X(s)\), which is the standing assumption in \cite{wang2015global}.
	In fact, we do not even require uniqueness of the minimizers.
}%

			\subsubsection{Smoothness of the image function}
				We now turn to the smoothness requirement of \(\epicomp Af\).
To this end, we introduce the following notion of \emph{smoothness with respect to a matrix}.
\begin{defin}[Smoothness relative to a matrix]\label{defin:Bsmooth}%
	We say that \(\func{h}{\R^n}{\R}\) is \DEF{smooth relative to a matrix \(C\in\R^{p\times n}\)}, and we write \(h\in C^{1,1}_C(\R^n)\), if \(h\) is differentiable and \(\nabla h\) satisfies the following Lipschitz condition: there exist \(L_{h,C}\) and \(\sigma_{h,C}\) with \(|\sigma_{h,C}|\leq L_{h,C}\) such that
	\begin{equation}\label{eq:sLC}
		\sigma_{h,C}\|C(x-y)\|^2
	{}\leq{}
		\innprod{
			\nabla h(x)-\nabla h(y)
		}{
			x-y
		}
	{}\leq{}
		L_{h,C}\|C(x-y)\|^2
	\end{equation}
	whenever \(\nabla h(x),\nabla h(y)\in\range\trans C\).
\end{defin}
This condition is similar to that considered in \cite{goncalves2017convergence}, where \(\proj_{\range\trans A}\nabla f\) is required to be Lipschitz.
The paper analyzes convergence of a proximal ADMM; standard ADMM can be recovered when matrix \(A\) is invertible, in which case both conditions reduce to Lipschitz differentiability of \(f\).
In general, our condition applies to a smaller set of points only, as it can be verified with \(f(x,y)=\tfrac12x^2y^2\) and \(A=[1~0]\).
In fact, \(\proj_{\range\trans A}\nabla f(x,y)=\binom{xy^2}{0}\) is \emph{not} Lipschitz continuous; however, \(\nabla f(x,y)\in\range\trans A\) iff \(xy=0\), in which case \(\nabla f\equiv 0\).
Then, \(f\) is smooth relative to \(A\) with \(L_{f,A}=0\).

To better understand how this notion of regularity comes into the picture, notice that if \(f\) is differentiable, then \(\nabla f(x)\in\range\trans A\) on some domain \(\mathcal U\) if there exists a differentiable function \(\func{q}{A\mathcal U}{\R}\) such that \(f(x)=q(Ax)\).
Then, it is easy to verify that \(f\) is smooth relative to \(A\) if the local ``reparametrization'' \(q\) is smooth (on its domain).
From an a posteriori perspective, if \(\epicomp Af\) is smooth, then due to the relation
\(
	\trans A\nabla\epicomp Af(Az_s)
{}={}
	\nabla f(z_s)
\)
holding for \(z_s\in\argmin_{z:Az=s}f(z)\) (cf. \cref{thm:epipartial}), it is apparent that \(q\) serves as \(\epicomp Af\).
Therefore, smoothness relative to \(A\) is somewhat a minimal requirement for ensuring smoothness of \(\epicomp Af\).

\begin{thm}[Smoothness of \(\epicomp Af\)]\label{prop:epicompsmooth}%
	Let \(A\in\R^{p\times n}\) be surjective and \(\func f{\R^n}{\R}\) be lsc.
	Suppose that there exists \(\beta\geq 0\) such that the function \(f+\tfrac\beta2\|A{}\cdot{}-s\|^2\) is level bounded for all \(s\in\R^p\).
	Then, the image function \(\epicomp Af\) is smooth on \(\R^p\), provided that either
	\begin{enumerate}
	\item\label{prop:epicompsmooth:CC}
		\(f\in C_A^{1,1}(\R^n)\), in which case
		\(
			L_{\epicomp Af}=L_{f,A}
		\)
		and
		\(
			\sigma_{\epicomp Af}=\sigma_{f,A}
		\),
	\item\label{prop:epicompsmooth:Lip}
		or \(f\in C^{1,1}(\R^n)\), and \(X(s)\coloneqq\argmin\set{f(x)}[Ax=s]\) is single valued and Lipschitz continuous with modulus \(M\), in which case
		\[
			L_{\epicomp Af}=L_fM^2
		\quad\text{and}\quad
			\sigma_{\epicomp Af}
		{}={}
			\begin{ifcases}
				\nicefrac{\sigma_f}{\|A\|^2}
			&
				\sigma_f\geq 0
			\\
				\sigma_fM^2
				\otherwise[\(\sigma_f<0\);]
			\end{ifcases}
		\]
	\item\label{prop:epicompsmooth:cvx}
		or \(f\in C^{1,1}(\R^n)\) is convex, in which case
		\(
			L_{\epicomp Af}
		{}={}
			\frac{L_f}{\sigma_+(\trans AA)}
		\)
		and
		\(
			\sigma_{\epicomp Af}
		{}={}
			\nicefrac{\sigma_f}{\|A\|^2}
		\).
	\end{enumerate}
	\begin{proof}
		As shown in \cref{thm:epiproper}, \(\epicomp Af\) is proper.
		The surjectivity of \(A\) and the level boundedness condition ensure that for all \(\alpha\in\R\) and \(s\in\R^p\) the set
		\(
			\set{x}[f(x)\leq\alpha,\,\|Ax-s\|<\varepsilon]
		\)
		is bounded for some \(\varepsilon>0\) (in fact, for all \(\varepsilon>0\)).
		Then, we may invoke \cite[Thm. 1.32]{rockafellar2011variational} to infer that \(\epicomp Af\) is lsc, that the set \(X(s)\coloneqq\argmin_x\set{f(x)}[Ax=s]\) is nonempty for all \(s\in\R^p\), and that the function
		\(
			H(z,s)
		{}\coloneqq{}
			f(x)+\indicator_{\set0}(Ax-s)
		\)
		is uniformly level bounded in \(x\) locally uniformly in \(s\), in the sense of \cite[Def. 1.16]{rockafellar2011variational}.
		Moreover, since \(f\) is differentiable, observe that
		\(
			\partial^\infty H(x,Ax)
		{}={}
			\range\binom{\trans A}{\I}
		\)
		for all \(x\in\R^m\).
		Hence, for all \(s\in\R^p\) it holds that
		\[
			\partial^\infty\epicomp Af(s)
		{}\subseteq{}
			\bigcup_{\mathclap{x\in X(s)}}{
				\set{y}[
					(0,y)
				{}\in{}
					\partial^\infty H(x,s)
				]
			}
		{}={}
			\ker{\trans A}
		{}={}
			\set0,
		\]
		where the inclusion follows from \cite[Thm. 10.13]{rockafellar2011variational}.
		By virtue of \cite[Thm. 9.13]{rockafellar2011variational}, we conclude that \(\epicomp Af\) is strictly continuous and has nonempty subdifferential on \(\R^p\).
		Fix \(s_i\in\R^p\) and \(y_i\in\partial\epicomp Af(s_i)\), \(i=1,2\), and let us proceed by cases.
		\begin{proofitemize}
		\item\ref{prop:epicompsmooth:CC} and \ref{prop:epicompsmooth:Lip}.
			It follows from \cref{thm:epipartial} and continuous differentiability of \(f\) that
			\(
				\trans Ay_i
			{}\in{}
				\partial f(x_i)
			{}={}
				\set{\nabla f(x_i)}
			\),
			for some \(x_i\in X(s_i)\), \(i=1,2\).
			We have
			\begin{align*}
				\innprod{y_1-y_2}{s_1-s_2}
			{}={} &
				\innprod{y_1-y_2}{Ax_1-Ax_2}
			{}={}
				\innprod{\trans Ay_1-\trans Ay_2}{x_1-x_2}
			\\
			\numberthis\label{eq:epicompInnprod}
			{}={} &
				\innprod{
					\nabla f(x_1)-\nabla f(x_2)
				}{x_1-x_2}.
			\end{align*}
			If \ref{prop:epicompsmooth:CC} holds, since \(\nabla f(x_i)=\trans Ay_i\in\range\trans A\), \(i=1,2\), smoothness of \(f\) relative to \(A\) implies
			\begin{align*}
				\sigma_{f,A}\|s_1-s_2\|^2
			{}={} &
				\sigma_{f,A}\|Ax_1-Ax_2\|^2
			\\
			{}\leq{} &
				\innprod{y_1-y_2}{s_1-s_2}
			{}\leq{}
				L_{f,A}\|Ax_1-Ax_2\|^2
			{}={}
				L_{f,A}\|s_1-s_2\|^2
			\end{align*}
			for all \(s_i\in\R^p\) and \(y_i\in\partial\epicomp Af(s_i)\), \(i=1,2\).
			Otherwise, if \ref{prop:epicompsmooth:Lip} holds, then
			\[
				\sigma_f\|x_1-x_2\|^2
			{}\leq{}
				\innprod{y_1-y_2}{s_1-s_2}
			{}\leq{}
				L_f\|x_1-x_2\|^2
			\]
			and from the bound
			\(
				\frac{1}{\|A\|}\|s_1-s_2\|
			{}\leq{}
				\|x_1-x_2\|
			{}\leq{}
				M\|s_1-s_2\|
			\)
			we obtain
			\[
				\sigma_{\epicomp Af}\|s_1-s_2\|^2
			{}\leq{}
				\innprod{y_1-y_2}{s_1-s_2}
			{}\leq{}
				L_{\epicomp Af}\|s_1-s_2\|^2
			\]
			with the constants \(\sigma_{\epicomp Af}\) and \(L_{\epicomp Af}\) as in the statement.
			Smoothness and hypoconvexity then follow by invoking \cref{thm:C11subdiff}.
		\item\ref{prop:epicompsmooth:cvx}.
			It follows from \cite[Thm. D.4.5.1 and Cor. D.4.5.2]{hiriarturruty2012fundamentals} that \(\epicomp Af\) is convex and differentiable, and satisfies
			\(
				\nabla\epicomp Af(s)
			{}={}
				y
			\),
			where for any \(x\in X(s)\), \(y\) is such that \(\trans Ay=\nabla f(x)\).
			For \(y_i=\nabla\epicomp Af(s_i)\) and \(x_i\in X(s_i)\), \(i=1,2\), the equalities in \eqref{eq:epicompInnprod} hold.
			In turn,
			\[
				\innprod{s_1-s_2}{y_1-y_2}
			{}\geq{}
				\tfrac{1}{L_f}\|\trans A(y_1-y_2)\|^2
			{}\geq{}
				\tfrac{\sigma_+(\trans AA)}{L_f}
				\|\proj_{\range A}(y_1-y_2)\|^2
			{}={}
				\tfrac{\sigma_+(\trans AA)}{L_f}
				\|y_1-y_2\|^2,
			\]
			where the first inequality is due to \(\nicefrac{1}{L_f}\)-cocoercivity of \(\nabla f\), see \cite[Thm. 2.1.5]{nesterov2003introductory}, the second inequality is a known fact (see \eg \cite[Lem. A.2]{goncalves2017convergence}), and the equality is due to the fact that \(A\) is surjective.
			We may again invoke \cite[Thm. 2.1.5]{nesterov2003introductory} to infer the claimed \(\tfrac{L_f}{\sigma_+(\trans AA)}\)-smoothness of \(\epicomp Af\).
			Since \(\epicomp Af\) is convex (thus \(0\)-hypoconvex), if \(\sigma_f=0\) there is nothing more to show.
			The case \(\sigma_f>0\) follows from \cref{thm:epicompStrCvx}.
		\qedhere
		\end{proofitemize}
	\end{proof}
\end{thm}
Notice that the condition in \Cref{prop:epicompsmooth:Lip} covers the case when \(f\in C^{1,1}(\R^n)\) and \(A\) has full column rank (hence is invertible), in which case \(M=\nicefrac{1}{\sigma_+(A)}\).
This is somehow trivial, since necessarily \(\epicomp Af(s)=f\circ A^{-1}\) in this case.

\grayout{%
	The Lipschitz continuity of the minimizers assumed in \Cref{prop:epicompsmooth:Lip} is required in \cite{wang2015global} for both the smooth and the nonsmooth function; on the contrary, in order to comply with \Cref{ass:FG} we need it (if ever) only for the smooth function \(f\), whereas for the nonsmooth term \(g\) much weaker sufficient conditions are available, as shown in \Cref{prop:flsc}.
}%

	\section{Conclusive remarks}
		\label{sec:Conclusion}%
		This paper provides new convergence results for nonconvex Douglas-Rachford splitting (DRS) and ADMM with an all-inclusive analysis of all possible relaxation parameters \(\lambda\in(0,4)\).
Under the only assumption of Lipschitz differentiability of one function, convergence is shown for larger prox-stepsizes and relaxation parameters than was previously known.
The results are tight when \(\lambda\in(0,2]\), covering in particular classical (non-relaxed) DRS and PRS, or when the differentiable function is nonconvex.
The necessity of \(\lambda<4\) and of a lower bound for the stepsize when \(\lambda>2\) is also shown.

Our theory is based on the Douglas-Rachford envelope (DRE), a continuous, real-valued, exact penalty function for DRS, and on a primal equivalence of DRS and ADMM that extends the well-known connection of the algorithms to arbitrary (nonconvex) problems.
The DRE is shown to be a better Lyapunov function for DRS than the augmented Lagrangian, due to its closer connections with the cost function and with DRS iterations.

	
	\begin{appendix}
		\proofsection{sec:Background}\label{sec:proof:Background}

\begin{appendixproof}{thm:C11subdiff}%
	The claimed hy\-po\-co\-nve\-xi\-ty follows from \cite[Ex. 12.28]{rockafellar2011variational}.
	It suffices to show that \(h\) is continuously differentiable, so that \(\partial h=\nabla h\) and the claim then follows from \eqref{eq:innprod}.
	To this end, without loss of generality we may assume that \(\sigma\geq0\), since \(h\) is continuously differentiable iff so is \(h-\tfrac\sigma2\|{}\cdot{}\|^2\).
	Thus, for all \(x_i\in\R^n\), \(v_i\in\partial h(x_i)\), \(i=1,2\), one has
	\begin{align*}
		h(x_1)
	{}\geq{} &
		h(x_2)
		{}+{}
		\innprod{v_2}{x_1-x_2}
	{}={}
		h(x_2)
		{}+{}
		\innprod{v_2-v_1}{x_1-x_2}
		{}+{}
		\innprod{v_1}{x_1-x_2}
	\\
	{}\geq{} &
		h(x_2)
		{}-{}
		L\|x_1-x_2\|^2
		{}+{}
		\innprod{v_1}{x_1-x_2},
	\end{align*}
	where the first inequality follows from convexity of \(h\) (being it \(0\)-hypoconvex by assumption).
	Rearranging,
	\[
		h(x_2)
	{}\leq{}
		h(x_1)
		{}+{}
		\innprod{v_1}{x_2-x_1}
		{}+{}
		L\|x_1-x_2\|^2
	\quad
		\forall x_i\in\R^n,~
		v_1\in\partial h(x_1),~
		i=1,2.
	\]
	Let \(\tilde h\coloneqq h-\innprod{v_1}{{}\cdot{}}\), so that \(0\in\partial h(x_1)\).
	Due to convexity, \(x_1\in\argmin\tilde h\), hence for all \(w\in\R^n\) and \(v_1'\in\partial h(x_1)\) one has
	\[
		\tilde h(x_1)
	{}\leq{}
		\tilde h(w)
	{}\leq{}
		h(x_1)
		{}+{}
		\innprod{v_1'}{w-x_1}
		{}+{}
		L\|w-x_1\|^2
		-\innprod vw
	{}={}
		\tilde h(x_1)
		{}+{}
		\innprod{v_1'-v_1}{w-x_1}
		{}+{}
		L\|w-x_1\|^2.
	\]
	By selecting \(w=x_1-\tfrac{1}{2L}(v_1'-v_1)\), one obtains \(\|v_1-v_1'\|^2\leq0\), hence necessarily \(v_1=v_1'\).
	From the arbitrarity of \(x_1\in\R^n\) and \(v_1,v_1'\in\partial h(x_1)\) it follows that \(\partial h\) is everywhere single valued, and the sought continuous differentiability of \(h\) then follows from \cite[Cor. 9.19]{rockafellar2011variational}.
\end{appendixproof}

\begin{appendixproof}{thm:rho}%
	\begin{proofitemize}
	\item\ref{thm:hypo:rho0}.
		This is the lower bound in \eqref{eq:LipBound}.
	\item\ref{thm:hypo:rho1}.
		Let \(L\geq L_h\) and \(\sigma\in(-L,\min\set{0,\sigma_h}]\) be fixed.
		Then, \(h\) is \(L\)-smooth and \(\sigma\)-hypoconvex, and from \cite[Thm. 2.1.12]{nesterov2003introductory} we obtain that
		\begin{equation}\label{eq:LipBound:mix}
			\innprod{
				\nabla h(y)-\nabla h(x)
			}{
				y-x
			}
		{}\geq{}
			\tfrac{\sigma L}{L+\sigma}\|x-y\|^2
			{}+{}
			\tfrac{1}{L+\sigma}\|\nabla h(x)-\nabla h(y)\|^2
		\end{equation}
		for all \(x,y\in\R^n\).
		(Although \cite[Thm. 2.1.12]{nesterov2003introductory} assumes \(\sigma>0\), the given proof does not necessitate this restriction).
		Moreover, \(\psi\coloneqq h-\tfrac\sigma2\|{}\cdot{}\|^2\) is convex and \(L_\psi\)-smooth, with \(L_\psi=L-\sigma\).
		Consequently, for all \(x,y\in\R^n\) one has
		\(
			\psi(y)
		{}\geq{}
			\psi(x)
			{}+{}
			\innprod{\nabla\psi(x)}{y-x}
			{}+{}
			\tfrac{1}{2L_\psi}\|\nabla\psi(y)-\nabla\psi(x)\|^2
		\),
		see \cite[Thm. 2.1.5]{nesterov2003introductory}, resulting in
		\begin{align*}
			h(y)
		{}\geq{} &
			h(x)
			{}+{}
			\innprod{\nabla h(x)}{y-x}
			{}+{}
			\tfrac{\sigma L}{2(L-\sigma)}
			\|y-x\|^2
			{}+{}
			\tfrac{1}{2(L-\sigma)}
			\|\nabla h(y)-\nabla h(x)\|^2
		\\
		&
			{}-{}
			\tfrac{\sigma}{2(L-\sigma)}
			\innprod{
				\nabla h(y)-\nabla h(x)
			}{
				y-x
			}.
		\end{align*}
		Since \(\sigma\leq0\), the coefficient of the scalar product in the second line is positive.
		We may thus invoke the inequality \eqref{eq:LipBound:mix} to arrive to
		\begin{align*}
			h(y)
		{}\geq{} &
			h(x)
			{}+{}
			\innprod{\nabla h(x)}{y-x}
			{}+{}
			\tfrac{\sigma L}{2(L-\sigma)}
			\|y-x\|^2
			{}+{}
			\tfrac{1}{2(L-\sigma)}
			\|\nabla h(y)-\nabla h(x)\|^2
		\\
		&
			{}-{}
			\tfrac{\sigma}{2(L-\sigma)}
			\left[
				\tfrac{\sigma L}{L+\sigma}\|x-y\|^2
				{}+{}
				\tfrac{1}{L+\sigma}\|\nabla h(x)-\nabla h(y)\|^2
			\right]
		\\
		{}={} &
			h(x)
			{}+{}
			\innprod{\nabla h(x)}{y-x}
			{}+{}
			\tfrac{\sigma L}{2(L+\sigma)}
			\|y-x\|^2
			{}+{}
			\tfrac{1}{2(L+\sigma)}
			\|\nabla h(y)-\nabla h(x)\|^2,
		\end{align*}
		hence the claimed inequality.
	\qedhere
	\end{proofitemize}
\end{appendixproof}

\begin{appendixproof}{thm:proxf}%
	Let \(\gamma\in(0,\nicefrac{1}{[\sigma_h]_-})\) be fixed, and let \(\psi\coloneqq\gamma h+\tfrac12\|{}\cdot{}\|^2\).
	Observe that \(\psi\in\cont^{1,1}(\R^n)\) is \(L_\psi\)-smooth and \(\sigma_\psi\)-strongly convex, with \(L_\psi=1+\gamma L_h\) and \(\sigma_\psi=1+\gamma\sigma_h\).
	In particular, due to strong convexity \(\inf\psi>-\infty\), and by definition of prox-boundedness it then follows that \(\gamma_h\geq\nicefrac{1}{[\sigma_h]_-}\).
	\begin{proofitemize}
	\item\ref{thm:proxfEquiv}.
		Follows from \eqref{eq:prox:subdiff}, by observing that \(h+\tfrac{1}{2\gamma}\|{}\cdot{}-s\|^2\) is strongly convex, hence that a minimizer is characterized by stationarity.
	\item\ref{thm:proxfLip}.
		For \(s,s'\in\R^n\), let \(u=\prox_{\gamma h}(s)\) and \(u'=\prox_{\gamma h}(s')\).
		Then,
		\[
			\innprod{s-s'}{u-u'}
		{}={}
			\innprod{
				\nabla\psi(u)
				{}-{}
				\nabla\psi(u')
			}{u-u'}
		{}\geq{}
			\sigma_\psi
			\|u-u'\|^2
		{}={}
			(1+\gamma\sigma_h)
			\|u-u'\|^2,
		\]
		where the first equality was shown in \ref{thm:proxfEquiv} and the inequality follows from \eqref{eq:innprod}.
		By using the \(\frac{1}{L_\psi}\)-cocoercivity of \(\nabla\psi\) \cite[Thm. 2.1.10]{nesterov2003introductory}, also the claimed strong monotonicity follows.
		In turn, the Cauchy-Schwartz inequality on the inner product yields \eqref{eq:biLip}.
	\item\ref{thm:MoreaufC1}.
		From \cite[Ex. 10.32]{rockafellar2011variational} it follows that \(h^\gamma\) is strictly continuous and that
		\(
			\partial h^\gamma(s)
		{}\subseteq{}
			\tfrac1\gamma(s-\prox_{\gamma h}(s))
		\).
		Because of single valuedness of \(\prox_{\gamma h}\), by invoking \cite[Thm. 9.18]{rockafellar2011variational} we conclude that \(h^\gamma\) is everywhere differentiable with
		\(
			\nabla h^\gamma(s)
		{}={}
			\tfrac1\gamma(s-\prox_{\gamma h}(s))
		\).
		Thus,
		\begin{align*}
			\innprod{\nabla h^\gamma(s)-\nabla h^\gamma(s')}{s-s'}
		{}={} &
			\tfrac{1}{\gamma}
			\left(
				\|s-s'\|^2
				{}-{}
				\innprod{s-s'}{u-u'}
			\right),
		\end{align*}
		and from the bounds in \ref{thm:proxfLip} we conclude that
		\[
			\tfrac{\sigma_{h}}{1+\gamma\sigma_{h}}
			\|s-s'\|^2
		{}\leq{}
			\innprod{\nabla h^\gamma(s)-\nabla h^\gamma(s')}{s-s'}
		{}\leq{}
			\tfrac{L_{h}}{1+\gamma L_{h}}
			\|s-s'\|^2.
		\]
		The claimed smoothness and hypoconvexity follow from the characterization of \eqref{eq:innprod}.
		\qedhere
	\end{proofitemize}
\end{appendixproof}

		\proofsection{sec:ADMM}\label{sec:proof:ADMM}
			\begin{appendixproof}{thm:lsc}
	\begin{proofitemize}
	\item\ref{thm:epiproper}.~%
		If \(\bar s\notin C\dom h\), then \(\epicomp Ch(\bar s)=\infty\).
		Otherwise, suppose \(\bar s=C\bar x\) for some \(\bar x\in\dom h\).
		Then,
		\[
			-\infty
		{}<{}
			\min_x\set{
				h(x)
				{}+{}
				\tfrac\beta2\|Cx-\bar s\|^2
			}
		{}\leq{}
			\inf_{x:\,Cx=\bar s}\set{
				h(x)
				{}+{}
				\tfrac\beta2\|Cx-\bar s\|^2
			}
		{}\defeq{}
			\epicomp Ch(\bar s),
		\]
		which is upper bounded by the finite quantity \(h(\bar x)\).
	\item\ref{thm:epiexact}.~%
		Since \(C(x_\beta+v)=Cx_\beta\) iff \(v\in\ker C\), for all \(s\in\R^p\) and \(x_\beta\in X_\beta(s)\) necessarily
		\(
			h(x_\beta)
		{}\leq{}
			h(x_\beta+v)
		\).
		Consequently,
		\[
			\epicomp Ch(Cx_\beta)
		{}\leq{}
			h(x_\beta)
		{}\leq{}
			\inf_{v\in\ker C}{
				h(x_\beta+v)
			}
		{}={}
			\inf_{
				x:\,Cx=Cx_\beta
			}{
				h(x)
			}
		{}={}
			\epicomp Ch(Cx_\beta).
		\]
	\item\ref{thm:AX}.~%
		Fix \(\bar s\in\R^p\), and let \(x_\beta\in X_\beta(\bar s)\).
		Then, from \ref{thm:epiexact} and the optimality of \(x_\beta\) we have
		\[
			\epicomp Ch(Cx_\beta)
			{}+{}
			\tfrac\beta2\|Cx_\beta-\bar s\|^2
		{}={}
			h(x_\beta)
			{}+{}
			\tfrac\beta2\|Cx_\beta-\bar s\|^2
		{}\leq{}
			h(x)
			{}+{}
			\tfrac\beta2\|Cx-\bar s\|^2
		\]
		for all \(x\in\R^n\).
		In particular, this holds for all \(s\in\R^p\) and \(x\) such that \(Cx=s\), hence
		\[
			\epicomp Ch(Cx_\beta)
			{}+{}
			\tfrac\beta2\|Cx_\beta-\bar s\|^2
		{}\leq{}\!
			\inf_{x:Cx=s}\set{
				h(x)
				{}+{}
				\tfrac\beta2\|Cx-\bar s\|^2
			}
		{}={}
			\epicomp Ch(s)
			{}+{}
			\tfrac\beta2\|s-\bar s\|^2,
		\]
		proving that
		\(
			Cx_\beta
		{}\in{}
			\prox_{\nicefrac{\epicomp Ch}{\beta}}(\bar s)
		\).
\grayout{%
			We now show local boundedness.
			Fix \(\beta'\in(\bar\beta,\beta)\).
			For \(\bar s\in\R^p\) and \(s\in\ball{\bar s}{\varepsilon}\) we have
			\begin{align*}
				h(x)
				{}+{}
				\tfrac\beta2\|Cx-s\|^2
			{}\geq{} &
				h(x)
				{}+{}
				\tfrac\beta2\|Cx-\bar s\|^2
				{}-{}
				\tfrac\beta2\varepsilon^2
				{}-{}
				\beta\varepsilon\|Cx-s\|
			\\
			{}={} &
				\smashunderbracket{
					h(x)
					{}+{}
					\tfrac{\beta'}{2}\|Cx-\bar s\|^2
				}{
					\text{lower bounded}
				}
				{}+{}
				\smashunderbracket{
					\tfrac{\beta-\beta'}{2}\|Cx-\bar s\|^2
				}{
					\text{quadratic in \(Cx\)}
				}
				{}-{}
				\smashunderbracket{
					\bigl(
						\beta\varepsilon\|Cx-\bar s\|
						{}+{}
						\tfrac{3\beta}{2}\varepsilon^2
					\bigr)
				}{
					\text{linear growth in \(Cx\).}
				}
			\end{align*}
			Therefore,
			\[
				\lim_{\|Cx\|\to\infty}{
					\inf_{s\in\ball{\bar s}{\varepsilon}}\set{
						h(x)+\tfrac\beta2\|Cx-s\|^2
					}
				}
			{}={}
				\infty
			\]
			and local boundedness of \(CX_\beta\) follows.
		\item\ref{thm:penalty}.~%
			This will follow immediately from \ref{thm:AX} once we show that for any proper function \(\func{q}{\R^n}{\Rinf}\) (not necessarily lsc) and \(\bar x\in\dom q\) it holds that
			\begin{equation}
				-\infty
			{}<{}
				\lim_{\gamma\to 0^+} q(x_\gamma)
			{}={}
				\lim_{\gamma\to 0^+} q(x_\gamma)+\tfrac{1}{2\gamma}\|x_\gamma-\bar x\|^2
			{}={}
				\liminf_{x\to\bar x} q(x)
			{}<{}
				\infty
			\end{equation}
			whenever \(x_\gamma\in\prox_{\gamma q}(\bar x)\).
			That the \(\liminf\) is finite is due to the fact that \(q\) is proper and that \(\bar x\in\dom q\).
			Notice that
			\[
				\liminf_{\gamma\to 0^+} q(x_\gamma)+\tfrac{1}{2\gamma}\|x_\gamma-\bar x\|^2
			{}\geq{}
				\liminf_{\gamma\to 0^+} q(x_\gamma)
			{}\geq{}
				\liminf_{x\to\bar x} q(x),
			\]
			where the last inequality follows from the fact that \(x_\gamma\to\bar x\) as \(\gamma\to 0^+\).
			Therefore, it suffices to prove that
			\[
				\limsup_{\gamma\to 0^+} q(x_\gamma)+\tfrac{1}{2\gamma}\|x_\gamma-\bar x\|^2
			{}\leq{}
				\liminf_{x\to\bar x} q(x).
			\]
			To avoid trivialities, let us assume that \(x_\gamma\neq\bar x\) for all \(\gamma\), and contrary to the claim suppose that for all \(\gamma\) there exists a point \(\omega_\gamma\) arbitrarily close to \(\bar x\) such that
			\(
				q(x_\gamma)+\tfrac{1}{2\gamma}\|x_\gamma-\bar x\|^2
			{}\geq{}
				q(\omega_\gamma)+\varepsilon
			\).
			In particular, we can take \(\omega_\gamma\in\ball{\bar x}{\sqrt{2\gamma\varepsilon}}\), so that
			\[
				q(\omega_\gamma)
				{}+{}
				\tfrac{1}{2\gamma}\|\omega_\gamma-\bar x\|^2
			{}\leq{}
				q(x_\gamma)
				{}+{}
				\tfrac{1}{2\gamma}\|x_\gamma-\bar x\|^2
				{}-{}
				\varepsilon
				{}+{}
				\tfrac{1}{2\gamma}\|\omega_\gamma-\bar x\|^2
			{}<{}
				q(x_\gamma)
				{}+{}
				\tfrac{1}{2\gamma}\|x_\gamma-\bar x\|^2,
			\]
			contradicting minimality of \(x_\gamma\in\prox_{\gamma q}(\bar x)\).
		\item\ref{thm:epilsc}.~%
			For \(\beta>\bar\beta\) let
			\(
				F_\beta(x,s)
			{}\coloneqq{}
				h(x)+\tfrac\beta2\|Cx-s\|^2
			\),
			so that \(\phi_\beta(s)=\min_xF_\beta(x,s)\).
			Clearly, \(\varphi_\beta\) are proper functions; we now show that they are also lsc.
			Fix \(\alpha\in\R\) and let \(\seq{s_k}\) be such that \(s_k\in\lev_{\leq\alpha}\phi_\beta\) for all \(k\).
			Suppose that \(s_k\to\bar s\) for some \(\bar s\), and for \(k\in\N\) let \(x_k\in X_\beta(s_k)\); then,
			\begin{align*}
				\phi_\beta(\bar s)
			{}\leq{} &
				h(x_k)
				{}+{}
				\tfrac\beta2
				\|Cx_k-\bar s\|^2
			\\
			{}={} &
				h(x_k)
				{}+{}
				\tfrac\beta2
				\|Cx_k-s_k\|^2
				{}+{}
				\tfrac\beta2
				\|s_k-\bar s\|^2
				{}+{}
				\beta
				\innprod{Cx_k-s_k}{s_k-\bar s}
			\\
			{}={} &
				\phi_\beta(s_k)
				{}+{}
				\smashunderbracket{
					\tfrac\beta2
					\|s_k-\bar s\|^2
					{}+{}
					\beta
					\innprod{Cx_k-s_k}{s_k-\bar s}
				}{
					\varepsilon_k
				}
			\\
			{}\leq{} &
				\alpha
				{}+{}
				\varepsilon_k.
			\end{align*}
			Since \(s_k\to\bar s\), from \ref{thm:AX} it follows that \(\seq{Cx_k}\) is bounded, and therefore \(\varepsilon_k\to0\).
			This proves that \(\bar s\in\lev_{\leq\alpha}\phi_\beta\), that is, that \(\phi_\beta\) is lsc.
			Since \(\seq{\phi_\beta}[\beta>\bar\beta]\) is an increasing sequence of lsc functions, we may invoke \cite[Prop.s 7.4(a) and 7.4(d)]{rockafellar2011variational} to infer that
			\(
				\phi_\star
			{}\coloneqq{}
				\sup_{\beta>\bar\beta}\phi_\beta
			\)
			is an lsc function.
			Clearly, \(\phi_\star\leq\epicomp Ch\) and \(\dom\phi_\star=C\dom h=\dom\epicomp Ch\).
			Suppose now that \(\epicomp Ch(\bar s)=\phi_\star(\bar s)\).
			Then,
			\[
				\liminf_{s\to\bar s}\epicomp Ch(s)
			{}\geq{}
				\liminf_{s\to\bar s}\phi_\star(s)
			{}\geq{}
				\phi_\star(\bar s)
			{}={}
				\epicomp Ch(\bar s)
			\]
			proving that \(\epicomp Ch\) is lsc at \(\bar s\).
			
			Conversely, suppose that \(\epicomp Ch\) is lsc at \(\bar s\).
			Then, since \(Cx_\beta\to\bar s\) we have
			\[
				\epicomp Ch(\bar s)
			{}\leq{}
				\liminf_{s\to\bar s}\epicomp Ch(s)
			{}\leq{}
				\liminf_{\beta\to\infty}\epicomp Ch(Cx_\beta)
			{}\overrel{\ref{thm:AX}}{} 
				\lim_{\beta\to\infty}h(x_\beta)
			{}\overrel{\ref{thm:penalty}}{}
				\phi_\star(\bar s)
			{}\leq{}
				\epicomp Ch(\bar s)
			\]
			proving the claimed necessary condition.
}%
	\qedhere
	\end{proofitemize}
\end{appendixproof}

\begin{appendixproof}{thm:epipartial}%
	Let \(\bar v\in\hat\partial\epicomp Ch(C\bar x)\).
	Then,
	\begin{align*}
	&
		\liminf_{\limsubstack{x&\to&\bar x\\x&\neq&\bar x}}{
			\frac{
				h(x)-h(\bar x)-\innprod{\trans C\bar v}{x-\bar x}
			}{
				\|x-\bar x\|
			}
		}
	\\
	{}={} &
		\liminf_{\limsubstack{x&\to&\bar x\\x&\neq&\bar x}}{
			\frac{
				h(x)-\epicomp Ch(C\bar x)-\innprod{\bar v}{C(x-\bar x)}
			}{
				\|x-\bar x\|
			}
		}
	\\
	{}\geq{} &
		\liminf_{\limsubstack{x&\to&\bar x\\x&\neq&\bar x}}{
			\frac{
				\epicomp Ch(Cx)-\epicomp Ch(C\bar x)-\innprod{\bar v}{C(x-\bar x)}
			}{
				\|x-\bar x\|
			}
		}
	\\
	{}={} &
		\liminf_{\limsubstack{x&\to&\bar x\\x&\neq&\bar x}}{
			\frac{
				\epicomp Ch(Cx)-\epicomp Ch(C\bar x)-\innprod{\bar v}{C(x-\bar x)}
			}{
				\|C(x-\bar x)\|
			}
			\frac{
				\|C(x-\bar x)\|
			}{
				\|x-\bar x\|
			}
		}
	{}\geq{}
		0,
	\end{align*}
	where the last inequality follows from the inclusion \(\bar v\in\hat\partial\epicomp Ch(C\bar x)\).
\end{appendixproof}

\begin{appendixproof}{thm:epicompStrCvx}%
	Convexity of the image function follows from \cite[Prop. 12.36(ii)]{bauschke2017convex}.
	Moreover, due to strong convexity, for every \(s\in C\dom h=\dom\epicomp Ch\) there exists a unique \(x_s\in\R^n\) such that \(Cx_s=s\) and \(\epicomp Ch(s)=h(x_s)\).
	Let \(v_s\in\partial\epicomp Ch(s)\).
	Then, it follows from \cref{thm:epipartial} that \(\trans Cv_s\in\partial h(x_s)\), hence, for all \(s'\in\dom\epicomp Ch\)
	\[
		h(x_{s'})
	{}\geq{}
		h(x_s)
		{}+{}
		\innprod{\trans Cv_s}{x_{s'}-x_s}
		{}+{}
		\tfrac{\sigma_h}{2}\|x_{s'}-x_s\|^2
	{}\geq{}
		h(x_s)
		{}+{}
		\innprod{v_s}{s'-s}
		{}+{}
		\tfrac{\sigma_h}{2\|C\|^2}\|s'-s\|^2.
	\]
	Strong convexity then follows by observing that \(h(x_s)=\epicomp Ch(s)\) and \(h(x_{s'})=\epicomp Ch(s')\).
\end{appendixproof}

	\end{appendix}


\bibliographystyle{plain}
\bibliography{TeX/Bibliography.bib}

\begin{thebibliography}{10}

\bibitem{attouch2013convergence}
Hedy Attouch, J\'er\^ome Bolte, and Benar~Fux Svaiter.
\newblock Convergence of descent methods for semi-algebraic and tame problems:
  proximal algorithms, forward-backward splitting, and regularized
  {G}auss-{S}eidel methods.
\newblock {\em Mathematical Programming}, 137(1):91--129, Feb 2013.

\bibitem{auslender2002asymptotic}
Alfred Auslender and Marc Teboulle.
\newblock {\em Asymptotic Cones and Functions in Optimization and Variational
  Inequalities}.
\newblock Springer Monographs in Mathematics. Springer New York, 2002.

\bibitem{bauschke2014local}
Heinz Bauschke and Dominikus Noll.
\newblock On the local convergence of the {D}ouglas-{R}achford algorithm.
\newblock {\em Archiv der Mathematik}, 102(6):589--600, Jun 2014.

\bibitem{bauschke2017convex}
Heinz~H. Bauschke and Patrick~L. Combettes.
\newblock {\em Convex analysis and monotone operator theory in {H}ilbert
  spaces}.
\newblock CMS Books in Mathematics. Springer, 2017.

\bibitem{bauschke2015projection}
Heinz~H. Bauschke and Valentin~R. Koch.
\newblock Projection methods: {S}wiss army knives for solving feasibility and
  best approximation problems with halfspaces.
\newblock In Simeon Reich and Alexander~J. Zaslavski, editors, {\em Infinite
  Products of Operators and Their Applications}, volume 636, pages 1--40.
  American Mathematical Society, 2015.

\bibitem{bauschke2014method}
Heinz~H. Bauschke, Hung~M. Phan, and Xianfu Wang.
\newblock The method of alternating relaxed projections for two nonconvex sets.
\newblock {\em Vietnam Journal of Mathematics}, 42(4):421--450, Dec 2014.

\bibitem{bertsekas2016nonlinear}
Dimitri~P. Bertsekas.
\newblock {\em Nonlinear Programming}.
\newblock Athena Scientific, 2016.

\bibitem{bochnak2013real}
Jacek Bochnak, Michel Coste, and Marie-Françoise Roy.
\newblock {\em Real Algebraic Geometry}.
\newblock A Series of Modern Surveys in Mathematics. Springer Berlin
  Heidelberg, 2013.

\bibitem{bolte2014proximal}
J{\'e}r{\^o}me Bolte, Shoham Sabach, and Marc Teboulle.
\newblock Proximal {A}lternating {L}inearized {M}inimization for nonconvex and
  nonsmooth problems.
\newblock {\em Mathematical Programming}, 146(1--2):459--494, 2014.

\bibitem{boyd2011distributed}
Stephen Boyd, Neal Parikh, Eric Chu, Borja Peleato, and Jonathan Eckstein.
\newblock Distributed optimization and statistical learning via the alternating
  direction method of multipliers.
\newblock {\em Found. Trends Mach. Learn.}, 3(1):1--122, January 2011.

\bibitem{douglas1956numerical}
Jim Douglas and Henry~H. Rachford.
\newblock On the numerical solution of heat conduction problems in two and
  three space variables.
\newblock {\em Transactions of the American Mathematical Society},
  82(2):421--439, 1956.

\bibitem{eckstein1992douglas}
Jonathan Eckstein and Dimitri~P. Bertsekas.
\newblock On the {D}ouglas-{R}achford splitting method and the proximal point
  algorithm for maximal monotone operators.
\newblock {\em Mathematical Programming}, 55(1):293--318, Apr 1992.

\bibitem{gabay1983chapter}
Daniel Gabay.
\newblock Chapter {IX} applications of the method of multipliers to variational
  inequalities.
\newblock In Michel F. and Roland G., editors, {\em Augmented {L}agrangian
  Methods: Applications to the Numerical Solution of Boundary-Value Problems},
  volume~15 of {\em Studies in Mathematics and Its Applications}, pages
  299--331. Elsevier, 1983.

\bibitem{gabay1976dual}
Daniel Gabay and Bertrand Mercier.
\newblock A dual algorithm for the solution of nonlinear variational problems
  via finite element approximation.
\newblock {\em Computers \& Mathematics with Applications}, 2(1):17--40, 1976.

\bibitem{glowinski2013numerical}
Roland Glowinski.
\newblock {\em Numerical Methods for Nonlinear Variational Problems}.
\newblock Scientific Computation. Springer, Berlin Heidelberg, 2013.

\bibitem{glowinski2014alternating}
Roland Glowinski.
\newblock On alternating direction methods of multipliers: A historical
  perspective.
\newblock In W.~Fitzgibbon, Y.~A. Kuznetsov, P.~Neittaanm{\"a}ki, and
  O.~Pironneau, editors, {\em Modeling, Simulation and Optimization for Science
  and Technology}, pages 59--82. Springer Netherlands, Dordrecht, 2014.

\bibitem{glowinski1975approximation}
Roland Glowinski and Americo Marrocco.
\newblock Sur l'approximation, par éléments finis d'ordre un, et la
  résolution, par pénalisation-dualité d'une classe de problèmes de
  dirichlet non linéaires.
\newblock {\em ESAIM: Mathematical Modelling and Numerical Analysis -
  Modélisation Mathématique et Analyse Numérique}, 9(R2):41--76, 1975.

\bibitem{goncalves2017convergence}
Max L.~N. Goncalves, Jefferson~G. Melo, and Renato D.~C. Monteiro.
\newblock Convergence rate bounds for a proximal {ADMM} with over-relaxation
  stepsize parameter for solving nonconvex linearly constrained problems.
\newblock {\em ArXiv e-prints}, February 2017.

\bibitem{guo2017convergence}
Ke~Guo, Deren Han, and Ting-Ting Wu.
\newblock Convergence of alternating direction method for minimizing sum of two
  nonconvex functions with linear constraints.
\newblock {\em International Journal of Computer Mathematics},
  94(8):1653--1669, 2017.

\bibitem{hesse2014alternating}
Robert Hesse, Russel Luke, and Patrick Neumann.
\newblock Alternating projections and {D}ouglas-{R}achford for sparse affine
  feasibility.
\newblock {\em IEEE Transactions on Signal Processing}, 62(18):4868--4881, Sept
  2014.

\bibitem{hesse2013nonconvex}
Robert Hesse and Russell Luke.
\newblock Nonconvex notions of regularity and convergence of fundamental
  algorithms for feasibility problems.
\newblock {\em SIAM Journal on Optimization}, 23(4):2397--2419, 2013.

\bibitem{hiriarturruty2012fundamentals}
Jean-Baptiste Hiriart-Urruty and Claude Lemaréchal.
\newblock {\em Fundamentals of Convex Analysis}.
\newblock Grundlehren Text Editions. Springer Berlin Heidelberg, 2012.

\bibitem{hong2016convergence}
Mingyi Hong, Zhi-Quan Luo, and Meisam Razaviyayn.
\newblock Convergence analysis of alternating direction method of multipliers
  for a family of nonconvex problems.
\newblock {\em SIAM Journal on Optimization}, 26(1):337--364, 2016.

\bibitem{li2017peaceman}
Guoyin Li, Tianxiang Liu, and Ting~Kei Pong.
\newblock Peaceman--{R}achford splitting for a class of nonconvex optimization
  problems.
\newblock {\em Computational Optimization and Applications}, 68(2):407--436,
  Nov 2017.

\bibitem{li2015global}
Guoyin Li and Ting~Kei Pong.
\newblock Global convergence of splitting methods for nonconvex composite
  optimization.
\newblock {\em SIAM Journal on Optimization}, 25(4):2434--2460, 2015.

\bibitem{li2016douglas}
Guoyin Li and Ting~Kei Pong.
\newblock Douglas-{R}achford splitting for nonconvex optimization with
  application to nonconvex feasibility problems.
\newblock {\em Mathematical Programming}, 159(1):371--401, Sep 2016.

\bibitem{liu2017further}
Tianxiang Liu and Ting~Kei Pong.
\newblock Further properties of the forward-backward envelope with applications
  to difference-of-convex programming.
\newblock {\em Computational Optimization and Applications}, 67(3):489--520,
  Jul 2017.

\bibitem{monteiro2018complexity}
Renato D.~C. Monteiro and Chee-Khian Sim.
\newblock Complexity of the relaxed {P}eaceman-{R}achford splitting method for
  the sum of two maximal strongly monotone operators.
\newblock {\em Computational Optimization and Applications}, 70(3):763--790,
  Jul 2018.

\bibitem{nesterov2003introductory}
Yurii Nesterov.
\newblock {\em Introductory lectures on convex optimization: A basic course},
  volume~87.
\newblock Springer, 2003.

\bibitem{patrinos2013proximal}
Panagiotis Patrinos and Alberto Bemporad.
\newblock Proximal {N}ewton methods for convex composite optimization.
\newblock In {\em 52nd IEEE Conference on Decision and Control}, pages
  2358--2363, 2013.

\bibitem{patrinos2014douglas}
Panagiotis Patrinos, Lorenzo Stella, and Alberto Bemporad.
\newblock {D}ouglas-{R}achford splitting: Complexity estimates and accelerated
  variants.
\newblock In {\em 53rd IEEE Conference on Decision and Control}, pages
  4234--4239, Dec 2014.

\bibitem{rockafellar2011variational}
R.~Tyrrell Rockafellar and Roger J.-B. Wets.
\newblock {\em Variational analysis}, volume 317.
\newblock Springer, 2011.

\bibitem{stella2017forward}
Lorenzo Stella, Andreas Themelis, and Panagiotis Patrinos.
\newblock Forward-backward quasi-{N}ewton methods for nonsmooth optimization
  problems.
\newblock {\em Computational Optimization and Applications}, 67(3):443--487,
  Jul 2017.

\bibitem{themelis2018forward}
Andreas Themelis, Lorenzo Stella, and Panagiotis Patrinos.
\newblock Forward-backward envelope for the sum of two nonconvex functions:
  Further properties and nonmonotone linesearch algorithms.
\newblock {\em SIAM Journal on Optimization}, 28(3):2274--2303, 2018.

\bibitem{wang2018global}
Yu~Wang, Wotao Yin, and Jinshan Zeng.
\newblock Global convergence of {ADMM} in nonconvex nonsmooth optimization.
\newblock {\em Journal of Scientific Computing}, Jun 2018.

\bibitem{yan2016self}
Ming Yan and Wotao Yin.
\newblock Self equivalence of the alternating direction method of multipliers.
\newblock In R.~Glowinski, S.~J. Osher, and W.~Yin, editors, {\em Splitting
  Methods in Communication, Imaging, Science, and Engineering}, pages 165--194.
  Springer International Publishing, Cham, 2016.

\end{thebibliography}

\grayout{%
	\section{Temp}
		\input{TeX/Text/Temp.tex}
}

\end{document}